\documentclass[twoside]{article}

\usepackage{fullpage}
\usepackage{jwn-article}
\usepackage{amsmath,amssymb}
\usepackage[all]{xy}            
\usepackage{lgreek}
\usepackage{mathptmx}

\newcommand{\exfmt}{\upshape}   

\newcommand{\goedel}{G\"odel}
\newcommand{\rqq}{\ensuremath{[R(q);q]}}
\newcommand{\sgr}{\ensuremath{17\:\mathrm{Gen}\:r}}
\newcommand{\val}[1]{\ensuremath{\mathsf{\MakeUppercase{#1}}}}
\renewcommand{\ge}{\geqslant}

\newcommand{\set}[1]{\mathcal{#1}}

\newtheorem{definition}{Definition}
\newtheorem{lem}[definition]{Lemma}
\newtheorem{cor}[definition]{Corollary}
\newtheorem{example}[definition]{Example}

\title{Resolving \goedel's Incompleteness Myth: \\ Polynomial Equations
  and Dynamical Systems for Algebraic Logic}

\author{Joseph W. Norman \\ 
  University of Michigan \\ 
  \texttt{jwnorman@umich.edu}}
\date{\today}

\begin{document}

\maketitle

\begin{abstract}
A new computational method that uses polynomial equations and
dynamical systems to evaluate logical propositions is introduced and
applied to G\"odel's incompleteness theorems. The truth value of a
logical formula subject to a set of axioms is computed from the
solution to the corresponding system of polynomial equations. A
reference by a formula to its own provability is shown to be a
recurrence relation, which can be either interpreted as such to
generate a discrete dynamical system, or interpreted in a static way
to create an additional simultaneous equation. In this framework the
truth values of logical formulas and other polynomial objectives have
complex data structures: sets of elementary values, or dynamical
systems that generate sets of infinite sequences of such
solution-value sets. Besides the routine result that a formula has a
definite elementary value, these data structures encode several
exceptions: formulas that are ambiguous, unsatisfiable, unsteady, or
contingent. These exceptions represent several semantically different
types of undecidability; none causes any fundamental problem for
mathematics. It is simple to calculate that G\"odel's formula, which
asserts that it cannot be proven, is exceptional in specific ways:
interpreted statically, the formula defines an inconsistent system of
equations (thus it is called unsatisfiable); interpreted dynamically,
it defines a dynamical system that has a periodic orbit and no fixed
point (thus it is called unsteady). These exceptions are not
catastrophic failures of logic; they are accurate mathematical
descriptions of G\"odel's self-referential construction. G\"odel's
analysis does not reveal any essential incompleteness in formal
reasoning systems, nor any barrier to proving the consistency of such
systems by ordinary mathematical means.
\end{abstract}

\clearpage

\tableofcontents

\clearpage

\section{Introduction}

In the last century, Kurt \goedel's incompleteness theorems
\cite{goedel} sent shockwaves through the world of mathematical logic.
The conventional wisdom is that \goedel's theorems and his
interpretations thereof are correct; the prevalent discussion concerns
what these results mean for logic, mathematics, computer science, and
philosophy \cite{snapper,nagel,goldstein}.  But as I shall demonstrate
here, \goedel's theorems are profoundly misleading and his
interpretations were incorrect: his analysis was corrupted by the
simplistic and flawed notions of truth value and proof that have
troubled logic since antiquity, compounded by his misapplication of a
static definition of consistency to a dynamical system.  Exposing
these errors reveals that reports of logic's demise have been greatly
exaggerated; we may yet realize the rationalist ideals of Leibniz and
complete the logicist and formalist programs of Frege, Russell, and
Hilbert.

There are two key principles here: first, that proof in formal
reasoning systems is an exercise in solving systems of polynomial
equations, yielding solutions that are sets of elementary truth
values; and second, that certain self-referential formula definitions
are recurrence relations that define discrete dynamical systems (and
in turn infinite sequences of basic solution sets).  It is a corollary
to these principles that all of the various syntactic results from
such calculations make semantic sense as the truth values of formulas,
including solution sets that are empty or have multiple members,
and dynamic solutions that change with each iteration and depend on
initial conditions.

These principles are familiar and uncontroversial in the contexts of
elementary algebra and dynamical systems; they apply just as well when
the basic mathematical objects are logical truth values instead of
ordinary numbers.  If you understand how to do arithmetic in different
number systems, what it means to solve equations, and how to deal with
recursive constructions like the Fibonacci sequence, then you can
understand \goedel's mistakes.  You will find the powerful new
paradigm of \emph{dynamic polynomial logic}, which is a continuation
of the pioneering 19th-century work of George Boole
\cite{boole-mal,boole}.  Dynamic polynomial logic is grounded in the
intrinsic unity of logic and mathematics.  

Relative to classical logic, dynamic polynomial logic is
paraconsistent, paracomplete, and modal.  In this algebraic framework
the misguided principle of explosion is corrected: inconsistent axioms
are shown to prove \emph{nothing} instead of \emph{everything}.
Moreover the principle of the excluded middle is clarified; in
classical logic this idea is applied incorrectly, reflecting confusion
between arithmetical and algebraic systems.  Dynamic polynomial logic
computes precise solutions that can be interpreted as alethic and
temporal modalities.

\subsection{\goedel's Argument}

\goedel{} considered formal reasoning systems as described in
Whitehead and Russell's \emph{Principia Mathematica} \cite{pm}, which
was an epic attempt to formalize the whole of mathematics.  \goedel's
basic argument was that every formal reasoning system powerful enough
to describe logical formulas, proof, and natural numbers (like PM)
must allow the construction of a special formula that is semantically
correct but syntactically undecidable: true by metalevel consideration
of its content, but impossible to prove or disprove by mathematical
calculation within the formal system itself.  This special formula,
denoted both \rqq{} and \sgr{} in \goedel's paper, asserts that the
formula itself cannot be proven within the system.  Thus follows the
apparent paradox, which \goedel{} described in this way (as translated
in \cite{heijenoort}):
\begin{quote}
  From the remark that \rqq{} says about itself that it is not
  provable it follows at once that \rqq{} is true, for \rqq{}
  \emph{is} indeed unprovable (being undecidable).  Thus, the
  proposition that is undecidable \emph{in the system PM} still was
  decided by metamathematical considerations.  The precise analysis of
  this curious situation leads to surprising results concerning
  consistency proofs for formal systems, results that will be
  discussed in more detail in Section~4 (Theorem~XI).
\end{quote}
\noindent \goedel's Theorem~VI states that there must exist (in a
formal system like PM) a formula such as his special \rqq{} which can
neither be proven nor disproven within the formal system that contains
it.  His claim that his special formula \rqq{} is semantically true is
presented in the text of his paper but is not called out as a theorem.
\goedel's Theorem~XI states that the existence of such an undecidable
formula renders the consistency of the enclosing formal system itself
an undecidable proposition.

\subsection{The Truth (Value) Is Complicated}

The surprising result from my analysis is that \goedel's special
formula is neither semantically correct nor semantically incorrect;
instead it is \emph{exceptional} in a particular way, relative to the
expectation that a formula should have a definite elementary value.
Such exceptions, which appeared to \goedel{} as `undecidability,' are
features not bugs in formal reasoning systems: it is appropriate that
some logical formulas cannot be proven to be simply true or false,
because they are in fact neither.  Evaluating \goedel's special
formula is analogous to taking the square root of a negative number:
such a square root cannot be proven to be any real number, because it
is not any real number.  More elaborate data structures are needed to
describe the values of some logical formulas, just as complex numbers
are needed to describe the values of some arithmetic formulas.

We shall consider four different exceptions that can be raised by
logical formulas: the static exceptions of \emph{unsatisfiability} and
\emph{ambiguity} and the dynamic exceptions of \emph{unsteadiness} and
\emph{contingency}.  These exceptions represent four distinct types of
`undecidability,' and each type has a particular semantic meaning and
gives characteristic syntactic results.  Ignorance of these exceptions
has plagued logic since ancient times.  One can already find in the
riddles of Aristotle's adversary Eubulides, who originated the liar
paradox in the 4th century \textsc{b.c.}, perfect demonstrations of
unsatisfiable and ambiguous logical propositions \cite{seuren}.  Yet
the same exceptions occur in algebra with ordinary numbers, and in
this context they are well-understood and not at all controversial.
With the appropriate data structures and algorithms, these exceptional
results become no more problematic for logic than $\sqrt{-1}$ is for
algebra.  These four exceptions are consequences of the two key
principles discussed next.

\subsection{Logical Equations Give Sets of Elementary Solutions}
\label{sec:static}

The first key principle is that the axioms in a formal system
(including the definitions of formulas) constitute a system of
simultaneous equations; hence the proper way to describe the truth
value of a formula is to give the \emph{solution-value set} to its
defining system of equations---the possible elementary values of the
formula when all the equations are satisfied.  As in general algebra,
a system of logical equations can have zero, one, or more solutions:
thus the solution-value set can be empty, have one member, or have
many members.

Using binary logic with the elementary values true and false (\val{t}
and \val{f}), there are four possible solution-value sets (hence four
different truth values): the set $\{\val{t}\}$; the set $\{\val{f}\}$;
the set $\{\val{t},\val{f}\}$; and the empty set $\{\}$.  The first
two sets are the definite or `unexceptional' results: if the feasible
set of values for some formula is $\{\val{t}\}$, then that formula is
\emph{necessarily true} (equivalently, it is a \emph{theorem}); if the
feasible set is $\{\val{f}\}$, then the formula is \emph{necessarily
  false} (its negation is a theorem).  The last two sets are
exceptional: if the solution-value set is $\{\val{t},\val{f}\}$, then
the formula is \emph{ambiguous}; if the solution-value set is $\{\}$,
then the formula is \emph{unsatisfiable} (and the underlying axioms
are inconsistent).  This style of categorization extends easily to
sets of elementary values with more than two values, including
infinite sets (like the natural numbers) and uncountable sets (like
the real numbers).  In general we consider truth values in the
\emph{power set} of the set of elementary values.

As Boole explained in his \emph{Laws of Thought} \cite{boole}, logical
formulas with binary truth values can be translated into polynomial
expressions: the coefficients $1$ and $0$ represent the elementary
values true and false; symbolic variables represent basic
propositions, with each variable $x$ subject to the constraint $x^2=x$
to ensure that its only feasible values are $0$ and $1$; logical
conjunction translates as multiplication; logical negation translates
as the difference from $1$; and logical disjunction translates as a
certain combination of addition and multiplication.  Boole presented a
complete algorithm to translate any logical formula into a polynomial
with ordinary integer coefficients.

For example, using Boole's original method the logical formula $x
\rightarrow y$ translates as the polynomial $xy-x+1$ and the logical
formula $y \rightarrow x$ translates as the polynomial $xy-y+1$.
Using these translations, constraining each formula to be true
(i.e.\ to equal the polynomial $1$), and separately constraining each
variable to be either $0$ or $1$ yields the following system of
polynomial equations:
\begin{equation}
  \label{eq:xy-impl}
  \begin{array}{rcl}
    xy - x + 1 & = & 1 \\
    xy - y + 1 & = & 1 \\
    x^2 & = & x \\
    y^2 & = & y
  \end{array}
\end{equation}
The only solutions to these equations for the variables $(x,y)$ are
$(0,0)$ and $(1,1)$.  You can see that $x$ and $y$ have the same value
in each solution, so you might expect the biconditional formula $x
\leftrightarrow y$ to be a theorem given the axioms $x \rightarrow y$
and $y \rightarrow x$.  Using Boole's original method, the logical
formula $x \leftrightarrow y$ is translated into the polynomial
$2xy-x-y+1$; substituting either solution for $(x,y)$ the value of
this polynomial is $1$.  Thus the solution-value set for the
polynomial $2xy-x-y+1$, subject to the constraints in
Equation~\ref{eq:xy-impl}, is the set $\{1\}$.

Translating back into logical notation, the solution-value set for the
formula $x \leftrightarrow y$ given the axioms $x \rightarrow y$ and
$y \rightarrow x$ is $\{\val{t}\}$; therefore $x \leftrightarrow y$ is
called a \emph{theorem} relative to these axioms.  However, neither
the formula $x$ nor the formula $y$ has a definite value given these
equations; each formula has the solution-value set $\{0,1\}$ in
polynomial notation, which is $\{\val{f},\val{t}\}$ in logical
notation.  Therefore each formula $x$ and $y$ is called
\emph{ambiguous} given the axioms $x \rightarrow y$ and $y \rightarrow
x$.  Note that neither $x$ nor $\neg x$ is a theorem given these
axioms; likewise neither $y$ nor $\neg y$ is a theorem.  This
`undecidability' for the formulas $x$ and $y$ is a correct
interpretation of the relevant equations, not a sign of pathological
incompleteness in formal reasoning.

Infeasible equations have empty solution sets and render unsatisfiable
all formulas that are subject to them.  For example the simultaneous
equations:
\begin{equation}
  z = 1-z, \quad z^2 = z
\end{equation}
have no solution; the first constraint is violated for both values
$z=0$ and $z=1$ that satisfy the second constraint.  Therefore the
formula $z$ is \emph{unsatisfiable} subject to these equations, as are
the formulas $0$, $1$, and $1-z$; all share the empty solution-value
set $\{\}$.  The equation $z=1-z$ is a polynomial translation of the
logical axiom $z = \neg z$ which states that the formula $z$ is
defined to be true exactly if it is not true.  This is one way to
model the liar paradox, and in this version the problem is no more
paradoxical than the unsatisfiable equation $1=0$.

\subsection{Self-Reference Gives Dynamical Systems}
\label{sec:intro-dynamic}

The second key principle is that, if it is permitted for a formula
definition (or any other axiom or equation) to refer to the solution
generated by evaluating its own system of equations, then the value of
every formula may gain a more complex data structure.  Such a solution
self-reference is a \emph{recurrence relation} that can be interpreted
in two distinct ways: in the \emph{dynamic} interpretation the
recurrence is used as such to define infinite sequences of solutions
governed by a \emph{discrete dynamical system}; and in the
\emph{static} interpretation the recurrence is used to generate an
additional simultaneous equation.  Both interpretations can be
reasonable and useful, but it is important not to confuse them during
analysis.

To illustrate, let us consider the quadratic equation $2x^2+3x+c=0$ in
which we define the coefficient $c$ to be the number of real solutions
(for $x$) to the equation in which it appears.  In other words we have
the following specification for a system of equations, using
real-valued variables $x,c\in\mathbb{R}$:
\begin{eqnarray}
  2x^2+3x+c & = & 0 \label{eq:cc-orig} \\
  c & := & \mbox{the number of real solutions for $x$ to
    Equation~\ref{eq:cc-orig}} 
  \label{eq:cc-orig-ref}
\end{eqnarray}
In this specification the value of $c$ depends on itself.  To model
this dependence, let us introduce an \emph{evolution function} denoted
$F(c)$ that gives the number of real solutions to
Equation~\ref{eq:cc-orig} when $c$ takes the value supplied as the
function's argument.  In the dynamic interpretation we take the
specification to define a system that changes over time according to
this evolution function $F(c)$, with the state $c_{t+1}$ at the next
time $t+1$ given by the value $F(c_t)$ of the evolution function
applied to the current state $c_t$.  Thus we have:
\begin{equation}
  2x_t^2+3x_t+c_t = 0, \quad c_{t+1} \Leftarrow F(c_t) \label{eq:ct};
  \quad t \in \{0,1,2,\ldots\}
\end{equation}
Alternatively, in the static interpretation we take the problem
specification to mean that the input value $c$ and the output value
$F(c)$ must agree simultaneously, as in:
\begin{equation}
  2x^2+3x+c = 0, \quad c = F(c) \label{eq:cc}
\end{equation}

We must now determine what this evolution function $F(c)$ is and which
particular values of $c$ need to be considered in solving the
equations.  For this problem the quadratic formula serves both needs;
a more general approach will be presented later.  As you may recall,
the quadratic formula states that the number of distinct real roots of
an equation $ax^2+bx+c=0$ depends on its determinant $b^2-4ac$: there
are no real solutions if the determinant is negative, one if it is
zero, and two if it is positive.  Therefore, the possible values of
$c$ are $\{0,1,2\}$; and using the quadratic-equation coefficients
$a=2$ and $b=3$ from Equation~\ref{eq:cc-orig}, the evolution function
$F(c)$ must satisfy:
\begin{eqnarray}
  F(c) & = & \left\{
  \begin{array}{rl}
    0, & 9-8c < 0 \\
    1, & 9-8c = 0 \\
    2, & 9-8c > 0
  \end{array}
  \right.
\end{eqnarray}
Therefore $F(0)=2$ since $9-8\cdot0>0$; $F(1)=2$ since $9-8\cdot1>0$;
and $F(2)=0$ since $9-8\cdot2<0$.  Using polynomial interpolation, a
closed-form function $F(c)$ can be constructed that performs exactly
these mappings $0 \mapsto 2$, $1 \mapsto 2$, and $2 \mapsto 0$:
\begin{eqnarray}
  F(c) & : & -c^2+c+2 \label{eq:fc}
\end{eqnarray}

For the static interpretation, substituting the evolution function
from Equation~\ref{eq:fc} into the system in Equation~\ref{eq:cc} and
then making explicit the domain of each variable produces the
following system of equations:
\begin{equation}
  2x^2+3x+c=0, \quad c = -c^2+c+2;
  \quad x \in \mathbb{R}, \quad c \in \{0,1,2\} \subset \mathbb{R}
  \label{eq:cc-static}
\end{equation}
After rearranging the second equation it is evident that this system
has no solution: $c^2=2$ implies $c=\sqrt{2}$ but $c$ is required to
be $0$, $1$, or $2$.  Thus in the static interpretation, the formula
$x$ is \emph{unsatisfiable} and the whole system of equations is
inconsistent.

For the dynamic interpretation, substituting the evolution function
from Equation~\ref{eq:fc} into the system in Equation~\ref{eq:ct} and
leads to the following system of equations, which includes the
recurrence $c_{t+1} \Leftarrow F(c_t)$:
\begin{equation}
  2x_t^2+3x_t+c_t=0, \quad c_{t+1} \Leftarrow -c_t^2+c_t+2;
  \quad x_t \in \mathbb{R}, \quad c_t \in \{0,1,2\} \subset
  \mathbb{R}, \quad t \in \{0,1,2,\ldots\}
  \label{eq:cc-dynamic}
\end{equation}
In this interpretation we consider $c_t$ to be the state of the
dynamical system at time $t$.  Each initial state $c_0\in\{0,1,2\}$
generates an infinite sequence $(c_0,c_1,\ldots)$ of states at
successive times $t$:
\begin{equation}
  \left[ \begin{array}{rcl}
      c_0 = 0 & \mapsto & (0,2,0,2,\ldots) \\
      c_0 = 1 & \mapsto & (1,2,0,2,\ldots) \\
      c_0 = 2 & \mapsto & (2,0,2,0,\ldots) \\
    \end{array} \right]
  \label{eq:cc-seq-c}
\end{equation}
In the state $c_t=0$ the main equation $2x_t^2+3x_t+c_t=0$ specializes
to $2x_t^2+3x_t+0=0$ which has the two solutions $x_t=-\frac{3}{2}$
and $x_t=0$.  Then given $c_t=1$ the main equation becomes
$2x_t^2+3x_t+1=0$ which has two solutions $x_t=-1$ and
$x_t=-\frac{1}{2}$.  Finally, given $c_t=2$ the main equation becomes
$2x_t^2+3x_t+2=0$ which has no real solutions.  Therefore
Equation~\ref{eq:cc-dynamic} also generates a collection of sequences
of solution-value sets for $x$, again depending on the initial
condition $c_0$:
\begin{equation}
  \left[ \begin{array}{rcl}
      c_0 = 0 & \mapsto & 
      (\{-\frac{3}{2},0\},\{\},\{-\frac{3}{2},0\},\{\},\ldots) \\
      c_0 = 1 & \mapsto & 
      (\{-1,-\frac{1}{2}\},\{\},\{-\frac{3}{2},0\},\{\},\ldots) \\
      c_0 = 2 & \mapsto & 
      (\{\},\{-\frac{3}{2},0\},\{\},\{-\frac{3}{2},0\},\ldots)
    \end{array} \right]
  \label{eq:cc-seq-x}
\end{equation}
These collections of infinite sequences are governed by a discrete
dynamical system.  This system can be displayed compactly as a graph
in which each node indicates a state of the parameter $c$ and each
edge shows the solution-value set for the formula $x$ that is
generated from assuming the state corresponding to the originating
node:
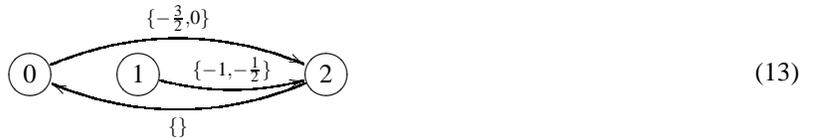
\begin{equation}
  \xymatrix{
    *++[o][F]{0} \ar@/^3ex/[rrr]^{\{-\frac{3}{2},0\}} & 
    *++[o][F]{1} \ar@/_/[rr]^{\{-1,-\frac{1}{2}\}} & &
    *++[o][F]{2} \ar@/^3ex/[lll]^{\{\}}
  }
  \label{eq:cc-graph}
\end{equation}
This graph, constructed from Equation~\ref{eq:cc-dynamic}, offers some
explanation as to \emph{why} the static system in
Equation~\ref{eq:cc-static} is unsatisfiable: the matching dynamical
system has no fixed points, only a periodic orbit that never reaches a
steady state.  You can read off the sequences of states or
solution-value sets generated by this dynamical system by following
the edges in its graph.

In general many different patterns could be discerned in the behaviors
of dynamical systems; but for the purpose of classifying formulas we
consider only the number of fixed points.  If a dynamical system has
exactly one fixed point then that system is called \emph{stable}; if
it has no fixed points then it is called \emph{unstable}; and if it
has more than one fixed point then it is called \emph{contingent}.  In
these terms, the dynamical system defined by
Equation~\ref{eq:cc-dynamic} and graphed in Equation~\ref{eq:cc-graph}
is unstable.  Unstable dynamical systems correspond to unsatisfiable
static equations, and contingent dynamical systems correspond to
ambiguous static equations.

\subsection{Anticipating \goedel's Error}

We shall see that \goedel's special formula \rqq{} behaves just like
the self-referential quadratic equation above.  Using $x$ to represent
\goedel's formula, it turns out that his definition specifies the
evolution function $F(x):1-x$.  When \goedel's reference to
provability is interpreted statically as the constraint $x=F(x)$, his
self-denying formula \rqq{} becomes an unsatisfiable system of
equations:
\begin{equation}
  x=1-x, \quad x \in \{0,1\}
\end{equation}
And when interpreted dynamically as $x_{t+1} \Leftarrow F(x_t)$,
\goedel's formula becomes the recurrence:
\begin{equation}
  x_{t+1} \Leftarrow 1-x_t, \quad x \in \{0,1\}
\end{equation}
This specifies a simple dynamical system with one periodic orbit and
no fixed points.  This dynamical system generates an alternating
sequence of values for its state $x_t$ at successive times $t$, for
each initial condition $x_0$:
\begin{equation}
  \left[ \begin{array}{rcl}
      x_0 = 0 & \mapsto & (0,1,0,1,\ldots) \\
      x_0 = 1 & \mapsto & (1,0,1,0,\ldots)
    \end{array} \right]
\end{equation}
By generating these results, dynamic polynomial analysis will show
that \goedel's `formula \rqq{} that is true if and only if is not
provable' is exceptional in precisely the same way as is `the
quadratic equation $2x^2+3x+c=0$ that has exactly $c$ solutions':
interpreted as static systems of equations, neither specification can
be satisfied; and interpreted as dynamical systems, both
specifications lead to sequences that oscillate infinitely and never
converge to a fixed value.

These types of `undecidability' are not syntactic aberrations; they
are semantically appropriate descriptions of the mathematical objects
specified by their respective self-referential definitions.  In fact,
these complex results only seem exceptional because of the misguided
expectation that the specified mathematical objects should have simple
elementary values.  Imagine the confusion that would result if we were
to speak of `the Fibonacci \emph{formula}' (expecting it to have a
definite numeric value) or `\emph{the} Fibonacci number' (expecting
there to be just one) instead of `the Fibonacci \emph{sequence}' and
`\emph{a} Fibonacci number.'  Although \goedel's self-denying
\emph{formula} may seem puzzling, \goedel's oscillating
\emph{sequence} and \goedel's inconsistent \emph{equation} are rather
less so.

It remains to be demonstrated that systems of logical axioms are
indeed equivalent to systems of polynomial equations, and that
self-reference of the type that \goedel{} described is accurately
translated using evolution functions, recurrence relations, and
dynamical systems.  But before discussing the details of translating
logic to algebra, let us review several useful kinds of results that
can be calculated by general methods of algebra from systems of
polynomial equations and from discrete dynamical systems.

\section{General Analysis of Polynomial Equations}
\label{sec:polynomial}

Polynomials are an important class of formulas in elementary and
abstract algebra.  In this section we introduce terminology and
notation to describe systems of polynomial equations and their
solutions (for the task of \emph{formula evaluation}), and we consider
ways to construct and to count the members of a polynomial ring that
meet certain desirable criteria (for the task of \emph{formula
  discovery}).  We develop simple algorithms to perform these tasks by
hand for small problems, and make reference to more general and
efficient methods from computational algebraic geometry.  For the
moment we are just discussing polynomials in the context of general
algebra, without any mention of their provenance in logic and
axiomatic formal reasoning.

\subsection{Formula Evaluation}
\label{sec:evaluation}

We consider systems of polynomial equations whose variables and
coefficients take values from an algebraic field $K$ such as the real
number system $\mathbb{R}$, the rational number system $\mathbb{Q}$,
or a finite field $\mathbb{F}_d$ of order $d$ (with $d$ a prime
number).  Note that in a finite field $\mathbb{F}_d$ we must use
integer arithmetic modulo $d$ for calculation; for example in the
binary finite field $\mathbb{F}_2$ the sum $1+1 \Rightarrow 0$ since
$2 \equiv 0 \pmod{2}$.  \cite{ideals-varieties} provides a good
general reference for polynomials and algebraic geometry; \cite{clark}
and \cite{eves} describe fields and other structures in abstract
algebra; modular arithmetic is discussed in \cite{concrete}.

\begin{definition}[System of Polynomial Equations]
\label{def:polynomial-system}
Given an algebraic field $K$ and a vector $\mathbf{x} :=
(x_1,x_2,\ldots,x_n)$ of variables, the polynomial ring over the
variables $\mathbf{x}$ with coefficients in $K$ is denoted
$K[\mathbf{x}]$.  We consider a set $A$ of simultaneous equations:
\begin{eqnarray*}
A & := & \left\{ \; q_1=0,\; q_2=0,\; \ldots,\; q_m=0 \; \right\}
\end{eqnarray*}
in which each polynomial $q_j$ is a member of the ring $K[\mathbf{x}]$
and also constrained to equal zero.  The \emph{solution set} to this
system of equations, denoted $\set{V}(A)$ or
$\set{V}(q_1=0,\;q_2=0,\;\ldots,\;q_m=0)$, is the set of values of
$\mathbf{x}$ for which these equality constraints are satisfied,
assuming that each variable $x_i$ takes a value in the field $K$:
\begin{eqnarray*}
  \set{V}(A) & := & 
  \left\{ \; \mathbf{x} \in K^n \; : \;
  q_1(\mathbf{x}) = 0, \;
  q_2(\mathbf{x}) = 0, \; \ldots, \;
  q_m(\mathbf{x}) = 0 \;
  \right\}
\end{eqnarray*}
By its construction the solution set $\set{V}(A) \subseteq K^n$
must be a subset of the affine space $K^n$.
\end{definition}
In the extreme case that the equations in $A$ are inconsistent then
$\set{V}(A)$ is the empty set; conversely if the equations are
tautological (or if there are no equations) then $\set{V}(A)$ is the
entire set $K^n$.  In algebraic geometry, the solution set
$\set{V}(A)$ is called the \emph{affine variety} defined by the
polynomials $\{q_1,q_2,\ldots,q_m\}$ used in the equations.  Note that
in the special case $n=0$ that there are no variables, we imagine that
the affine space $K^0$ has one member $()$ that is the unique
zero-length tuple.  Thus with $n=0$ the solution set $\set{V}(0=0)
\Rightarrow \{()\}$ has one member since the equation $0=0$ holds;
however $\set{V}(1=0) \Rightarrow \{\}$ is the empty set since the
equation $1=0$ does not hold.  Additionally, we take the polynomial
ring $K[]$ with no variables to be the same as the original field $K$.

We next consider the evaluation of a polynomial function subject to
the constraints in a system of polynomial equations.  Recall that the
power set $2^s$ of any set $s$ is the set of all possible subsets of
$s$ (including by necessity the empty set and the original set $s$
itself).

\begin{definition}[Solution-Value Set]
\label{def:objective}
Consider a system of polynomial equations as described in
Definition~\ref{def:polynomial-system}, along with an \emph{objective
  formula} $p \in K[\mathbf{x}]$ from the same polynomial ring.  The
\emph{solution-value set} $\set{S}_A(p)$ of the objective $p$ subject to
the equations $A$ is defined as the set of its feasible values when
the equations in the system are satisfied:
\begin{eqnarray*}
  \set{S}_A(p) & := &
  \left\{ \;
  p(\mathbf{x}) \in K \; : \; \mathbf{x} \in \set{V}(A) \;
  \right\}
\end{eqnarray*}
in other words the image under $p$ of the solution set
$\set{V}(A)$.  By its construction the solution-value set
$\set{S}_A(p) \subseteq K$ must be a subset of the set of elementary
values in the original field $K$.  Equivalently the solution-value set
must be a member of the power set of $K$:
\begin{eqnarray*}
  \set{S}_A(p) & \in & 2^K
\end{eqnarray*}
Here using $K$ to stand for the set of elementary values in the field
as well as the field itself.  

If the solution-value set $\set{S}_A(p)$ for a polynomial $p$ is a
singleton $\{k\}$ that contains just one member $k \in K$, then we
describe $p$ as \emph{necessarily $k$} subject to $A$; this is the
unexceptional result.  If the solution-value set is empty, then we
describe $p$ as \emph{unsatisfiable} subject to $A$; this means that
the equations $A$ in the system are inconsistent.  Otherwise the
solution-value set must have more than one member (including the
special cases that it is infinite or uncountable) and we describe $p$
as \emph{ambiguous} subject to $A$.
\end{definition}
At the extremes the solution-value set $\set{S}_A(p)$ may be the empty
set or the entire set of elementary values from the underlying
algebraic structure $K$.

Moving on, when using coefficients from the binary finite field
$\mathbb{F}_2$ it is possible to simplify systems of polynomial
equations in the following way.

\begin{lem}[Conjunction of Binary Constraints]
\label{lem:conjunction}
Any system of polynomial equations with coefficients in the binary
finite field $\mathbb{F}_2$ or with each polynomial constrained to the
values $0$ and $1$ can be simplified to a single equation.  For any
$a,b\in\{0,1\}$ the product $ab=1$ if and only if both factors $a=1$
and $b=1$.  From this it follows that any set $A$ of equations:
\begin{displaymath}
  \left\{ \;
  q_1(\mathbf{x}) = 0, \;
  q_2(\mathbf{x}) = 0, \; \ldots, \;
  q_m(\mathbf{x}) = 0
  \; \right\}
\end{displaymath}
with each $q_j\in\mathbb{F}_2[\mathbf{x}]$ or each
$q_j(\mathbf{x})\in\{0,1\}$ for all feasible values of $\mathbf{x}$,
must have the same solution set as the single equation:
\begin{eqnarray*}
  (q_1(\mathbf{x})+1) (q_2(\mathbf{x})+1) \cdots (q_m(\mathbf{x})+1)
  & = & 1
\end{eqnarray*}
In this case the set $A$ of equations given in
Definition~\ref{def:polynomial-system} can be replaced with the
singleton $\{q^*=0\}$ or the equivalent $\{q^*+1=1\}$ whose solitary
member uses the \emph{conjunction polynomial} $q^*$ derived from the
original equations:
\begin{eqnarray*}
  q^* & := &
  (q_1(\mathbf{x})+1) (q_2(\mathbf{x})+1) \cdots (q_m(\mathbf{x})+1)
  -1
\end{eqnarray*}
The conjunction polynomial produces the same solution set as the
original equations:
\begin{displaymath}
  \set{V}(q^*=0) \quad = \quad 
  \set{V}(q^*+1=1) \quad = \quad 
  \set{V}(q_1=0,\;q_2=0,\;\ldots,\;q_m=0)
\end{displaymath}
and consequently identical solution-value sets for any objective
formula.
\end{lem}
Note that in the special case $m=0$ that there are no constraints we
take the product of zero factors to be the multiplicative identity.
Hence from $A=\{\}$ we derive the conjunction polynomial $q^*=1-1$
which gives the tautological constraint $0=0$.


\begin{example}\exfmt
\label{ex:xy}

To illustrate these definitions, consider the following system of two
simultaneous equations in the variables $x$ and $y$ (thus $n=2$,
$x_1=x$, $x_2=y$, $\mathbf{x}=(x,y)$, and $m=2$):
\begin{eqnarray}
x (y+1) & = & 0 \label{eq:eqx} \\
y (x+1) & = & 0 \label{eq:eqy}
\end{eqnarray}
Expanding these polynomials yields the constraint set:
\begin{eqnarray}
A & := & \left\{ xy+x=0,\; xy+y=0 \right\} \label{eq:eqa1}
\end{eqnarray}
Let us proceed to evaluate this system of polynomial equations using
two different algebraic structures: the real numbers $\mathbb{R}$ and
the binary finite field $\mathbb{F}_2$.

First we assume that the variables $x$ and $y$ and their coefficients
take real-number values.  Let us use $\set{V}(A)_\mathbb{R}$ to
denote the solution set to Equations \ref{eq:eqx} and \ref{eq:eqy}
using the algebraic structure $\mathbb{R}$.  According to
Definition~\ref{def:polynomial-system}:
\begin{eqnarray}
  \set{V}(A)_\mathbb{R} & := & \left\{ (x,y) \in \mathbb{R}^2: 
  \begin{array}{l} 
    xy+x = 0 \\ 
    xy+y = 0 
  \end{array} \right\}
\end{eqnarray}
You can see from inspection that only two pairs $(x,y)$ of real
numbers satisfy the equations, namely $(-1,-1)$ and $(0,0)$.  Thus we
have the solution set:
\begin{eqnarray}
\set{V}(A)_\mathbb{R} & \Rightarrow & \{ (-1,-1), (0,0) \}
\end{eqnarray}
This variety $\set{V}(A)_\mathbb{R}$ is a subset of the affine space
$\mathbb{R}^2$ of all pairs of real numbers.  Given this solution set,
it is straightforward to identify the solution-value sets for the
objective formulas $x$, $y$, and $x-y+1$ by substituting the values
from $\set{V}(A)_\mathbb{R}$ into each objective:
\begin{eqnarray}
\set{S}_A(x)_\mathbb{R} & \Rightarrow & \{ -1, 0 \} \\
\set{S}_A(y)_\mathbb{R} & \Rightarrow & \{ -1, 0 \} \\
\set{S}_A(x-y+1)_\mathbb{R} & \Rightarrow & \{ 1 \}
\end{eqnarray}

Next if we use the binary finite field $\mathbb{F}_2$ as the algebraic
structure for the variables $x$ and $y$ and their coefficients, then
the polynomial system defined by Equations \ref{eq:eqx} and
\ref{eq:eqy} has the solution set:
\begin{equation}
  \begin{array}{rcl}
    \set{V}(A)_{\mathbb{F}_2} & := & 
    \left\{ (x,y) \in (\mathbb{F}_2)^2: 
    \begin{array}{l} 
      xy+x = 0 \\ 
      xy+y = 0 
    \end{array} \right\} \\
    & \Rightarrow & \{ (0,0), (1,1) \} \label{eq:xysol}
\end{array}
\end{equation}
Recall that in the finite field $\mathbb{F}_2$ the sum $1+1
\Rightarrow 0$ using modular arithmetic.  Based on this solution set
$\set{V}(A)_{\mathbb{F}_2}$, the solution-value sets for the
objectives $x$, $y$, and $x-y+1$ are:
\begin{eqnarray}
\set{S}_A(x)_{\mathbb{F}_2} & \Rightarrow & \{ 0, 1 \} \\
\set{S}_A(y)_{\mathbb{F}_2} & \Rightarrow & \{ 0, 1 \} \\
\set{S}_A(x-y+1)_{\mathbb{F}_2} & \Rightarrow & \{ 1 \} \label{eq:xy1}
\end{eqnarray}
In either polynomial ring $\mathbb{R}[x,y]$ or $\mathbb{F}_2[x,y]$ the
formula $x-y+1$ has only one feasible value given the constraints in
Equations \ref{eq:eqx} and \ref{eq:eqy}: it is necessarily $1$.
However, given these same constraints each formula $x$ and $y$ is
ambiguous in either polynomial ring: the solution-value sets
$\{-1,0\}$ and $\{0,1\}$ each have two members.

\end{example}

\subsection{Formula Discovery}
\label{sec:discovery}

We now consider the inverse problem: subject to a set of polynomial
equations, instead of computing the solution-value set for some
polynomial objective formula, we must find the set of polynomials that
yield a given solution-value set.

\begin{definition}[Inverse-Value Set]
\label{def:inverse}
Consider a system of polynomial equations as specified in
Definition~\ref{def:polynomial-system}, along with a query set $s
\subseteq K$ of elementary values from the algebraic field $K$.  The
\emph{inverse-value set} $\set{S}_A^{-1}(s)$ is the set of
polynomials in the ring $K[\mathbf{x}]$ whose solution-value set is
exactly $s$:
\begin{eqnarray*}
  \set{S}_A^{-1}(s) & := & \left\{ \;
  p \in K[\mathbf{x}] \; : \; \set{S}_A(p) = s
  \; \right\}
\end{eqnarray*}
When working with coefficients in the finite field $\mathbb{F}_2$,
replacing the original set $A$ of equations with the set $\{q^*=0\}$
using the conjunction polynomial from Lemma~\ref{lem:conjunction} does
not change the inverse-value sets.  In that case, for every possible
query $s \subseteq K$:
\begin{displaymath}
  \set{S}_{\{q^*=0\}}^{-1}(s) \quad = \quad
  \set{S}_{\{q_1=0,\;q_2=0,\;\ldots,\;q_m=0\}}^{-1}(s)
\end{displaymath}
If the set $A$ of equations is inconsistent then according to
Definition~\ref{def:objective} the solution-value set $\set{S}_A(p)$
for every polynomial $p \in K[\mathbf{x}]$ must be the empty set.
Hence when $A$ is inconsistent the inverse-value set
$\set{S}_A^{-1}(\{\})$ for the empty-set query yields the entire
polynomial ring $K[\mathbf{x}]$, and the inverse-value set
$\set{S}_A^{-1}(s)$ with any non-empty query $s$ yields the empty set
of polynomials.  In other words, given infeasible constraints, every
polynomial has no feasible value and no polynomial has any feasible
value.
\end{definition}
Considering a fixed set of equations $A$, the inverse-value-set
function $\set{S}_A^{-1}$ is the \emph{inverse image} or
\emph{preimage} of the solution-value-set function $\set{S}_A$.  Thus
if $\set{S}_A^{-1}(s) \Rightarrow \{p_1,p_2,\ldots\}$ then for every
polynomial $p_i$ in this inverse-value set it holds that
$\set{S}_A(p_i)=s$; and furthermore the members $\{p_i\}$ of this
inverse-value set are the \emph{only} polynomials in the ring
$K[\mathbf{x}]$ that have the solution-value set $s$.  In algebraic
geometry, the inverse-value set $\set{S}_A^{-1}(\{0\})$ with the query
$\{0\}$ is called the \emph{ideal} of the affine variety $\set{V}(A)$:
the set of every polynomial that vanishes on $\set{V}(A)$.  The
characteristics of ideals \cite{ideals-varieties} lead to the
following special case.

\begin{lem}[The Ideal from a Single Equation]
  \label{lem:singleton}
  In the special case that the set $A$ of equations contains a single
  equation $q=0$ that is satisfiable, then every member of the
  inverse-value set $\set{S}_{\{q=0\}}^{-1}(\{0\})$ must be the
  product of some polynomial $p$ from $K[\mathbf{x}]$ with the
  polynomial $q$ from the equation; conversely every such product $p
  \times q$ must be a member of the inverse-value set:
  \begin{eqnarray*}
    \set{S}_{\{q=0\}}^{-1}(\{0\}) & = & 
    \left\{ p \times q : p \in K[\mathbf{x}] \right\}
  \end{eqnarray*}
  The products $p \times q$ are not necessarily distinct for different
  values of $p$.  Recall from Definition~\ref{def:inverse} that
  infeasible equations yield empty inverse-value sets for non-empty
  queries: thus $\set{S}_{\{q=0\}}^{-1}(\{0\}) \Rightarrow \{\}$ if
  the equation $q=0$ is infeasible.
\end{lem}

Since there are dedicated methods in algebraic geometry to compute the
ideal generated by a set of polynomials (in my terminology the
inverse-value set for the query $\{0\}$), it is useful to define
inverse-value sets for queries other than $\{0\}$ in terms of the
ideal.
\begin{cor}[Incrementing Inverse Polynomials]
  \label{cor:increment}
  For any query $\{k\}$ that contains a solitary value $k \in K$, the
  inverse-value set $\set{S}_A^{-1}(\{k\})$ can be calculated by
  adding that value $k$ to every polynomial in the inverse-value set
  for the query $\{0\}$:
  \begin{eqnarray*}
    \set{S}_A^{-1}(\{k\}) & = & 
    \left\{ p+k : p \in \set{S}_A^{-1}(\{0\}) \right\}
  \end{eqnarray*}
  In the case that the inverse-value set $\set{S}_A^{-1}(\{0\})$ was
  generated from a single constraint as in Lemma~\ref{lem:singleton},
  then the same incremented inverse-value set is also given by:
  \begin{eqnarray*}
    \set{S}_{\{q=0\}}^{-1}(\{k\}) & = & 
    \left\{ p \times q + k : p \in K[\mathbf{x}] \right\}
  \end{eqnarray*}
  This corollary allows closed-form description of inverse-value sets
  for singleton queries.
\end{cor}

Continuing on, it is useful to appreciate that the set of polynomials
over a finite number of variables with coefficients in a finite field
is itself finite.

\begin{lem}[Counting Polynomials in Finite Fields]
\label{lem:count}
The polynomial ring $\mathbb{F}_d[\mathbf{x}]$ over $n$ variables
$\mathbf{x}:=(x_1,x_2,\ldots,x_n)$ with coefficients in a finite field
$\mathbb{F}_d$ of order $d$ (with $d$ a prime number) contains
exactly \[d^{d^n}\] distinct polynomials, each of which can be
expressed as a unique sum $\sum_j c_j t_j$ of coefficients and
monomials.  In this sum each coefficient $c_j$ is a member of the
field $\mathbb{F}_d$ and each monomial (power product) $t_j$ is a
product of the variables $x_1$ through $x_n$ in which each variable is
raised to a power between $0$ and $d-1$.  There are $d^n$ possible
monomials.
\end{lem}
Since (by Fermat's Little Theorem) the identity $a^d=a$ holds for
every element $a\in\mathbb{F}_d$ of a finite field of order $d$ (with
$d$ a prime number), exponents higher than $d-1$ are not required in
monomials that use $\mathbb{F}_d$.  For example, using the finite
field of order $2$ and two variables $x$ and $y$, the resulting
polynomial ring $\mathbb{F}_2[x,y]$ has a set of $2^2=4$ possible
monomials:
\begin{eqnarray}
\left\{ \; x^1 y^1,\; x^1 y^0,\; x^0 y^1,\; x^0 y^0 \; \right\}
& \Rightarrow & \left\{ xy, x, y, 1 \right\}
\end{eqnarray}
In any polynomial in this ring, each of these $4$ monomials must be
assigned one of the $2$ coefficients in $\{0,1\}$.  Thus there are
$2^4=16$ possible polynomials in $\mathbb{F}_2[x,y]$; they are listed
in Table~\ref{tbl:worksheet}.  Alternatively, using coefficients from
the finite field of order $3$ there would be $3^{3^2} = 19,683$
possible polynomials in the corresponding ring $\mathbb{F}_3[x,y]$.

\begin{example}\exfmt

We can enumerate all sixteen polynomials in the ring
$\mathbb{F}_2[x,y]$ with binary finite-field coefficients, and
partition them into the four possible inverse-value sets (for the
queries $s=\{\}$, $s=\{0\}$, $s=\{1\}$, and $s=\{0,1\}$)
given the constraints $A=\{xy+x=0,\;xy+y=0\}$ from Equations
\ref{eq:eqx} and \ref{eq:eqy}.  Table~\ref{tbl:worksheet} shows these
sixteen polynomials and their solution-value sets;
Definition~\ref{def:inference} below explains how the table was
created.  From the table you can see that there are four polynomials
whose only feasible value is zero:
\begin{equation}
\label{eq:theta-0}
\begin{array}{rcl}
\set{S}_A^{-1}(\{0\}) & \Rightarrow & 
\left\{ \; p_{1},\; p_{7},\; p_{11},\; p_{13} \; \right\} \\
& \Rightarrow & \left\{ \; 0,\; x+y,\; xy+y,\; xy+x \; \right\}
\end{array}
\end{equation}
There are also four polynomials whose only feasible value is one:
\begin{equation}
\label{eq:theta-1}
\begin{array}{rcl}
\set{S}_A^{-1}(\{1\}) & \Rightarrow & 
\left\{ \; p_{2},\; p_{8},\; p_{12},\; p_{14} \; \right\} \\
& \Rightarrow & \left\{ \; 1,\; x+y+1,\; xy+y+1,\; xy+x+1 \; \right\}
\end{array}
\end{equation}
There are no polynomials whose solution-value
set is empty (because the constraints are consistent):
\begin{eqnarray}
\set{S}_A^{-1}(\{\}) & \Rightarrow & \{\}
\end{eqnarray}
The remaining eight polynomials share the solution-value set $\{0,1\}$
of both binary finite-field values:
\begin{eqnarray}
\set{S}_A^{-1}(\{0,1\}) & \Rightarrow &
\left\{ \;
p_{3},\; p_{4},\; p_{5},\; p_{6},\; p_{9},\; p_{10},\; p_{15},\; p_{16} \;
\right\}
\end{eqnarray}
Thus, using the constraints in Equations \ref{eq:eqx} and
\ref{eq:eqy}, every polynomial $p_i \in \mathbb{F}_2[x,y]$ is assigned
to one of the four possible inverse-value sets for polynomials with
coefficients in the binary finite field.

Some of these inverse-value sets can be derived in another way.  Since
we are using polynomials with coefficients in the binary finite field
$\mathbb{F}_2$, Lemma~\ref{lem:conjunction} states that the original
constraints $xy+x=0$ and $xy+y=0$ in Equations \ref{eq:eqx} and
\ref{eq:eqy} are equivalent to the single constraint:
\begin{eqnarray}
(xy+x+1)(xy+y+1) & = & 1
\end{eqnarray}
Thus the conjunction polynomial $q^*$ is given by:
\begin{equation}
q^* \quad := \quad (xy+x+1) (xy+y+1) - 1 \quad \Rightarrow \quad x+y
\end{equation}
using modular arithmetic.  Definition~\ref{def:inverse} assures us
that for any query $s$ the inverse-value set computed from the
conjunction polynomial is the same as the original:
\begin{eqnarray}
  \set{S}_{\{x+y=0\}}^{-1}(s) & = & 
  \set{S}_{\{xy+x=0,\;xy+y=0\}}^{-1}(s)
\end{eqnarray}
Now, Lemma~\ref{lem:singleton} says that every polynomial in the
inverse-value set $\set{S}_{\{x+y=0\}}^{-1}(\{0\})$ must be the
product of some polynomial $p$ in $\mathbb{F}_2[x,y]$ with the
conjunction polynomial $q^*=x+y$ just computed:
\begin{eqnarray}
\set{S}_{\{x+y=0\}}^{-1}(\{0\}) & = & 
\left\{ p \times (x+y) : p \in \mathbb{F}_2[x,y] \right\}
\end{eqnarray}
And indeed, multiplying any polynomial in the ring $\mathbb{F}_2[x,y]$
by the sum $x+y$ will yield one of the four polynomials in
$\set{S}_A^{-1}(\{0\})$ listed in Equation~\ref{eq:theta-0}.  For
example we have the following products (keeping in mind that the
finite field $\mathbb{F}_2$ uses integer arithmetic modulo $2$ in
which $x^2=x$, $1+1=0$, etc.):
\begin{eqnarray}
(y+1) (x+y) & \Rightarrow & xy + x \\
(xy) (x+y) & \Rightarrow & 0 \\
(xy+y)(x+y) & \Rightarrow & xy+y \\
(xy+x+y)(x+y) & \Rightarrow & x+y
\end{eqnarray}
Furthermore, Corollary~\ref{cor:increment} states that for every
polynomial $p$ in the inverse-value set
$\set{S}_{\{x+y=0\}}^{-1}(\{0\})$, the sum $p+1$ must be a member of
the related inverse-value set $\set{S}_{\{x+y=0\}}^{-1}(\{1\})$.  You
can see this illustrated by comparing Equations \ref{eq:theta-0} and
\ref{eq:theta-1}.

\end{example}

\subsection{How to Solve It}
\label{sec:solve}

The calculations required for formula evaluation and formula discovery
using systems of polynomial equations with coefficients in finite
fields can be performed by an extension of the truth-table methods of
\cite{post} and \cite{wittgenstein}.  It is practical to carry out
these calculations by hand with pen and paper for small problems.
Alternatively, there are sophisticated algebraic geometry methods
implemented in several widely-available computer algebra systems that
accomplish the necessary calculations in a more efficient, robust, and
scalable manner (using coefficients from finite fields as easily as
rational or real coefficients).  These computational methods are
derived from the Gr\"obner-basis algorithms for solving polynomial
equations that were invented by Buchberger in the 1970s
\cite{buchberger}.  Commands to perform the requisite calculations in
a computer algebra system are included in
Section~\ref{sec:translation}.  Table-based inference is explained
presently.

For a system of $m$ polynomial equations with $n$ variables and
coefficients in a finite field with $d$ elements, naive table-based
inference requires the construction of a table with $d^{d^n}$ rows and
$d^n$ columns.  The resulting $d^{d^n+n}$ entries must be exhaustively
enumerated.  The following definition introduces the \emph{value
  worksheet} data structure and an inference algorithm based on it.
Each row of the value worksheet contains the same information as a
traditional logical truth table, but in a flattened form.

\begin{definition}[Table-Based Inference with Finite Fields]
\label{def:inference}
For a system of polynomial equations with coefficients in a finite
field $\mathbb{F}_d$ of some prime order $d$, the solution set,
solution-value sets, and inverse-value sets described in Definitions
\ref{def:polynomial-system}, \ref{def:objective}, and
\ref{def:inverse} can be computed using the following algorithm.
\begin{enumerate}
\item Construct the value worksheet $\mathbf{W}(A)$:
  \begin{enumerate}
  \item Make a table $\mathbf{W}(A)$ with a row for every polynomial
    $p_i \in \mathbb{F}_d[\mathbf{x}]$ and a column for every point
    $\mathbf{x}_k \in (\mathbb{F}_d)^n$.  Following
    Lemma~\ref{lem:count}, this table $\mathbf{W}(A)$ will have
    $d^{d^n}$ rows and $d^n$ columns (excepting the row and column
    labels).
  \item For every polynomial $p_i \in \mathbb{F}_d[\mathbf{x}]$ and
    every point $\mathbf{x}_k \in (\mathbb{F}_d)^n$, use the value
    $p_i(\mathbf{x}_k)$ of that polynomial evaluated at that point
    as the entry at row $i$ and column $k$ of the table.
  \end{enumerate}
\item Compute the solution set $\set{V}(A)$ from the constraint
  polynomials $A$:
  \begin{enumerate}
  \item If $d=2$ and $m \neq 1$, replace the set $A$ of equations with
    the singleton $\{q^*=0\}$ using the conjunction polynomial
    described in Lemma~\ref{lem:conjunction} (this step may be omitted
    for pedagogic purposes).
  \item For each constraint $(q_j=0) \in A$ find the matching
    polynomial $p_i=q_j$ at row $i$ in the table $\mathbf{W}(A)$.
    Examining this row $i$ of the table, mark as infeasible every
    column $k$ for which the value $p_i(\mathbf{x}_k)$ at row
    $i$ and column $k$ is not zero.  If there are no constraint
    polynomials then all columns remain unmarked.
  \item The unmarked columns identify the feasible points that
    constitute the solution set $\set{V}(A)$.  If all columns are
    marked as infeasible then the solution set is empty.
  \end{enumerate}
\item Compute the solution-value set $\set{S}_A(p_i)$ for every possible
  polynomial objective $p_i \in \mathbb{F}_d[\mathbf{x}]$:
  \begin{enumerate}
  \item Add a column labeled $\set{S}_A(p_i)$ to the table
    $\mathbf{W}(A)$ for solution-value sets.
  \item For each polynomial $p_i$ at row $i$, the solution-value set
    is the set of table entries at the unmarked columns.  Thus the
    solution-value set for each polynomial $p_i$ is the set of every
    value $p_i(\mathbf{x}_k)$ taken by that polynomial at a
    feasible point $\mathbf{x}_k$ indicated as an unmarked column
    $k$.
  \item If there are no unmarked columns then every solution-value set
    is empty.
  \end{enumerate}
\item Compute the inverse-value set $\set{S}_A^{-1}(s)$ for every
  possible solution-value set query $s \subseteq \mathbb{F}_d$:
  \begin{enumerate}
  \item Examine the entries in the solution-value set column
    $\set{S}_A(p_i)$.  For every unique entry $s$ in that column,
    select the rows $i_1$, $i_2$, etc.\ such that
    $\set{S}_A(p_{i_1})=s$, $\set{S}_A(p_{i_2})=s$, and so on.  The
    corresponding polynomials $p_{i_1}$, $p_{i_2}$, etc.\ constitute
    the inverse-value set $\set{S}_A^{-1}(s)$.
  \item The inverse-value set for every query $s \subseteq
    \mathbb{F}_d$ that does not appear as an entry in the
    solution-value set column $\set{S}_A(p_i)$ is empty.
  \end{enumerate}
\end{enumerate}
\end{definition}

\begin{example}\exfmt

\begin{table}
\caption{The value worksheet $\mathbf{W}(A)$ for the equations $A
  = \{xy+x=0,\;xy+y=0\}$ using the polynomial ring
  $\mathbb{F}_2[x,y]$.  The points $(0,1)$ and $(1,0)$ are marked
  infeasible; the solution set $\set{V}(A)_{\mathbb{F}_2}$ contains
  the remaining points $(0,0)$ and $(1,1)$.}
  \label{tbl:worksheet}
  \begin{tabular}{r@{\qquad}l@{\qquad}c@{\quad}c@{\quad}c@{\quad}c@{\qquad}c}
    \hline\hline
    $i$ & $p_i(x,y)$ & $p_i(0,0)$ & $p_i(0,1)$ & $p_i(1,0)$ & $p_i(1,1)$ & 
    $\set{S}_{\{xy+x=0,\;xy+y=0\}}(p_i)$ \\ \hline
    $1$ & $0$ & $0$ & $0$ & $0$ & $0$ & $\{0\}$ \\
    $2$ & $1$ & $1$ & $1$ & $1$ & $1$ & $\{1\}$ \\
    $3$ & $y$ & $0$ & $1$ & $0$ & $1$ & $\{0,1\}$ \\
    $4$ & $y+1$ & $1$ & $0$ & $1$ & $0$ & $\{0,1\}$ \\
    $5$ & $x$ & $0$ & $0$ & $1$ & $1$ & $\{0,1\}$ \\
    $6$ & $x+1$ & $1$ & $1$ & $0$ & $0$ & $\{0,1\}$ \\
    $7$ & $x+y$ & $0$ & $1$ & $1$ & $0$ & $\{0\}$ \\
    $8$ & $x+y+1$ & $1$ & $0$ & $0$ & $1$ & $\{1\}$ \\
    $9$ & $xy$ & $0$ & $0$ & $0$ & $1$ & $\{0,1\}$ \\
    $10$ & $xy+1$ & $1$ & $1$ & $1$ & $0$ & $\{0,1\}$ \\
    $11$ & $xy+y$ & $0$ & $1$ & $0$ & $0$ & $\{0\}$ \\
    $12$ & $xy+y+1$ & $1$ & $0$ & $1$ & $1$ & $\{1\}$ \\
    $13$ & $xy+x$ & $0$ & $0$ & $1$ & $0$ & $\{0\}$ \\
    $14$ & $xy+x+1$ & $1$ & $1$ & $0$ & $1$ & $\{1\}$ \\
    $15$ & $xy+x+y$ & $0$ & $1$ & $1$ & $1$ & $\{0,1\}$ \\
    $16$ & $xy+x+y+1$ & $1$ & $0$ & $0$ & $0$ & $\{0,1\}$ \\ \hline
    \multicolumn{2}{c}{Infeasible?} &
    $\Box$ & $\boxtimes$ & $\boxtimes$ & $\Box$ \\ \hline\hline
  \end{tabular}
\end{table}

Table~\ref{tbl:worksheet} shows the value worksheet $\mathbf{W}(A)$
constructed according to Definition~\ref{def:inference} from the
system of constraints $A = \{xy+x=0,\;xy+y=0\}$ in
Equation~\ref{eq:eqa1}; the polynomial ring $\mathbb{F}_2[x,y]$ with
binary finite-field coefficients is used.  There are two ways to mark
infeasible points.  Using the original constraint polynomials we note
that $xy+x$ appears as $p_{13}$ and $xy+y$ appears as $p_{11}$ in
$\mathbf{W}(A)$.  Examining row~11 we see that $p_{11}(0,1) \neq 0$,
therefore we mark the column for $(0,1)$ infeasible.  In row~13 the
entry $p_{13}(1,0) \neq 0$ hence we mark the column for $(1,0)$
infeasible.  Alternatively, using the conjunction polynomial $q^* = x+y$
gives equivalent results.  The polynomial $x+y$ appears at row~7 in
$\mathbf{W}(A)$ and in this row both entries $p_7(0,1)$ and $p_7(1,0)$
have nonzero values; hence the corresponding columns are marked
infeasible.

Each solution-value set $\set{S}_A(p_i)$ at row~$i$ is the set
$\{p_i(0,0),p_i(1,1)\}$ of values taken by the polynomial $p_i$ at the
feasible points (the unmarked columns in $\mathbf{W}(A)$).  For
example in row~8 we have $p_8=x+y+1$ and both $p_8(0,0)=1$ and
$p_8(1,1)=1$; thus the solution-value set $\set{S}_A(x+y+1)
\Rightarrow \{1\}$.  The inverse-value sets are given by selected
rows: for example rows 1, 7, 11, and 13 share the solution-value set
$\{0\}$ so the inverse-value set $\set{S}_A^{-1}(\{0\})$ is
$\{\;p_{1},\;p_{7},\;p_{11},\;p_{13}\;\}$.  The inverse-value sets for
the queries $\{1\}$ and $\{0,1\}$ are constructed in a similar way.
However, since no polynomial has an empty solution-value set (the
entry $\{\}$ does not appear in the last column of the value
worksheet), the inverse-value set of the empty set is itself empty:
$\set{S}_A^{-1}(\{\}) \Rightarrow \{\}$.

\end{example}

\section{Dynamical Systems from References to Solutions}
\label{sec:dynamic}

Now let us consider the complexity added to a system of equations when
it is permitted to refer to the solution of that system of equations
within the system itself.  This type of self-reference is an elaborate
\emph{recurrence relation}, which can either be interpreted as such to
define a \emph{discrete dynamical system}, or interpreted in a static
way to provide an additional simultaneous equation.  Although much of
the terminology and notation for `dynamical systems' is relatively new
(following a resurgence of interest in the 1970s, especially in
nonlinear and chaotic dynamical systems), the mathematical treatment
of recursion is quite old.  For example, what we know as the Fibonacci
sequence has been studied in various guises since at least the Middle
Ages; and methods for \emph{finite differences} and \emph{difference
  equations} have been developed since the work of Newton and then
Taylor around the turn of the 18th century \cite{bell}.  The modern
treatment of dynamical systems \cite{galor,holmgren} dates from
Poincar\'e's work at the end of the 19th century.  It happens that the
dynamical systems that we will encounter in the study of logic are
quite simple: discrete time, finite phase-space, first-order,
autonomous, and usually linear.

\subsection{Extended Systems of Polynomial Equations}

In order to develop a computable representation of solution
self-reference, we introduce several new features to the systems of
polynomial equations described in Section~\ref{sec:polynomial}:
parameters, equation templates, iterative assignments, and solution
references.  The idea is that each constraint may be specified as a
template instead of as a simple polynomial equation; the parameters of
the templates are allowed to refer to solution sets and solution-value
sets from the systems of equations in which they reside.

\begin{definition}[Parametric System of Polynomial Equations]
  \label{def:parametric}
  Consider a system of polynomial equations as specified in
  Definition~\ref{def:polynomial-system}, with variables $\mathbf{x} :=
  (x_1,x_2,\ldots,x_n)$ and coefficients in an algebraic field $K$;
  the polynomials in the system are members of the ring
  $K[\mathbf{x}]$.  We introduce a tuple $\Theta :=
  (\theta_1,\theta_2,\ldots,\theta_d)$ of parameters with the
  requirement that each parameter $\theta_j$ must take a value in some
  specified set $U_j$.  Thus the set $U$ of possible values for
  $\Theta$ is given by the Cartesian product:
  \begin{eqnarray*}
    U & := & U_1 \otimes U_2 \otimes \cdots \otimes U_\ell
  \end{eqnarray*}
  For each variable $x_i$ and each parameter
  $\theta_j$ we add a data-type constraint $\tau$ to identify the
  appropriate set of possible values:
  \begin{displaymath}
  \begin{array}{rrcl}
    \tau_{x_i} : & x_i & \in & K \\
    \tau_{\theta_j} : & \theta_j & \in & U_j
  \end{array}
  \end{displaymath}
  Instead of using a simple set $A$ of polynomial equations as
  described in Definition~\ref{def:polynomial-system}, an extended
  polynomial system is specified using a set $A(\Theta)$ of parametric
  \emph{constraint templates} and a set $D$ of \emph{assignment
    templates}, in addition to the above data-type constraints.  In
  the set $A(\Theta) := \{\alpha_1, \alpha_2, \ldots, \alpha_m \}$ each
  constraint template $\alpha_i$ is a function of the parameters
  $\Theta$:
  \begin{eqnarray*}
    \alpha_i : \quad q_i(\Theta) & = & 0
  \end{eqnarray*}
  such that for any instantiation of the parameter values $\Theta$ the
  template $\alpha_i$ reduces to an ordinary polynomial equation with
  $q_i(\Theta) \in K[\mathbf{x}]$.  In the set $D :=
  \{\delta_{\theta_1},\delta_{\theta_2},\ldots,\delta_{\theta_\ell}\}$
  each assignment template $\delta_j$ specifies the value to be
  assigned to the corresponding parameter $\theta_j$ at every
  iteration, using some \emph{parameter-updating function} $\lambda_j$
  whose arguments may include the solution set $\set{V}(A(\Theta))$ to
  the extended system of equations being defined, as well as any
  parameter including $\theta_j$ itself (here the double left arrow
  $\Leftarrow$ denotes assignment):
  \begin{eqnarray*}
    \delta_{\theta_j} : \quad 
    \theta_j & \Leftarrow & 
    \lambda_j(\;
    \theta_1,\theta_2,\ldots,\theta_\ell\;,\;\set{V}(A(\Theta))\;)
  \end{eqnarray*}
  In particular a parameter-updating function $\lambda_j$ may use the
  solution-value set $\set{S}_{A(\Theta)}(p)$ for some objective
  formula $p$ or the cardinality $|\set{S}_{A(\Theta)}(p)|$ of such a
  solution-value set.  Although every parameter-updating function
  $\lambda_j$ that assigns a value to a parameter $\theta_j$ must
  return a value in the designated set $U_j$ of possible values for
  that parameter, the updating functions are not required to be
  polynomial.  If an explicit updating function for any parameter
  $\theta_j$ is omitted, then the identity function $\lambda_j(\theta)
  : \theta$ is used as a default.
\end{definition}

According to Definition~\ref{def:parametric} the self-referential
quadratic system in Equations \ref{eq:cc-orig} and
\ref{eq:cc-orig-ref} can be specified as the following constraint and
assignment templates:
\begin{equation}
  \begin{array}{rrcl}
    \tau_x : & x & \in & \mathbb{R} \\
    \tau_c : & c & \in & \{0,1,2\} \subset \mathbb{R} \\
    \alpha_1 : & 2x^2 + 3x + c & = & 0 \\
    \delta_c : & c & \Leftarrow & \left|\set{S}_{\{\alpha_1\}}(x)\right|
  \end{array}
  \label{eq:c-template}
\end{equation}
Here the variable $x$ takes real values and the possible values of the
single parameter $c$ are the set $U:=\{0,1,2\}$ of integers between
$0$ and $2$.  This parameter $c$ is assigned the cardinality of the
solution-value set for the objective $x$ (subject to the equation
$\alpha_1$) by the updating function $\lambda_1(c)$ specified by the
assignment template $\delta_c$.

For a different example, a system to generate Fibonacci-like sequences
can be specified with two parameters and a pair of assignment
templates, without the use of any conventional variables or any
constraint templates:
\begin{equation}
  \begin{array}{rrcl}
    \tau_a : & a & \in & \mathbb{Z} \\
    \tau_b : & b & \in & \mathbb{Z} \\
    \delta_a : & a & \Leftarrow & b \\
    \delta_b : & b & \Leftarrow & a+b
  \end{array}
  \label{eq:fib-template}
\end{equation}

We next develop an algorithm to make explicit the discrete dynamical
system implied by a parametric system of polynomial equations.

\begin{definition}[Extracting the Evolution Function]
  \label{def:extract}
  A parametric system of polynomial equations as described in
  Definition~\ref{def:parametric} encodes a functional relationship
  between the value of the parameter $\Theta$ and itself.  Here we
  extract this \emph{evolution function}, denoted $F(\Theta)$, in two
  useful special cases.  First, in the case that each
  parameter-updating function $\lambda_j$ is a simple algebraic
  function of the parameters in $\Theta$---without reference to the
  solution set $\set{V}(A(\Theta))$---then the evolution function $F$
  is given by the tuple of parameter-updating functions:
  \begin{eqnarray*}
    F(\theta_1,\theta_2,\ldots,\theta_\ell) & : & \left( 
    \lambda_1(\Theta), \lambda_2(\Theta), \ldots, \lambda_\ell(\Theta) \right)
  \end{eqnarray*}
  This is the usual situation with recurrence relations in general
  algebra.  

  Second, in the case that the parameter domain $U$ is countable and
  finite then we can proceed with \emph{hypothetico-deductive
    analysis} to derive $F$ by computing a mapping $\Theta_i \mapsto
  \Theta'_i$ for every parameter value $\Theta_i \in U$.  To express
  the hypothesis that the parameters $\Theta$ have some particular
  value $\Theta_i$, each parameter $\theta_j$ is assigned a constant
  value $\widehat{\theta_j}$ from its domain $U_j$.  Using these
  values the constraint templates are instantiated into an ordinary
  set $A(\Theta_i)$ of polynomial equations:
  \begin{eqnarray*}
    A(\Theta_i) & := & \{ \; q_1(\Theta_i)=0,\;
    q_2(\Theta_i)=0,\; \ldots, \; q_m(\Theta_i)=0 \; \}
  \end{eqnarray*}
  The solution set $\set{V}(A(\Theta_i))$ for these instantiated
  equations is computed by the table-based algorithm in
  Definition~\ref{def:inference} or by a general algebraic geometry
  method as appropriate.  Following this solution subroutine the
  assignment templates in the set $D$ are processed: an updated value
  $\theta'_j$ is computed for every parameter via its
  parameter-updating function, using as arguments the solution-value
  set $\set{V}(A(\Theta_i))$ just computed along with the input values
  of the parameters:
  \begin{eqnarray*}
    \theta'_j & \Leftarrow & 
    \lambda_j(\; 
    \widehat{\theta_1}, \widehat{\theta_2}, \ldots, \widehat{\theta_\ell},\;
    \set{V}(A(\Theta_i))
    \;)
  \end{eqnarray*}
  The desired evolution function $F$ must map the hypothesized
  parameter values to their updated counterparts:
  \begin{eqnarray*}
    (\widehat{\theta_1}, \widehat{\theta_2}, \ldots, \widehat{\theta_\ell})
    & \mapsto &
    (\theta'_1, \theta'_2, \ldots, \theta'_\ell)
  \end{eqnarray*}
  The set of mappings $\Theta_i \mapsto \Theta_i'$ for every parameter
  value $\Theta_i \in U$ completely defines the evolution function
  $F$.  When the parameter domain $U$ is a subset of the field $K$
  used as polynomial coefficients, then the evolution function $F$ can
  be expressed in closed form using polynomial interpolation.
\end{definition}
To illustrate the first case in Definition~\ref{def:extract}, the
evolution function for the Fibonacci-like parametric system in
Equation~\ref{eq:fib-template} is simply:
\begin{eqnarray}
  F(a,b) & : & (b,a+b)
  \label{eq:fib-transit}
\end{eqnarray}
using the parameter-updating functions $\lambda_1(a):b$ and
$\lambda_2(b):a+b$ specified by the assignment templates $\delta_a$
and $\delta_b$ in Equation~\ref{eq:fib-template}.

\begin{table}
  \caption{The state-transition worksheet for the dynamical system
    defined by the quadratic equation $2x^2+3x+c=0$ with $c$
    solutions, as modeled by the parametric templates in
    Equation~\ref{eq:c-template}.  Each value $c_i$ specifies a state
    of the system.  The set $A(c_i)$ contains the instantiation of the
    constraint template $\alpha_1$ at that state $c_i$.  The
    solution-value set $\set{S}_{A(c_i)}(x)$ gives the solutions for
    the objective $x$ subject to the instantiated constraint.  The
    successor state $F(c_i)$ is calculated by processing the
    assignment template $\delta_c$ which asks the size of this
    solution-value set $\set{S}_{A(c_i)}(x)$.}
  \label{tbl:sst-c}
  \begin{tabular}{r@{\qquad}c@{\quad}c@{\qquad}c@{\qquad}c} 
    \hline\hline
    $i$ & $c_i$ & $A(c_i)$ & $\set{S}_{A(c_i)}(x)$ & $F(c_i)$ \\ \hline
    $1$ & $0$ & $\left\{2x^2+3x+0=0\right\}$ 
    & $\{-\frac{3}{2},0\}$ & $2$ \\
    $2$ & $1$ & $\left\{2x^2+3x+1=0\right\}$ 
    & $\{-1,-\frac{1}{2}\}$ & $2$ \\
    $3$ & $2$ & $\left\{2x^2+3x+2=0\right\}$ & $\{\}$ & $0$ \\ \hline\hline
  \end{tabular}
\end{table}

The requisite calculations for the second case in
Definition~\ref{def:extract} can be organized in a
\emph{state-transition worksheet} constructed as follows.  In each row
of the state-transition worksheet we record an index $i$, some value
$\Theta_i$ from the set $U$, the instantiation $A(\Theta_i)$ of the
constraint templates at that value, the relevant feature of the
solution-value set $\set{V}(A(\Theta_i))$ given those instantiated
equations, and the value $F(\Theta_i)$ of the successor that was
computed from $\set{V}(A(\Theta_i))$ and $\Theta_i$ by processing the
assignment templates according to Definition~\ref{def:extract}.  Each
row of the state-transition worksheet gives a specific value
$F(\Theta_i)$ of the evolution function for the argument $\Theta_i$.
The complete evolution function $F(\Theta)$ can be represented as a
simple transition matrix or as a state-transition table; or using the
method described in Lemma~\ref{lem:matrix} below, $F(\Theta)$ can be
specified as a polynomial with coefficients in a finite field.

Table~\ref{tbl:sst-c} shows the state-transition worksheet for the
parametric system in Equation~\ref{eq:c-template}, which uses a single
parameter $c$ with domain $U=\{0,1,2\}$.  In this case the relevant
feature of the solution set is the solution-value set
$\set{S}_{A(c_i)}(x)$ for the objective $x$.  This worksheet
constructs the following evolution function:
\begin{equation}
  \begin{array}{r@{\:}c@{\:}l}
    F & : & \{0,1,2\} \rightarrow \{0,1,2\} \\
    && 0 \mapsto 2 \\ && 1 \mapsto 2 \\ && 2 \mapsto 0
  \end{array}
  \label{eq:fc-map}
\end{equation}
The function-development method of Lemma~\ref{lem:matrix} below
enables us to construct a polynomial function that matches any
transition function $F:U \rightarrow U$ for a finite phase space $U$
(after assigning integers to identify the states if they are not
already numeric); it happens that $F(c):-c^2+c+2$ agrees with the
mappings in Equation~\ref{eq:fc-map} (using the polynomial ring
$\mathbb{R}[c]$ with real coefficients).

\begin{definition}[The Derived Dynamical System]
\label{def:dynamical}
A parametric system of polynomial equations as described in
Definition~\ref{def:parametric} in turn defines a discrete-time
dynamical system with state $\Theta$ whose phase space is the set $U$
of parameter values and whose evolution function is $F(\Theta)$ from
Definition~\ref{def:extract}.  The set $\{0,1,2,\ldots\}$ of
nonnegative integers is used as the domain for the evolution (time)
parameter $t$.  The state $\Theta_t$ of the dynamical system at any
time $t$ is the result of the $t$-fold composition of the evolution
function $F$ applied to the initial state $\Theta_0$:
\begin{eqnarray*}
  \Theta_t & \Leftarrow & F^{(t)}(\Theta_0)
\end{eqnarray*}
For example, given state $\Theta_0$ at time $t=0$ the state of the
system at time $t=3$ is the composed value
$F^{(3)}(\Theta_0)=F(F(F(\Theta_0)))$.  Following the usual
conventions for discrete dynamical systems, a sequence of successive
states $\Theta_t$, $F(\Theta_t)$, $F(F(\Theta_t))$, and so on
constitutes an \emph{orbit} \cite{galor,holmgren}.  Furthermore a
state $\Theta_{t^*}$ is a \emph{fixed point} in the dynamical system
exactly if it is its own successor: $\Theta_{t^*}=F(\Theta_{t^*})$.
We categorize a dynamical system by the number of fixed points it has.
Let us say that a dynamical system is \emph{steady} if it has one
fixed point; \emph{unsteady} if it has no fixed points; and
\emph{contingent} if it has more than one fixed point.
\end{definition}
Note that this categorization concerns the state $\Theta$ of the
dynamical system indicated by the parameters
$(\theta_1,\theta_2,\ldots,\theta_\ell)$ rather than the solution-value
set for any objective formula \emph{per se}; it could happen in an
unsteady system that the sequences of solution-value sets for some
formulas are nonetheless monotonous.  Note also that if a system of
equations has no parameters, then the dynamical system extracted
according to Definition~\ref{def:dynamical} will have the phase space
$U=\{()\}$ containing one state which is the empty tuple $()$.  In
this case the state-transition function $F():()$ is the identity
function and the solitary state must be a fixed point; the system is
trivially steady.

By the categorization scheme in Definition~\ref{def:dynamical} the
dynamical system derived from Equation~\ref{eq:c-template} for the
quadratic equation $2x^2+3x+c=0$ is \emph{unsteady}; as the graph in
Equation~\ref{eq:cc-graph} shows there is a periodic cycle and there
are no fixed points.  In contrast the Fibonacci-like system derived
from Equation~\ref{eq:fib-template} is \emph{steady}; as we shall see
there is one fixed point at $(a,b)=(0,0)$ even though the orbits
through the other points do not converge.

\subsection{Static and Dynamic Interpretations}

A dynamical system from Definition~\ref{def:dynamical} can be
interpreted in two different ways to evaluate an objective formula.
In the dynamic interpretation, the orbits in the dynamical system are
used to generate infinite sequences of solution-value sets for the
objective (a sequential report).  In the static interpretation, the
evolution function is used to provide an additional static constraint
limiting attention to the fixed points; what is reported is the union
of solution-value sets for the objective from these fixed points (a
stationary report).  Such static and dynamic interpretations of a
self-referential system of equations reflect subtly different views;
it is fine to try either or both for any given problem.

\begin{definition}[Collected Sequences of Solution-Value Sets]
  \label{def:sequence}
  Consider a dynamical system derived as in Definitions
  \ref{def:extract} and \ref{def:dynamical} from a parametric
  polynomial system of equations as described in
  Definition~\ref{def:parametric}.  For any given objective formula $p
  \in K[\mathbf{x}]$ or $p \in \Theta$, the dynamical system specified
  by $U$ and $F(\Theta)$ encodes an infinite sequence of
  solution-value sets for $p$ using each state $\Theta_0 \in U$ as an
  initial condition.  Each infinite sequence $\vec{\set{S}}$ of
  solution-value sets results from solving the systems of equations
  $A(\Theta_t)$ instantiated from the constraint templates at
  successive states $\Theta_0, \Theta_1, \Theta_2, \ldots$:
  \begin{eqnarray*}
    \vec{\set{S}}_{A(\Theta)}(p \:|\: \Theta_0) & := & \left( \;
    \set{S}_{A(\Theta_0)}(p), \;
    \set{S}_{A(\Theta_1)}(p), \;
    \set{S}_{A(\Theta_2)}(p), \; \ldots
    \; \right)
  \end{eqnarray*}
  With reference to Definition~\ref{def:extract}, when the objective
  $p$ refers to a parameter $\theta_j$ then the initial (hypothesized)
  value $\widehat{\theta_j}$ for the current state should be used to
  evaluate that objective, instead of the updated parameter value
  $\theta'_j$.

  The \emph{collection} $\vec{\set{S}}^*$ of solution-value set
  sequences is defined as a relation between the states in $U$ and the
  infinite sequences that arise from them.  For every state $\Theta_i
  \in U$ the collection contains a mapping to the sequence with the
  respective initial condition:
  \begin{eqnarray*}
    \vec{\set{S}}_{A(\Theta)}^*(p) & := & \left\{ \;
    \Theta_t \; \mapsto \;
    \vec{\set{S}}_{A(\Theta)}(p \:|\: \Theta_t) \; : \; \Theta_t \in U^*
    \; \right\}
  \end{eqnarray*}
  Such a collection is usually written in square brackets as
  illustrated in the examples.  If every solution-value set in a
  sequence has exactly one member (which is always the case when the
  objective is a parameter), then the braces around those sets may be
  omitted to simplify notation: thus the sequence
  $(\{k_0\},\{k_1\},\{k_2\},\ldots)$ could instead be written
  $(k_0,k_1,k_2,\ldots)$.
\end{definition}
For the parametric system in Equation~\ref{eq:c-template},
Equation~\ref{eq:cc-seq-c} in the introduction already illustrated the
collection $\vec{\set{S}}_{\{\alpha_1\}}^*(c)$ of solution-value-set
sequences for the objective formula $c$ (with the set brackets omitted
since each solution-value set is a singleton).  Likewise
Equation~\ref{eq:cc-seq-x} shows the collection
$\vec{\set{S}}_{\{\alpha_1\}}^*(x)$ of solution-set sequences for the
objective $x$.  Note that you can read the sequences in either
collection from the graph in Equation~\ref{eq:cc-graph}.

Considering the Fibonacci example, the parametric system in
Equation~\ref{eq:fib-template} specifies a dynamical system with phase
space $U=\mathbb{Z}^2$, state $\Theta=(a,b)$, and evolution function
$F(a,b):(b,a+b)$.  Using this dynamical system the collection of
solution-value-set sequences for the objective parameter $a$ is given
by:
\begin{eqnarray}
  \vec{\set{S}}_{\{\}}^*(a) & \Rightarrow &
  \left[ \begin{array}{rcl}
      (a,b)_0=(0,0) & \mapsto & (0,0,0,0,0,0,0,\ldots) \\
      (a,b)_0=(0,1) & \mapsto & (0,1,1,2,3,5,8,\ldots) \\
      (a,b)_0=(2,1) & \mapsto & (2,1,3,4,7,11,18,\ldots) \\
      & \vdots &
  \end{array} \right]
\end{eqnarray}
The initial condition $(a,b)_0=(0,1)$ yields the familiar Fibonacci
sequence.  Again since each solution-value set is a singleton the
braces around sequence elements have been omitted.

The dynamical system can instead be interpreted in a static way to
examine its fixed points.

\begin{definition}[Simultaneous Solution for Fixed Points]
  \label{def:flat}
  Consider a parametric system of polynomial equations as described in
  Definition~\ref{def:parametric} and its extracted evolution function
  $F(\Theta)$ derived by Definition~\ref{def:extract}, in the special
  case that the domain $U_j$ for each parameter $\theta_j$ in
  $\Theta:=(\theta_1,\theta_2,\ldots,\theta_\ell)$ is a subset of the
  algebraic structure $K$ used for the conventional variables and
  polynomial coefficients in the problem.  In this special case the
  solutions from the fixed points in the dynamical system can be
  computed in an alternative way by solving an extended system of
  static equations in which each parameter $\theta_j$ is treated as an
  indeterminate along with the conventional variables in
  $\mathbf{x}:=(x_1,x_2,\ldots,x_n)$.  The extended static system
  includes the data-type constraints $\tau$ and a polynomial equation
  transcribed from each each constraint template $\alpha$ and from
  each assignment template $\delta$:
  \begin{displaymath}
    \begin{array}{rrcl}
      \tau_{x_j} : & x_j & \in & K \\
      \tau_{\theta_j} : & \theta_j & \in & U_j \subseteq K \\
      \alpha_i : & q_i & = & 0 \\
      \delta_{\theta_j} : & \theta_j & = & \lambda_j( \Theta, \set{V}(A) )
    \end{array}
  \end{displaymath}
  Each parametric constraint polynomial $q_i$ is now treated as a
  member of the extended polynomial ring:
  \[K[x_1,x_2,\ldots,x_n,\theta_1,\theta_2,\ldots,\theta_\ell]\]
\end{definition}
Note that in this alternative formulation, fixed points that give
empty solution-value sets are not reported.

By Definition~\ref{def:flat} the parametric quadratic system in
Equation~\ref{eq:c-template} yields the following system of static
equations:
\begin{equation}
  \begin{array}{rrcl}
    \tau_x : & x & \in & \mathbb{R} \\
    \tau_c : & c & \in & \{0,1,2\} \subset \mathbb{R} \\
    \alpha_1 : & 2x^2 + 3x + c & = & 0 \\
    \delta_c : & c & = & -c^2+c+2
  \end{array}
\end{equation}
using the polynomial interpolation of the extracted evolution function
$F(c)$ shown in Equation~\ref{eq:fc-map}.  This system of equations is
infeasible; thus the solution-value set
$\set{S}_{\{\tau_x,\tau_c,\alpha_1,\delta_c\}}(x)$ evaluates to the
empty set, as does the solution-value set for any other objective
formula subject to these unsatisfiable constraints.

For the Fibonacci-like system in Equation~\ref{eq:fib-template},
transcribing the assignment templates into static equations yields the
following set $A$ of equations (using the extended polynomial ring
$\mathbb{Z}[a,b]$ with integer coefficients):
\begin{equation}
  \begin{array}{rrcl}
    \tau_a   : & a & \in & \mathbb{Z} \\
    \tau_b   : & b & \in & \mathbb{Z} \\
    \delta_a : & a & = & b \\
    \delta_b : & b & = & a+b
  \end{array}
\end{equation}
The unique solution $(a,b)=(0,0)$ identifies the fixed point noted
earlier.  Consequently for the objective $a$ the solution-value set
$\set{S}_{\{\tau_a,\tau_b,\delta_a,\delta_b\}}(a) \Rightarrow \{0\}$.
Interpreting its familiar recurrence as a static constraint, it is
appropriate to say that `\emph{the} Fibonacci number' is zero.

\begin{example}\exfmt
Consider a new quadratic equation $\frac{1}{2}x^2+3bx+\frac{11}{2}=0$
in which the coefficient $b$ is defined to be the number of real
solutions to the equation in which it appears.  This problem is
specified as the following parametric system, using the conventional
variable $x$ and the parameter $b$:
\begin{equation}
  \begin{array}{rrcl}
    \tau_x : & x & \in & \mathbb{R} \\
    \tau_b : & b & \in & \{0,1,2\} \subset \mathbb{R} \\
    \alpha_1 : & \frac{1}{2}x^2 + 3bx + \frac{11}{2} & = & 0 \\
    \delta_b : & b & \Leftarrow & \left| \set{S}_{\{\alpha\}}(x) \right|
  \end{array}
  \label{eq:b-template}
\end{equation}
By Definition~\ref{def:dynamical} these templates yield a dynamical
system with the phase space $U=\{0,1,2\}$ and an evolution function
$F(b)$ that must satisfy the following criteria:
\begin{equation}
  F(0) = 0, \quad F(1) = 0, \quad F(2) = 2
\end{equation}
The polynomial $F(b):b^2-b$ (constructed by Lemma~\ref{lem:matrix})
meets these criteria for $b\in\{0,1,2\}$.  The derived dynamical
system, labeled with the solution-value sets discussed next, is
displayed in this graph:
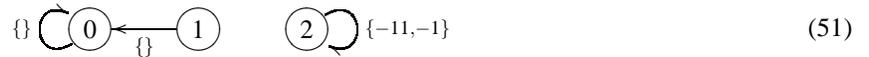
\begin{equation}
  \xymatrix{
    *++[o][F]{0} \ar@(dl,ul)[]^{\{\}} & 
    *++[o][F]{1} \ar[l]^{\{\}} &
    *++[o][F]{2} \ar@(ur,dr)[]^{\{-11,-1\}}
  }
  \label{eq:b-graph}
\end{equation}
There are fixed points at $b=0$ and $b=2$ and no nonconvergent orbits;
in the terminology of Definition~\ref{def:dynamical} this
dynamical system \emph{contingent}.  By Definition~\ref{def:sequence}
the dynamical system yields the following collection of solution-value
sets for the objective $x$:
\begin{eqnarray}
  \vec{\set{S}}_{\{\alpha_1\}}^*(x) & \Rightarrow &
  \left[ \begin{array}{rcl}
      b_0 = 0 & \mapsto & (\{\},\{\},\ldots) \\
      b_0 = 1 & \mapsto & (\{\},\{\},\ldots) \\
      b_0 = 2 & \mapsto & (\{-11,-1\},\{-11,-1\},\ldots)
    \end{array} \right]
\end{eqnarray}
Reviewing the fixed points and their associated solution-value sets
shown in the graph in Equation~\ref{eq:b-graph}, the dynamic
interpretation shows that there are two cases in which the equation
$\frac{1}{2}x^2+3bx+\frac{11}{2}=0$ has exactly $b$ real roots: when
$b=0$ and there are no solutions, and when $b=2$ and there are two
solutions (namely $-11$ and $-1$).

Finally, using Definition~\ref{def:flat} and the polynomial
interpolation $F(b):b^2-b$ of the evolution function, the templates in
Equation~\ref{eq:b-template} can be transcribed into a static set of
equations using polynomials in the extended ring $\mathbb{R}[x,b]$:
\begin{equation}
  \begin{array}{rrcl}
    \tau_x : & x & \in & \mathbb{R} \\
    \tau_b : & b & \in & \{0,1,2\} \subset \mathbb{R} \\
    \alpha_1 : & \frac{1}{2}x^2 + 3bx + \frac{11}{2} & = & 0 \\
    \delta_b : & b & = & b^2-b
  \end{array}
  \label{eq:b-flat}
\end{equation}
The static polynomial system in Equation~\ref{eq:b-flat} has two
solutions $(x,b)=(-11,2)$ and $(x,b)=(-1,2)$; these identify the
fixed-point solutions.  Consequently in this static formulation the
solution-value set for $x$ is given by
$\set{S}_{\{\tau_x,\tau_b,\alpha_1,\delta_b\}}(x) \Rightarrow
\{-11,-1\}$.  Note that this static interpretation does not reveal
that $b=0$ gives a valid instantiation of the template equation
$\alpha_1$, since in that instantiation there is no real solution for
$x$.

\end{example}

\section{Translation from Logic into Algebra}
\label{sec:translation}

Let us now visit the foundations of logic and explore the relationship
between logical reasoning with truth values and mathematical reasoning
with ordinary numbers.  What we call \emph{logic}---reasoning with
binary truth values, formulas, axioms, theorems, and proof---is
nothing more than algebra presented in peculiar notation.  This basic
equivalence was asserted by Boole in his innovative mathematical
treatment of logic \cite{boole-mal,boole}.  Indeed, it seems
apparent that in both logic and algebra we find common concepts of
number, operation and formula; variable, function, equation, and
solution to equations; recursion and infinite sequence.  However,
through a cascade of historical misunderstanding and obfuscation (much
of it associated with \goedel's incompleteness argument), the
intrinsic unity of logic and mathematics has not achieved universal
acceptance.  Especially, recurrence relations and dynamical systems
have not been recognized as such in the contexts of logic and set
theory.  

Considering these things I have extended Boole's ideas, and the two
key principles presented in the introduction of this essay, into
several \emph{Articles of Algebraic Translation for Logic}:
\begin{enumerate}
  \raggedright
  \renewcommand{\theenumi}{\Roman{enumi}}
  \renewcommand\labelenumi{\theenumi.}
  \item Logical formula $\leadsto$ Polynomial formula
  \item Axiom $\leadsto$ Polynomial equation
  \item Solution-value set $\leadsto$ Truth value (when
    self-reference is forbidden)
  \item Polynomial with solution-value set $\{1\}$ $\leadsto$
    Theorem
  \item Reference to provability $\leadsto$ Recurrence relation
    \label{art:recurrence}
  \item Features of discrete dynamical system $\leadsto$ Truth
    value (when self-reference is allowed)
    \label{art:recurrent-value}
\end{enumerate}
Here the wavy arrow $\leadsto$ should be read `translates as,' with
the idea that these general principles subsume many specific
functional mappings.  By analogy the general principle:
\[\textrm{Roman number} \leadsto \textrm{Hindu-Arabic number}\] 
includes the specific mappings $\textrm{IIII} \leadsto 4$ and
$\textrm{IV} \leadsto 4$.  \goedel's formula \rqq, which denies its
own provability, is analyzed during the presentation of Articles
\ref{art:recurrence} and \ref{art:recurrent-value}.  For concreteness,
the examples in this section are accompanied by commands to perform
their calculations in the computer algebra system \emph{Mathematica}
\cite{mathematica}.  Carnielli has studied Boole's polynomial
formulation of mathematical logic and uses related ideas in his
Polynomial Ring Calculus
\cite{carnielli-polynomizing,agudelo-carnielli}.

\subsection{Logical Formulas as Polynomials}

Let us consider two methods to represent logical formulas as
polynomial expressions: the original method created by Boole, and a
revised method using finite fields and modular arithmetic.  In either
case we start with a logical formula in the propositional calculus
built from some atomic formulas $(x_1,x_2,\ldots,x_n)$ and the usual
unary and binary operations (negation, conjunction, disjunction,
material implication, and so on).  Any such formula has a particular
truth table, and it is truth tables that are directly translated into
polynomial functions; in fact, each translated polynomial can be
considered a closed-form representation of a truth table.  Many
different logical formulas might share a given truth table; hence many
different logical formulas might have the same polynomial translation.
Translation is not limited to $2$-valued logic; in general any finite
number of elementary truth values can be accommodated.  

Although I discuss polynomial \emph{translation}, for Boole
polynomials were not used to translate from some other symbolic
notation for logical formulas; they were his only mathematical
representation.  Notably, Frege's \emph{Begriffsschrift} was published
many years after Boole's lifetime, as were the subsequent works from
Peano and Hilbert in which the ideography of contemporary logical
notation was developed.

The detailed translation methods described below have been used to
generate a table of polynomials representing common logical formulas,
including the elementary truth values, the unary negation of a
formula, and the common binary operations applied to a pair of
formulas.  These translations, which are often called the Stone
isomorphisms after \cite{stone}, appear in
Table~\ref{tbl:translation}.  A practical way to translate a simple
logical formula is to apply the substitutions listed in the table
recursively, until all traditional logical operations have been
converted into polynomial form.  Using Table~\ref{tbl:translation} it
is possible to translate any formula from the propositional calculus
into a polynomial with either real-number coefficients, or
coefficients in the binary finite field $\mathbb{F}_2$.

For example the logical formula $y \wedge (z \oplus w)$, which is the
logical conjunction of the atomic proposition $y$ with the exclusive
disjunction (\textsc{xor}) of the atomic propositions $z$ and $w$, can
be translated to a polynomial in the ring $\mathbb{R}[w,y,z]$ by first
substituting the polynomial form of the inner disjunction, then
substituting the outer conjunction, and finally simplifying the whole
expression using standard algebra:
\begin{eqnarray}
  y \wedge (z \oplus w) & \leadsto & y \wedge (z+w-2zw) \\
  & \leadsto & (z+w-2zw) \\
  & \Rightarrow & yz + yw - 2yzw
\end{eqnarray}
Besides translating the entire formula $y \wedge (z \oplus w)$, each
variable in it must be constrained to limit its possible values to $0$
and $1$.  This is accomplished by the equations $w^2=w$, $y^2=y$, and
$z^2=z$ (which would be tautological if we were using constraints in
$\mathbb{F}_2$).

\begin{table}
  \caption{Translation from logical to polynomial notation, to map
    logical formulas $p$ and $q$ into polynomials with real or
    finite-field coefficients.  The atomic formulas are represented by
    variables $\mathbf{x}:=(x_1,x_2,\ldots,x_n$).  When using real
    coefficients each $x_i$ is subject to the constraint $x_i^2=x_i$
    to limit its possible values to $\{0,1\}$.}
  \label{tbl:translation}
  \begin{tabular}{l@{\quad}l@{\quad}l@{\quad}l} \hline\hline
    \bfseries Logical & 
    \bfseries Polynomial in $\mathbb{R}[\mathbf{x}]$ & 
    \bfseries Polynomial in $\mathbb{F}_2[\mathbf{x}]$ & 
    \bfseries Description \\ 
    \hline
    $\val{t}$ & $1$ & $1$ & True \\
    $\val{f}$ & $0$ & $0$ & False \\
    $\neg p$ & $1-p$ & $1+p$ & Logical negation (\textsc{not}) \\
    $p \wedge q$ & $pq$ & $pq$ & Conjunction (\textsc{and}) \\
    $p \oplus q$ & $p+q-2pq$ & $p+q$ & 
    Exclusive disjunction (\textsc{xor}) \\
    $p \vee q$ & $p+q-pq$ & $p+q+pq$ & 
    Nonexclusive disjunction (\textsc{or}) \\
    $p \rightarrow q$ & $1-p+pq$ & $1+p+pq$ & Material implication \\
    $p \leftrightarrow q$ & $1-p-q+2pq$ & $1+p+q$ & 
    Biconditional (\textsc{xnor}) \\
    $p \uparrow q$ & $1-pq$ & $1+pq$ & Nonconjunction (\textsc{nand}) \\
    $p \downarrow q$ & $1-p-q+pq$ & $1+p+q+pq$ & 
    Nondisjunction (\textsc{nor}) \\ 
    \hline\hline
  \end{tabular}
\end{table}

\subsubsection{Boole's Original Representation Scheme}
\label{sec:original}

Boole's method to represent logical formulas as polynomial formulas
was presented it in preliminary form in \emph{Mathematical Analysis of
  Logic} and in more complete form in \emph{Laws of Thought}
\cite{boole-mal,boole}.  This method produces polynomials with integer
coefficients that can be manipulated using ordinary arithmetic.  Boole
used the number $1$ to represent the logical value true, the number
$0$ to represent the logical value false, and symbolic variables such
as $x$, $y$, and $z$ to represent atomic logical propositions.  Boole
noted that if a variable $x$ is limited to the values $0$ and $1$,
then the equation $x^2=x$ must hold; he identified this equation as a
`special law.'  Based on this equation and its rearranged form
$x(1-x)=0$, Boole devised a clever method to develop polynomial
functions from arbitrary truth tables in two-valued logic.  Each
developed polynomial has the property that it agrees in value with its
logical predecessor for any combination of truth values of the atomic
propositions.

Boole's function-development method is illustrated here for functions
of two variables, which is the most useful case since traditional
logical notation uses unary and binary operations (negation,
conjunction, material implication, and so on).  Suppose that we have
two real-valued variables $x$ and $y$, each restricted to values in
the set $\{0,1\}$.  We require a polynomial function $p(x,y)$ that
yields some specified value $z_{1}$ (also either $0$ or $1$) when
$x=1$ and $y=1$; likewise we require $p(1,0)=z_{2}$, $p(0,1)=z_{3}$,
and $p(0,0)=z_{4}$, with each value $z_{i} \in \{0,1\}$.  These
required values can be arranged as illustrated on page~76 of
\cite{boole}, perhaps the earliest specimen of a logical truth table:
\begin{equation}
  \begin{array}{cc|c}
    x & y & p(x,y) \\ \hline
    1 & 1 & z_{1} \\
    1 & 0 & z_{2} \\
    0 & 1 & z_{3} \\
    0 & 0 & z_{4}
  \end{array}
  \label{eq:tbl-boole}
\end{equation}
Boole proved that the requisite function $p(x,y)$ can always be
calculated as the following polynomial, using the required values
$z_1$ through $z_4$ as coefficients:
\begin{eqnarray}
  p(x,y) & := & 
  z_{1} xy \;+\;
  z_{2} x (1-y) \;+\; 
  z_{3} (1-x) y \;+\; 
  z_{4} (1-x)(1-y)
  \label{eq:dev-boole}
\end{eqnarray}
Following Boole's special law we add the constraints $x^2=x$ and
$y^2=y$ to require that each variable must be either $0$ or $1$.  Note
that for any $x\in\{0,1\}$ and $y\in\{0,1\}$ only one of the four
terms Equation~\ref{eq:dev-boole} attains a nonzero value.  This
feature, which holds for functions of any number of binary variables,
is the basis of Boole's polynomial translation method.

For example, the translation of the logical exclusive disjunction $x
\oplus y$ uses its truth table:
\begin{equation}
  \begin{array}{cc|c}
    x & y & x \oplus y \\ \hline
    1 & 1 & 0 \\
    1 & 0 & 1 \\
    0 & 1 & 1 \\
    0 & 0 & 0 \\
  \end{array}
  \label{eq:tbl-or}
\end{equation}
which produces the following polynomial according to
Equation~\ref{eq:dev-boole}:
\begin{eqnarray}
  p(x,y) & := & x (1-y) \;+\; (1-x) y \\
  & \Rightarrow & x + y - 2xy
\end{eqnarray}
along with the constraints $x^2=x$ and $y^2=y$.  It can be verified by
simple calculations that $p(x,y)$ satisfies the stated criteria:
$p(1,1)=0$; $p(1,0)=1$; $p(0,1)=1$; and $p(0,0)=0$.  Note that in
Boole's original method there is no cause to invoke unusual rules of
arithmetic such as $1+1=1$.  Boole's arithmetic was emphatically the
standard fare, not what we now call `Boolean algebra'.  However, as we
shall soon see it is helpful to use instead of real numbers
coefficients in the finite field $\mathbb{F}_2$ and thus adopt integer
arithmetic modulo $2$ (in which $1+1=0$).

Table~\ref{tbl:translation} shows the polynomial representations of
several logical formulas using Boole's original method, alongside
their conventional forms.  In general, for any formula of the
propositional calculus in a logical system with atomic formulas
$\mathbf{x}:=(x_1,x_2,\ldots,x_n)$, Boole's representation method
yields a polynomial in the ring $\mathbb{R}[\mathbf{x}]$ with
real-number coefficients.  The coefficients are more specifically
always integers; hence the Boolean polynomials can also be considered
members of the ring $\mathbb{Z}[\mathbf{x}]$ with integer coefficients
or the ring $\mathbb{Q}[\mathbf{x}]$ with rational coefficients as
convenient.

Note that there are some polynomials in the ring
$\mathbb{Z}[\mathbf{x}]$ (hence also $\mathbb{Q}[\mathbf{x}]$ and
$\mathbb{R}[\mathbf{x}]$) that do not correspond to any well-formed
logical formulas at all using Boole's representation; Boole recognized
this fact and called such polynomials `not interpretable.'
Specifically, in Boole's original scheme it is not the case that
polynomial addition corresponds directly to either kind of logical
disjunction (exclusive or nonexclusive): the polynomial $x+y$ is
neither the translation of $x \oplus y$ (using \textsc{xor}) nor the
translation of $x \vee y$ (using \textsc{or}).  The polynomial
difference $x-y$ is similarly not interpretable as a logical formula
in Boole's original translation scheme.

There were some minor flaws in Boole's presentation of his polynomial
representation method.  Although the general algorithm presented in
Chapter~V of his \emph{Laws of Thought} is correct, Boole did not
always apply his own algorithm correctly.  For example on page~95 of
\cite{boole} he translated $y\wedge(z \oplus w)$ as $yz+yw$, having
omitted the last term $-2yzw$.  It is also confusing that Boole used
the signs for addition and subtraction in two different senses
(sometimes as the elementary arithmetic operators and sometimes to
signify set union and set difference); either usage is fine but it
takes careful reading to disambiguate the overloading.  Finally, Boole
did not have very robust methods to solve the multivariate polynomial
equations that he had formulated.

\subsubsection{Revised Translation Using Linear Algebra and Finite
  Fields} 
\label{sec:field}

Instead of using an ordinary number system such as the real numbers
$\mathbb{R}$, logical formulas can be translated into polynomials with
coefficients and variable values in a finite field $\mathbb{F}_d$ of
(prime) order $d$.  At the same time the range of acceptable input
formulas can be widened from those using binary logic to those using
multivalued logic (with $d \ge 2$ elementary truth values).  There are
three main benefits to using coefficients in a finite field
$\mathbb{F}_d$ instead of the real numbers $\mathbb{R}$ (or the
integers, rationals, or complex numbers).  First, the polynomial ring
$\mathbb{F}_d[\mathbf{x}]$ over a finite set $\mathbf{x} :=
(x_1,x_2,\ldots,x_n)$ of variables is itself countable and finite.
Having the set of possible polynomials thus limited allows a simple
tabular approach to solving systems of equations.  In contrast there
are infinitely many polynomials in any ring $\mathbb{R}[\mathbf{x}]$
(as well as $\mathbb{Z}[\mathbf{x}]$ etc.).  Second, polynomial
translations of logical formulas are a bit simpler and more intuitive
using finite-field coefficients; for example the operation of addition
in a polynomial ring $\mathbb{F}_2[\mathbf{x}]$ maps directly to
exclusive disjunction in $2$-valued logic.  Third, every polynomial in
a ring $\mathbb{F}_d[\mathbf{x}]$ corresponds to a well-formed formula
in $d$-valued logic; there are no longer any uninterpretable
polynomials.

First let us generalize the notion of a logical truth table into a
\emph{finite-integer function}.

\begin{definition}[Finite-Integer Function]
  \label{def:finite-function}
  A finite-integer function $T:(\mathbb{Z}_d)^n \rightarrow
  \mathbb{Z}_d$ of order $d$ and arity $n$ maps from the
  $n$-dimensional affine space $(\mathbb{Z}_d)^n$ to the set
  $\mathbb{Z}_d$, where $\mathbb{Z}_d := \{0,1,2,\ldots,d-1\}$ is
  the ring of integers modulo $d$ (with $d \ge 2$).  Such a function
  $T$ can be visualized as a table in which each row $i$ describes the
  corresponding mapping $(a_{i,1},a_{i,2},\ldots,a_{i,n}) \mapsto
  z_i$, such that $T(x_1,x_2,\ldots,x_n)=z_i$ when each $x_j=a_{i,j}$.
  \begin{equation}
    \begin{array}{c|cccc|c} \hline
      i & x_1 & x_2 & \cdots & x_n & T(\mathbf{x}_i) \\ \hline
      1 & a_{1,1} & a_{1,2} & \cdots & a_{1,n} & z_1 \\
      2 & a_{2,1} & a_{2,2} & \cdots & a_{2,n} & z_2 \\
      \vdots & \multicolumn{4}{c|}{\ddots} & \vdots \\
      d^n & a_{d^n,1} & a_{d^n,2} & \cdots & a_{d^n,n} & z_{d^n} \\ \hline
    \end{array}
  \end{equation}
\end{definition}

\begin{lem}[Polynomial Encoding by Linear Algebra]
  \label{lem:matrix}
  Consider a finite-integer function $T:(\mathbb{Z}_d)^n \rightarrow
  \mathbb{Z}_d$ of arity $n$ and order $d$ as described in
  Definition~\ref{def:finite-function}, with the additional
  restriction that $d$ is prime.  This function $T$ can be represented
  in closed form as a polynomial $p$ over the variables $\mathbf{x} :=
  (x_1,x_2,\ldots,x_n)$, with coefficients that are either rational
  numbers or integers in the modular ring $\mathbb{Z}_d$, such that
  $p(\mathbf{x})=T(\mathbf{x})$ for every point
  $\mathbf{x}\in(\mathbb{Z}_d)^n$ (using modular arithmetic to
  evaluate $p$ if its coefficients are in $\mathbb{Z}_d$, otherwise
  ordinary arithmetic). The developed polynomial $p$ requires $d^n$
  coefficients, which can be computed as the solution of a system of
  $d^n$ linear equations in $d^n$ variables.

  Referring to Definition~\ref{def:finite-function}, let us designate
  a point $\mathbf{a}_i := (a_{i,1},a_{i,2},\ldots,a_{i,n}) \in
  (\mathbb{F}_d)^n$ as an \emph{index vector}; we consider the list
  $(\mathbf{a}_1,\mathbf{a}_2,\ldots,\mathbf{a}_{d^n})$ of all $d^n$
  possible index vectors to be arranged in some (arbitrary) order
  which is used throughout this computation.  By Lemma~\ref{lem:count}
  there are $d^n$ possible monomials in the polynomial ring
  $\mathbb{F}_d[x_1,x_2,\ldots,x_n]$; each of them can be generated by
  using an index vector to supply exponents for the variables in the
  tuple $\mathbf{x} := (x_1,x_2,\ldots,x_n)$.  From each index
  vector $\mathbf{a}_j$ we generate the corresponding monomial $t_j$
  as the following product:
  \begin{eqnarray*}
    t_j & := & (x_1)^{a_{j,1}} (x_2)^{a_{j,2}} \cdots (x_n)^{a_{j,n}}
  \end{eqnarray*}
  The vector $\mathbf{t} := (t_1,t_2,\ldots,t_{d^n})$ contains all
  the monomials for the polynomial ring.  An index vector
  $\mathbf{a}_i$ can also be used as coordinates for some point
  $\mathbf{x}_i \in (\mathbb{Z}_d)^n$, as given by:
  \begin{eqnarray*}
    \mathbf{x}_i & := & (a_{i,1},a_{i,2},\ldots,a_{i,n})
  \end{eqnarray*}
  Combining these uses, we construct a matrix $\mathbf{M}$ of size
  $d^n \times d^n$ in which each row $\mathbf{M}_i$ is the monomial
  vector $\mathbf{t}$ evaluated at the point
  $\mathbf{x}_i=\mathbf{a}_i$ defined by the index vector
  $\mathbf{a}_i$:
  \begin{eqnarray*}
    \mathbf{M}_i & := & \mathbf{t}(\mathbf{a}_i)
  \end{eqnarray*}
  Equivalently the element at each row $i$ and each column $j$ of the
  matrix $\mathbf{M}$ is the value of the monomial $t_j$ defined by the
  exponents $\mathbf{a}_j$ evaluated at the point $\mathbf{x}_i$
  defined by the coordinates $\mathbf{a}_i$ (taking $0^0=1$ for this
  calculation):
  \begin{eqnarray*}
    \mathbf{M}_{i,j} & := & 
    (a_{i,1})^{a_{j,1}} 
    (a_{i,2})^{a_{j,2}} \cdots 
    (a_{i,n})^{a_{j,n}}
  \end{eqnarray*}
  Each matrix element is roughly $(\mathbf{a}_i)^{\mathbf{a}_j}$,
  where exponentiation works element-wise and the resulting products
  are themselves multiplied together.  

  Next we introduce the column vector $\mathbf{z} :=
  (z_1,z_2,\ldots,z_{d^n})$ in which each $z_i := T(\mathbf{x}_i)$ is
  the value of the finite-integer function $T$ evaluated at the
  corresponding point.  Finally we use the column vector $\mathbf{c}
  := (c_1,c_2,\ldots,c_{d^n})$ for the coefficients of the polynomial
  $p$ which are to be determined.  The system of linear equations
  $\mathbf{M} \mathbf{c} = \mathbf{z}$ is a restatement of the
  original function $T$, amended to include the computed monomial
  values.  The desired coefficients $\mathbf{c}$ of the new polynomial
  $p$ are given by the solution $\mathbf{c} = \mathbf{M}^{-1}
  \mathbf{z}$ to this linear system.  The resulting polynomial $p =
  \mathbf{c} \cdot \mathbf{t} =
  \left(\mathbf{M}^{-1}\mathbf{z}\right)\cdot \mathbf{t}$ is the
  scalar product of the solved coefficient vector and the monomials.
  The matrix inversion and scalar product can be computed using
  ordinary arithmetic, in which case the inferred polynomial $p$ will
  be a member of the ring $\mathbb{Q}[\mathbf{x}]$ with rational
  coefficients (thus also a member of $\mathbb{R}[\mathbf{x}]$ with
  real coefficients); alternatively, the requisite calculations can be
  performed using integer arithmetic modulo $d$ in which case the
  inferred polynomial $p$ will be a member of the polynomial ring
  $\mathbb{Z}_d[\mathbf{x}]$ with coefficients in the ring of integers
  modulo $d$ (which is the finite field $\mathbb{F}_d$ when $d$ is
  prime).
\end{lem}
Lemma~\ref{lem:matrix} assumes that the matrix $\mathbf{M}$
constructed by its directions is invertible.

\begin{cor}[Encoding Indeterminate Polynomials]
  \label{cor:indeterminate}
  The polynomial construction in Lemma~\ref{lem:matrix} works
  essentially unchanged if the values in the finite-integer function
  $T$ are left indeterminate.  Using the elements of the vector
  $\mathbf{z} := (z_1,z_2,\ldots,z_{d^n})$ of function values as
  symbolic variables, the result $p =
  \left(\mathbf{M}^{-1}\mathbf{z}\right)\cdot \mathbf{t}$ inferred by
  Lemma~\ref{lem:matrix} is a member of the extended polynomial ring
  $\mathbb{F}_d[x_1,x_2,\ldots,x_n;z_1,z_2,\ldots,z_{d^n}]$ that
  includes these $d^n$ new variables.  This allows a fully-parametric
  closed-form representation of any finite-integer function $T$ of
  order $d$ and arity $n$ as described in
  Definition~\ref{def:finite-function}, at the cost of introducing
  $d^n$ new variables.
\end{cor}

\begin{example}\exfmt
  Using $2$-valued logic (thus $d=2$) and the variables $x$ and $y$
  (thus $\mathbf{x}=(x,y)$ and $n=2$) leads to the following $d^n=4$
  index vectors according to Lemma~\ref{lem:matrix}:
  \begin{equation}
    \mathbf{a}_1 = (1,1), \quad 
    \mathbf{a}_2 = (1,0), \quad 
    \mathbf{a}_3 = (0,1), \quad 
    \mathbf{a}_4 = (0,0)
    \label{eq:vec-a}
  \end{equation}
  Using these index vectors as exponents generates the following
  monomials:
  \begin{equation}
    t_1 = x^1 y^1 = xy, \quad 
    t_2 = x^1 y^0 = x, \quad 
    t_3 = x^0 y^1 = y, \quad 
    t_4 = x^0 y^0 = 1
  \end{equation}
  and thus the monomial vector $\mathbf{t}=(xy,x,y,1)$.  The following
  \emph{Mathematica} commands compute these index vectors and
  monomials:
\begin{verbatim}
d = 2; xs = {x, y}; n = Length[xs];                (* variables *)
as = Tuples[Reverse[Range[0, d - 1]], n];          (* index vectors *)
t = Table[Apply[Times, xs^as[[i]]], {i, 1, d^n}];  (* monomials *)
\end{verbatim}
Continuing on, we construct the $4 \times 4$ matrix $\mathbf{M}$ according to
Lemma~\ref{lem:matrix}:
\begin{eqnarray}
  \mathbf{M} & = & \left(
  \begin{array}{cccc}
    1 & 1 & 1 & 1 \\
    0 & 1 & 0 & 1 \\
    0 & 0 & 1 & 1 \\
    0 & 0 & 0 & 1 \\
  \end{array}
  \right)
\end{eqnarray}
For example the second row $\mathbf{M}_2=(0,1,0,1)$ is the value of
the monomials $\mathbf{t}=(xy,x,y,1)$ at the point $(x,y)=(1,0)$
defined by the second index vector $\mathbf{a}_2$.  In
\emph{Mathematica} the matrix \verb|M| is computed by:
\begin{verbatim}
M = Table[ReplaceAll[t, Table[xs[[j]] -> as[[i]][[j]], {j, 1, n}]], 
       {i, 1, d^n}];
\end{verbatim}

Now we add a vector $\mathbf{z}=(0,1,1,0)$ of truth values for the
logical formula $x \oplus y$ arranged in the same order as the index
vectors in Equation~\ref{eq:vec-a}.  These correspond to the truth
table shown in Equation~\ref{eq:tbl-or}.  The desired coefficients in
the inferred polynomial are given by the linear system
$\mathbf{M}\mathbf{c}=\mathbf{z}$.  Using ordinary arithmetic, this
system yields the solution coefficients $\mathbf{c}=(-2,1,1,0)$.
Combining this solution $\mathbf{c}$ with the generated monomials
$\mathbf{t}$, the inferred polynomial $p$ in the ring
$\mathbb{R}[x,y]$ is given by $\mathbf{c} \cdot \mathbf{t}$; thus
$p(x,y):x+y-2xy$.  In \emph{Mathematica} this same solution is
accomplished by:
\begin{verbatim}
z = {0, 1, 1, 0};  
c = LinearSolve[M, z];  
p = c . t;
\end{verbatim}
For completeness we add the constraints $x^2=x$ and $y^2=y$ to limit
the values of these variables appropriately.  This first
linear-algebra solution agrees with the result from Boole's original
translation method.

Alternatively, solving $\mathbf{M}\mathbf{c}=\mathbf{z}$ using integer
arithmetic modulo $2$ yields the solution vector
$\mathbf{c}=(0,1,1,0)$ and consequently the inferred polynomial
$p(x,y):x+y$ in the ring $\mathbb{F}_2[x,y]$.  No additional
constraints are necessary in this case.  In \emph{Mathematica} this
same modular-arithmetic solution is computed by:
\begin{verbatim}
c2 = LinearSolve[M, z, Modulus -> d];
p2 = PolynomialMod[c2 . t, d];
\end{verbatim}

Leaving the desired values $(z_1,z_2,z_3,z_4)$ indeterminate and
computing $(\mathbf{M}^{-1} \mathbf{z}) \cdot \mathbf{t}$ produces the
following parametric polynomial $p_1$ in the ring
$\mathbb{R}[x,y,z_1,z_2,z_3,z_4]$:
  \begin{eqnarray}
    p_1 & : &
    x y z_1-x y z_2-x y z_3+x y z_4-x z_1+x z_3-y z_1+y z_2+z_1
  \end{eqnarray}
  In factored form this is exactly Boole's polynomial shown in
  Equation~\ref{eq:dev-boole}.  Using modular arithmetic instead
  yields a similar parametric polynomial $p_2$ in the ring
  $\mathbb{F}_2[x,y,z_1,z_2,z_3,z_4]$:
  \begin{eqnarray}
    p_2 & : &
    x y z_1+x y z_2+x y z_3+x y z_4+x z_1+x z_3+y z_1+y z_2+z_1
  \end{eqnarray}
  In \emph{Mathematica} these parametric solutions are produced by:
\begin{verbatim}
Clear[z];  (* remove old numeric values *)
zs = Table[Subscript[z, i], {i, 1, d^n}];  (* list of variables *)
p1 = Expand[(Inverse[M] . zs) . t];
p2 = PolynomialMod[(Inverse[M, Modulus -> d] . zs) . t, d];
\end{verbatim}

\end{example}

Using Lemma~\ref{lem:matrix} and modular arithmetic, the familiar
values and operations of $2$-valued logic can be rewritten as
polynomial formulas with coefficients in the binary finite field
$\mathbb{F}_2$; the results are summarized in
Table~\ref{tbl:translation}.  The basic values true (\val{t}) and
false (\val{f}) map to the respective elementary values $1$ and $0$.
Other logical functions of atomic formulas $x$ and $y$ are translated
as follows.  Logical conjunction (\textsc{and}) translates directly as
multiplication: for $x \wedge y$ we substitute the product $x \times
y$ (abbreviated $xy$ or $x \cdot y$ in the usual algebraic fashion).
Exclusive logical disjunction (\textsc{xor}) translates directly as
addition: for $x \oplus y$ we substitute the sum $x+y$.  Logical
negation translates as incrementation by $1$ (or equivalently the
difference from $1$): for the negation $\neg x$ of any formula $x$ we
substitute $1+x$ (which is the same as $x-1$ or $x+1$ or $1-x$ using
integer arithmetic modulo $2$; but not the same as $-x$ with a unary
minus, which is just $x$ itself).  Nonexclusive disjunction (the usual
\textsc{or}) is a polynomial sum: for $x \vee y$ we substitute
$x+y+xy$.  For the material implication $x \rightarrow y$ we
substitute $1+x+xy$.  The biconditional $x \leftrightarrow y$ (also
written $x \equiv y$ or designated \textsc{xnor}) is translated into
the polynomial $1+x+y$.  Finally, the nonconjunction $x \uparrow y$
(alternative denial, \textsc{nand}, Sheffer stroke, $x|y$) is
translated into the polynomial $1+xy$, and the nondisjunction $x
\downarrow y$ (joint denial, \textsc{nor}, Pierce arrow, Quine dagger)
is translated into the polynomial $1+x+y+xy$.  

There is already a \verb|BooleanTable| command in \emph{Mathematica}
to generate truth tables for binary-logical functions, and a
\verb|Boole| function to convert the logical values \verb|True| and
\verb|False| to the numbers $0$ and $1$.  The \verb|Tuples[]|
statement above was written so that the index vectors are generated in
the same order that the \verb|BooleanTable| function uses (the
variables in \verb|xs| also need to be in alphabetical order, as the
\verb|Sort| command would produce).  With these caveats it is simple
to convert logical expressions into polynomials within
\emph{Mathematica}.  For example the following commands:
\begin{verbatim}
z = Boole[BooleanTable[Implies[x,y]]];  (* truth table {1,0,1,1} *)
p1 = (Inverse[M] . z) . t;
p2 = PolynomialMod[(Inverse[M, Modulus -> d] . z) . t, d];
\end{verbatim}
produce the values $1-x+xy$ for \verb|p1| and $1+x+xy$ for \verb|p2|
which you can see listed in Table~\ref{tbl:translation} as the
translations for $x \rightarrow y$ using real and finite-field
coefficients.

\subsubsection{Translation as a Function}

The translation from a logical formula (in the propositional calculus)
to a polynomial can be considered a function; every logical formula
maps to a unique polynomial.  But in general there are infinitely many
distinct well-formed logical formulas that map to any given polynomial
in a ring $\mathbb{F}_d[\mathbf{x}]$ with finite-field coefficients;
in other words the translation function is surjective but not
injective.  By comparison, translation into $\mathbb{Z}[\mathbf{x}]$
(hence also $\mathbb{Q}[\mathbf{x}]$ or $\mathbb{R}[\mathbf{x}]$) is
neither injective nor surjective---some polynomials with integer
(hence also rational or real) coefficients correspond to many
well-formed logical formulas and others correspond to none.

\begin{definition}[Logical Preimage of a Polynomial Function]
Let us designate the set of well-formed logical formulas that map to a
given polynomial as the \emph{logical preimage} of that polynomial;
these formulas constitute an equivalence class in the sense that they
always have the same value as one another (their truth tables are
identical).
\end{definition}
For example the logical formulas $x \rightarrow y$ and $\neg x \vee y$
and $(x \rightarrow y)\wedge(\neg x \vee y)$ are each translated into
the polynomial $1+x+xy$ in the ring $\mathbb{F}_2[x,y]$; therefore
they are members of its logical preimage.

\subsection{Equations from Axioms}
\label{sec:axioms-equations}

An axiom in a logical system has the same meaning as an equation in an
algebraic system: each is an assertion made for the purpose of
deducing what follows.  Equations and axioms are \emph{constraints}
rather than \emph{commandments}: it is generally possible to write
systems of axioms or other equations that cannot be satisfied.  The
essential nature of an equation does not change because its formulas
originated in logical instead of polynomial notation.

Logical axioms translate directly as polynomial equations, when those
axioms are simple formulas from the propositional calculus (not
involving references to provability, quantifiers, or indeterminate
predicates over infinite domains; any of which complicate matters as
addressed later).  Using the turnstile $\vdash$ to mark an axiom as in
Frege's \emph{Begriffsschrift} \cite{frege}, an axiom $\vdash q$ is
the assertion that the included formula $q$ is true; in Frege's
terminology the whole axiom is a `judgment' and the included formula
is its `content.'  Let us say that the logical system of interest has
a set $A := \{\vdash q_1,\vdash q_2,\ldots,\vdash q_m\}$ of axioms, in
which each content $q_i$ is a formula built from some atomic formulas
$\mathbf{x} := (x_1,x_2,\ldots,x_n)$ and the usual logical operations.
For concreteness we assume $2$-valued logic, though this is not
required for algebraic analysis.

Using Table~\ref{tbl:translation} or the general method in
Lemma~\ref{lem:matrix} each axiom $\vdash q_i$ can be translated into
a polynomial equation with either real-number coefficients or binary
finite-field coefficients.  Each content $q_i$ is first translated
from logical to polynomial notation using $\mathbb{R}[\mathbf{x}]$ or
$\mathbb{F}_2[\mathbf{x}]$ as desired.  Then each judgment is
translated to a polynomial constraint in the obvious way: the axiom
$\vdash q_i$ becomes the equation $q_i=1$.  The set $A$ of translated
axiom-equations has a conjunction polynomial according to
Lemma~\ref{lem:conjunction}:
\begin{eqnarray}
  \label{eq:qstar-raw}
  q^* & := & (q_1)(q_2) \cdots (q_m) - 1
\end{eqnarray}
with a slightly different form since each right-hand-side is $1$
instead of $0$.  The constraint $q^*=0$ using this conjunction
polynomial corresponds to the assertion that the logical conjunction
$q_1 \wedge q_2 \wedge \cdots \wedge q_m$ of the axiom contents is
true.  Therefore using binary finite-field coefficients a set $A$ of
axioms translates to a single equation featuring its conjunction
polynomial:
\begin{eqnarray}
  A & \leadsto & \left\{ q^* = 0 \right\}, \quad q^* \in
  \mathbb{F}_2[\mathbf{x}] 
\end{eqnarray}

\begin{example}\exfmt
\label{ex:carroll}

Lewis Carroll \cite{carroll-barbers} provided this succinct
version of his delightful logic puzzle about the barbers Allen, Brown,
and Carr:
\begin{em}
\begin{itemize}
\item[] There are two Propositions, $A$ and $B$.
\item[] It is given that
  \begin{enumerate}
  \item[(1)] If $C$ is true, then, if $A$ is true, $B$ is not true;
  \item[(2)] If $A$ is true, $B$ is true.
  \end{enumerate}
\item[] The question is, can $C$ be true?
\end{itemize}
\end{em}
Using the variables $a$, $b$, and $c$ for the respective propositions
$A$ (that Allen is out of the shop), $B$ (that Brown is out), and $C$
(that Carr is out), we model the problem as the following axioms:
\begin{eqnarray}
& & \vdash \; c \rightarrow (a \rightarrow \neg b) \label{eq:ax-abc-1}\\
& & \vdash \; a \rightarrow b \label{eq:ax-abc-2}
\end{eqnarray}
Translating the content of each axiom into a polynomial in the ring
$\mathbb{F}_2[a,b,c]$ with binary finite-field coefficients according
to Table~\ref{tbl:translation} and Lemma~\ref{lem:matrix} yields
the following polynomial equations:
\begin{eqnarray}
  \vdash \; c \rightarrow (a \rightarrow \neg b) & \leadsto &  
  abc+1=1 \label{eq:ax-abc-p1} \\ 
  \vdash \; a \rightarrow b & \leadsto & ab+a+1=1 \label{eq:ax-abc-p2}
\end{eqnarray}
By Lemma~\ref{lem:conjunction} the conjunction polynomial for this
pair of equations is given by the following expression (calculated by
integer arithmetic modulo $2$ since the polynomials use coefficients
in $\mathbb{F}_2$):
\begin{equation}
  \label{eq:abc-conjunction}
  q^*  \quad := \quad (abc+1)(ab+a+1) - 1 \quad \Rightarrow \quad abc+ab+a
\end{equation}
Thus the axioms $A$ in Equations \ref{eq:ax-abc-1} and
\ref{eq:ax-abc-2} are translated into to the following set of a
solitary polynomial equation:
\begin{eqnarray}
  A & \leadsto & \{abc+ab+a = 0\}
\end{eqnarray}
using polynomial coefficients in $\mathbb{F}_2$.
\end{example}

\subsection{Solution-Value Sets as Truth Values}
\label{sec:truth-values}

Translating axioms from logical to polynomial notation is only the
first step; next it is necessary to solve these equations and
interpret their solution.  The \emph{truth value} of a logical formula
$p$ subject to a set $A$ of axioms is given by the solution-value set
$\set{S}_A(p)$ of that formula subject to those axiom-equations.  In
order to compute the solution-value set as specified in
Definition~\ref{def:objective}, it is first necessary to compute the
solution set $\set{V}(A)$ to the equations $A$ that is specified in
Definition~\ref{def:polynomial-system}, using either a general method
from computational algebraic geometry or the manual table-based method
given in Definition~\ref{def:inference}.  Calculations can be
performed with real or finite-field coefficients.

According to Definition~\ref{def:objective} the solution-value set
$\set{S}_A(p) \subseteq K$ is a subset of the set of elementary values
in the algebraic field $K$ that contains the value of the objective
formula $p$; therefore the possible values of $\set{S}_A(p)$ are the
members of the power set $2^K$.  Thus using the binary finite field
$\mathbb{F}_2$ with elementary values $\{0,1\}$ there are four
possible solution-value sets:
\begin{equation}
  \{ \{\}, \{0\}, \{1\}, \{0,1\} \}
  \label{eq:powerset}
\end{equation}
It happens that these are the only four solution-value sets
encountered even when logical formulas have been translated into the
polynomial ring $\mathbb{R}[\mathbf{x}]$ with real coefficients, since
the range of each polynomial function is still limited to $\{0,1\}$ by
its construction in Lemma~\ref{lem:matrix}.  Anyway each set of
elementary values in Equation~\ref{eq:powerset} makes sense as the
truth value of a logical formula.  As described in
Definition~\ref{def:objective}, these solution-value sets allow
logical formulas to be categorized into those that are necessarily
true ($\set{S}_A(p)=\{1\}$), necessarily false ($\set{S}_A(p)=\{0\}$),
ambiguous ($\set{S}_A(p)=\{0,1\}$), and unsatisfiable
($\set{S}_A(p)=\{\}$).

Note that an objective formula is unsatisfiable if and only if the
whole set of axiom-equations is inconsistent; in that case every
possible formula is unsatisfiable subject to those same equations.
For example given the unsatisfiable constraint $0=1$ there is no
feasible value for any objective formula: $\set{S}_{\{0=1\}}(0)
\Rightarrow \{\}$, $\set{S}_{\{0=1\}}(1) \Rightarrow \{\}$,
$\set{S}_{\{0=1\}}(x) \Rightarrow \{\}$, and so on.  In contrast to
the principle of explosion that an inconsistent set of axioms allows
\emph{any} formula to be proven true (\emph{ex contradictione sequitur
  quodlibet}), in the algebraic paradigm an inconsistent set of axioms
allows \emph{no} formulas to be proven true (nor false, nor
ambiguous).  Moreover, such undecidability by reason of inconsistency
is distinguished from the undecidability by reason of ambiguity that
results from trivial, tautological, or otherwise inconclusive axioms.
For example $\set{S}_{\{0=0\}}(x) \Rightarrow \{0,1\}$ differs from
$\set{S}_{\{0=1\}}(x) \Rightarrow \{\}$ since each is undecidable by a
different mechanism.

\begin{example}\exfmt
Returning to Carroll's barbershop, it was established in
Example~\ref{ex:carroll} that axioms about Allen, Brown, and Carr
translate as the equation $abc+ab+a=0$ using polynomials in the binary
finite field $\mathbb{F}_2$.  The solution set
$\set{V}(abc+ab+a=0)$ is easily calculated using the appropriate
modular arithmetic in \emph{Mathematica}:
\begin{verbatim}
In[4]:= Solve[{a b c + a b + a == 0}, {a, b, c}, Modulus -> 2]

Out[4]= {{a -> 0, b -> 0, c -> 0}, {a -> 0, b -> 0, c -> 1}, {a -> 0, 
  b -> 1, c -> 0}, {a -> 0, b -> 1, c -> 1}, {a -> 1, b -> 1, c -> 0}}
\end{verbatim}
The result contains five solutions for $(a,b,c)$:
\begin{eqnarray}
  \label{eq:abc-sol}
  \set{V}(abc+ab+a=0) & \Rightarrow & 
  \{ (0,0,0), (0,0,1), (0,1,0), (0,1,1), (1,1,0) \}
\end{eqnarray}
Using the constraints in Equations \ref{eq:ax-abc-p1} and
\ref{eq:ax-abc-p2} directly (instead of their conjunction polynomial)
gives the same solution set:
\begin{verbatim}
In[5]:= Solve[{a b c + 1 == 1, a b + a + 1 == 1}, {a, b, c}, 
 Modulus -> 2]

Out[5]= {{a -> 0, b -> 0, c -> 0}, {a -> 0, b -> 0, c -> 1}, {a -> 0, 
  b -> 1, c -> 0}, {a -> 0, b -> 1, c -> 1}, {a -> 1, b -> 1, c -> 0}}
\end{verbatim}
You can verify that the logical back-translation of each solution
satisfies the axioms in Carroll's problem stated in Equations
\ref{eq:ax-abc-1} and \ref{eq:ax-abc-2}.  For example with
$a=\val{t}$, $b=\val{t}$, and $c=\val{f}$ as in the last solution
$(a,b,c)=(1,1,0)$: the material implication $a \rightarrow b$ is true;
the inner material implication $a \rightarrow \neg b$ is false; but
since $c$ is false the outer material implication $c \rightarrow (a
\rightarrow \neg b)$ is nonetheless true.

Carroll's problem statement requests the truth value of the logical
formula $c$.  It is evident from the solution set
$\set{V}(abc+ab+a=0)$ shown in Equation~\ref{eq:abc-sol} that both
$0$ and $1$ are feasible solutions for $c$.  Thus we have:
\begin{eqnarray}
  \set{S}_{\{abc+ab+a=0\}}(c) & \Rightarrow & \{0,1\}
\end{eqnarray}
and it is neither the case that $c$ is necessarily $0$ nor the case
that $c$ is necessarily $1$.  The formula $c$ is
\emph{ambiguous}---given Carroll's axioms the barber Carr could be
either in or out of the shop.

Let us also calculate the truth value of the proposition $abc$.  You
can see by inspection of Equation~\ref{eq:abc-sol} that every solution
in $\set{V}(A)$ contains the value zero for at least one variable
in $(a,b,c)$.  Thus the product $abc$ must always be zero using these
solutions; in other words the solution-value set for the objective
formula $abc$ is given by:
\begin{eqnarray}
  \set{S}_{\{abc+ab+a=0\}}(abc) & \Rightarrow & \{0\}
\end{eqnarray}
This result indicates that the polynomial $abc$ is necessarily $0$,
and therefore that the corresponding logical formula is necessarily
false---given Carroll's axioms all three barbers cannot be out of the
shop simultaneously.

\end{example}

\subsection{Discovering All Theorems}
\label{sec:theorem}

In logic there is a special name for a formula that is necessarily
true: it is a \emph{theorem}.  Therefore a logical formula $p$ is a
theorem given some set $A$ of axioms exactly if its solution-value set
$\set{S}_A(p)$ has the value $\{1\}$ (that is, the set containing
exactly the number one).  Furthermore the set of all theorems entailed
by the set $A$ of axioms (using some finite list
$\mathbf{x}:=(x_1,x_2,\ldots,x_n)$ of propositional variables) is
given by the inverse-value set $\set{S}_A^{-1}(\{1\})$ described in
Definition~\ref{def:inverse}, Lemma~\ref{lem:singleton}, and
Corollary~\ref{cor:increment}.  Every logical formula (whose
propositional variables are $\mathbf{x}$) that is a theorem given
these axioms must translate as one of the polynomials in this
inverse-value set (which is empty if the axioms are unsatisfiable).
Using polynomials with finite-field coefficients, this algebraic
formulation transforms the discovery of theorems from a search through
an infinite set of logical formulas to a search among a finite set of
polynomials (each of which corresponds to a unique truth table).

It is left as a separate exercise to choose which member of the
logical preimage of each polynomial theorem should be used to
represent it; this is properly an optimization problem in the area of
\emph{logic synthesis}.  For a polynomial with coefficients in the
binary field $\mathbb{F}_2$, a reasonable default choice for its
representative logical formula is the transliteration of the
polynomial that maps multiplication back to logical conjunction
($\wedge$, \textsc{and}) and that maps addition back to logical
exclusive disjunction ($\oplus$, \textsc{xor}).  But note that it is
not appropriate to use what is commonly called `disjunctive normal
form' since as Table~\ref{tbl:translation} indicates, non-exclusive
disjunction ($\vee$, \textsc{or}) does not map directly to polynomial
addition (using either real or finite-field coefficients).  Note also
that an exclusive disjunction $x_1 \oplus x_2 \oplus \cdots \oplus
x_n$ in $2$-valued logic is true if and only if an odd number of its
arguments are true; otherwise the disjunction is false.  This
corresponds to integer arithmetic modulo $2$ which is appropriate for
polynomials with coefficients in the binary finite field
$\mathbb{F}_2$.

The concept underlying theorems is not exclusive to logic; the more
primitive notion is that a formula might have a definite solution
value given some set of equations or other constraints.  In general
algebra the definite solution value zero is given special status.  As
mentioned after Definition~\ref{def:inverse}, what is defined herein
as the inverse-value set $\set{S}_A^{-1}(\{0\})$ is in algebraic
geometry called the \emph{ideal} generated by the polynomials in the
equations $A$.  In this sense a theorem is almost ideal!  Others have
made the connection between logical theorems and polynomial ideals
\cite{kapur-narendran,roanes-lozano}.

\begin{example}\exfmt
  Consider the axioms $\vdash x$ and $\vdash x \rightarrow y$.  Using
  Table~\ref{tbl:translation} each axiom translates as an equation
  using polynomials with binary finite-field coefficients in the ring
  $\mathbb{F}_2[x,y]$:
  \begin{eqnarray}
    \vdash x & \leadsto & x = 1 \\
    \vdash x \rightarrow y & \leadsto & 1+x+xy = 1
  \end{eqnarray}
  By Lemma~\ref{lem:conjunction} these two equations yield the
  conjunction polynomial $q^* = (x)(1+x+xy)-1$ which evaluates to
  $1+xy$ using modular arithmetic.  Thus the axioms are translated as
  the set $A=\{1+xy=0\}$ containing one polynomial equation.  Using the
  table-based inference method in Definition~\ref{def:inference} or
  the \emph{Mathematica} command:
\begin{verbatim}
Solve[{1+ x y == 0}, {x, y}, Modulus -> 2]
\end{verbatim}
  reveals that his set $A$ of equations has the solution set
  $\set{V}(A) = \{(1,1)\}$ for $(x,y)$.  Evaluating each of the $16$
  polynomials in the ring $\mathbb{F}_2[x,y]$ at this unique solution
  reveals $8$ polynomials whose only feasible value is $1$ (those
  whose listed value is $1$ in the column labeled $p_i(1,1)$ in
  Table~\ref{tbl:worksheet}, namely $p_2$, $p_3$, $p_5$, and so on).
  These $8$ polynomials, which comprise the inverse-value set
  $\set{S}_{\{xy+1=0\}}^{-1}(\{1\})$, are the translations of all the
  theorems entailed by the stated axioms:
  \begin{equation}
    1, \quad y, \quad x \quad 1+x+y, \quad 
    xy, \quad 1+y+xy, \quad 1+x+xy, \quad x+y+xy
  \end{equation}
  Following Lemma~\ref{lem:singleton} and
  Corollary~\ref{cor:increment}, this inverse-value set can also be
  described in closed form as:
  \begin{eqnarray}
    \set{S}_{\{xy+1=0\}}^{-1}(\{1\}) & \Rightarrow &
    \left\{
    p \times (1+xy) + 1 : p \in \mathbb{F}_2[x,y]
    \right\}
  \end{eqnarray}
  where modular arithmetic is to be used to evaluate the included
  polynomial expression.  Each polynomial theorem has infinitely many
  formulas in its logical preimage; for example the last polynomial
  $x+y+xy$ corresponds to the logical formula $(x \wedge y) \oplus x
  \oplus y$ as well as to the logical formula $x \vee y$.  Using the
  matching formulas in Table~\ref{tbl:translation} provides one choice
  of logical back-translation for each of the $8$ theorems entailed by
  the axioms $\vdash x$ and $\vdash x \rightarrow y$:
  \begin{equation}
    \val{t}, \quad y, \quad x, \quad x \leftrightarrow y, 
    \quad x \wedge y, \quad y \rightarrow x, \quad x \rightarrow y,
    \quad x \vee y
  \end{equation}
  Included in this complete set of theorems is the sole formula $y$
  that a direct application of \emph{modus ponens} would prove from
  the axioms $\vdash x$ and $\vdash x \rightarrow y$.  You can see
  that algebraic analysis provides more comprehensive results.
\end{example}

\begin{example}\exfmt

Returning once again to Carroll's barbershop, we can use algebraic
analysis to discover the set of all theorems entailed by the axioms
$\vdash c\rightarrow(a \rightarrow \neg b)$ and $\vdash a \rightarrow
b$.  The conjunction polynomial $q^* = a+ab+abc$ in
Equation~\ref{eq:abc-conjunction} is a translation of these logical
axioms.  It follows from Definition~\ref{def:inverse},
Lemma~\ref{lem:singleton}, and Corollary~\ref{cor:increment} that the
set of theorem polynomials is given by:
\begin{eqnarray}
  \label{eq:barber-theorems}
  \set{S}_{\{a+ab+abc=0\}}^{-1}(\{1\}) & \Rightarrow &
  \{ p \times (a+ab+abc) + 1 : p \in \mathbb{F}_2[a,b,c] \}
\end{eqnarray}
By Lemma~\ref{lem:count} there are $2^{2^3}=256$ distinct polynomials
in the ring $\mathbb{F}_2[a,b,c]$.  Evaluating all of these
polynomials with \emph{Mathematica} reveals that there are exactly
$8$ theorems in the polynomial ring $\mathbb{F}_2[a,b,c]$ entailed
by Carroll's barbershop axioms:
\begin{eqnarray}
  \set{S}_{\{a+ab+abc=0\}}^{-1}(\{1\}) & \Rightarrow & 
  \left\{
  \begin{array}{l}
    1 \\ 
    1 + a + a b + a b c \\ 
    1 + a c \\
    1 + a + a b + a c + a b c \\
    1 + a b c \\ 
    1 + a + a b \\ 
    1 + a c + a b c \\
    1 + a + a b + a c
  \end{array}
  \right\}
\end{eqnarray}
The logical negation of each of these $8$ theorems has the
solution-value set $\{0\}$; thus each negated theorem is necessarily
false.  The remaining $240$ polynomials in $\mathbb{F}_2[a,b,c]$ are
ambiguous, with the common solution-value set $\{0,1\}$.  Note that
Carroll's original query $c$ and its negation $c+1$ are both members
of this set $\set{S}_{\{a+ab+abc=0\}}^{-1}(\{0,1\})$ of ambiguous
polynomials.

For concreteness, the following \emph{Mathematica} commands were used
to generate the members of the inverse-value set
$\set{S}_{\{a+ab+abc=0\}}^{-1}(\{1\})$ shown above:
\begin{verbatim}
d = 2; xs = {a, b, c}; n = Length[xs];             (* variables *)
as = Tuples[Reverse[Range[0, d - 1]], n];          (* index vectors *)
t = Table[Apply[Times, xs^as[[i]]], {i, 1, d^n}];  (* monomials *)
all = Tuples[Range[0, d - 1], d^n] . t;            (* all polys in ring *)
(* construct inverse-value set from closed-form expression: *)
theorems = Table[all[[i]] * (a + a b + a b c) + 1, {i, 1, Length[all]}];
(* substitute to remove squares by Fermat's little theorem: *)
theorems = Expand[theorems] /. {a^2 -> a, b^2 -> b, c^2 ->c };
theorems = DeleteDuplicates[PolynomialMod[theorems, 2]]  (* modulo 2 *)
\end{verbatim}
The last command gives the output displayed above:
\begin{verbatim}
{1, 1 + a + a b + a b c, 1 + a c, 1 + a + a b + a c + a b c, 
 1 + a b c, 1 + a + a b, 1 + a c + a b c, 1 + a + a b + a c}
\end{verbatim}
\end{example}

\subsection{Dynamical Systems from References to Provability}

A logical formula that refers to its own provability can be translated
into algebraic form using a parametric system of polynomial equations
as described in Definition~\ref{def:parametric}.  The parameters and
parameter-updating functions introduced in the definition are
expressly permitted to use features of solution sets and
solution-value sets; this allows references to provability to be
modeled.  From the specified parametric system of polynomial equations
is extracted by Definition~\ref{def:dynamical} an evolution function
$F(\Theta)$ that describes how the parameters change as a function of
themselves.  This evolution function can be used in one of two ways:
either as a static constraint $\Theta=F(\Theta)$ that is added to the
system of equations (which is then instantiated from parametric form
into an ordinary system of equations; or as a recurrence relation
$\Theta_{t+1} \Leftarrow F(\Theta_t)$ that is used to extend the
static equations into a dynamical system.

\begin{example}\exfmt

Let us proceed now with the analysis of \goedel's formula \rqq{} that
asserts its own unprovability.  Using $x$ to denote the formula and
$\mathsf{Bew}(x)$ to denote the proposition that $x$ is provable
(`beweis\-bar' in German), we begin with the declaration
$x\in\{\val{t},\val{f}\}$ that the possible values of the formula $x$
are the elementary logical values true and false.  We define the
formula $x$ using the axiom:
\begin{equation}
  \vdash x = \neg \mathsf{Bew}(x)
  \label{eq:nbew}
\end{equation}
Using Definition~\ref{def:parametric} and Table~\ref{tbl:translation}
the type declaration for $x$ and the axiom in Equation~\ref{eq:nbew}
translate as the following parametric system of polynomial equations:
\begin{eqnarray}
  \begin{array}{rrcl}
    \tau_x: & x & \in & \{0,1\} \subset \mathbb{R} \\
    \delta_x : & x & \Leftarrow & 
    1 - \left( \set{S}_{\{\}}(x)=\{1\} \right)
  \end{array}
  \label{eq:nbew-templ}
\end{eqnarray}
For concreteness Equation~\ref{eq:nbew-templ} uses real-number
coefficients; however for this problem the calculations would be
identical using coefficients in $\mathbb{F}_2$ instead.

Hypothetico-deductive analysis according to
Definition~\ref{def:dynamical} now reveals the evolution function
$F(x)$ for the parametric system in Equation~\ref{eq:nbew-templ}.
Since the parameter $x$ has the domain $\{0,1\}$, we note that the
function $F$ must map from the set $\{0,1\}$ back to itself.  We
consider the two possible values of the parameter.  In the case $x=0$
the parameter-updating function is instantiated as
$\lambda_1(0):1-(\set{S}_{\{\}}(0)=\{1\})$.  Absent any constraints,
the only feasible value of the objective formula $0$ is $0$; in
general for any constant $k$ the solution-value set
$\set{S}_{\{\}}(k)$ simply yields $\{k\}$.  Hence the value of
$\lambda_1(0)$ is $1-(\{0\}=\{1\})$, which is $1$ since the false
comparison statement yields the value $0$.  Similarly, in the case
$x=1$, the parameter-updating function is instantiated as
$\lambda_1(1):1-(\set{S}_{\{\}}(1)=\{1\})$.  The solution-value set
$\set{S}_{\{\}}(1)\Rightarrow\{1\}$ gives $\lambda_1(1) \Rightarrow
1-(\{1\}=\{1\})$ which evaluates to $0$.  These mappings $0 \mapsto 1$
and $1 \mapsto 0$ define the evolution function $F(x)$ encoded by the
parametric system in Equation~\ref{eq:nbew-templ}.  It is evident from
inspection that the interpolated polynomial function $F(x):1-x$
provides the appropriate values.

There are static and dynamic ways to interpret
Equation~\ref{eq:nbew-templ}.  In the static interpretation we replace
the assignment template $\delta_x$ with the static constraint $x=F(x)$
using the extracted evolution function $F(x):1-x$.  In this static
interpretation \goedel's axiom in Equation~\ref{eq:nbew} is translated
via Equation~\ref{eq:nbew-templ} into this system of equations (the
constraint $\tau_x$ could be expressed in polynomial form as $x^2=x$
using Boole's representation scheme):
\begin{eqnarray}
  \label{eq:nbew-inst}
  \begin{array}{rrcl}
    \tau_x: & x & \in & \{0,1\} \subset \mathbb{R} \\
    \delta_x : & x & = & 1 - x \\
  \end{array}
\end{eqnarray}
It is evident from inspection that Equation~\ref{eq:nbew-inst} has no
solution; the solution-value set $\set{S}_{\{\tau_x,\delta_x\}}(x)$
evaluates to $\{\}$.  In fact simple substitution reveals that
Equation~\ref{eq:nbew-inst} is equivalent to the constraint $1=0$.
Therefore, in the static interpretation, the truth value of \goedel's
formula \rqq{} is that it is \emph{unsatisfiable}.  In other words,
there is no such thing as a logical formula $x$ that is true if and
only if $x$ is simultaneously not provable, just as there is no such
thing as a quadratic equation $2x^2+3x+c=0$ that has exactly $c$ real
solutions for $x$.  Each specification is internally inconsistent, and
by algebraic analysis each inconsistency has been exposed as an
infeasible system of equations.

The dynamic interpretation of the parametric system in
Equation~\ref{eq:nbew-templ} explains the mechanism of the
inconsistency in \goedel's formula.  In this interpretation the
parameter $x_t \in\{0,1\}$ gives the state of the dynamical system at
each time $t$; the recurrence $x_{t+1} \Leftarrow 1-x_t$ from the
evolution function $F(x):1-x$ creates a transition from the state
$x=0$ to the state $x=1$ and vice versa.  The solution-value set for
the objective $x$ computed for each time $t$ simply contains the
respective value of the state $x_t$.  Thus the dynamical system
derived from \goedel's formula is represented by the following graph:
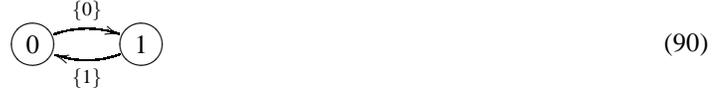
\begin{equation}
  \label{eq:bew-graph}
  \xymatrix{
    *++[o][F]{0} \ar@/^/[r]^{\{0\}} &
    *++[o][F]{1} \ar@/^/[l]^{\{1\}}
  }
\end{equation}
Following Definition~\ref{def:sequence} either initial value
$x_0$ gives an alternating infinite sequence of solution-value
sets for $x$ at successive times $t$:
\begin{eqnarray}
  \vec{\set{S}}_{\{\}}^*(x) & \Rightarrow &
  \left[ \begin{array}{rcl}
      x_0 = 0 & \mapsto & 
      (\{0\},\{1\},\{0\},\{1\},\ldots) \\
      x_0 = 1 & \mapsto & 
      (\{1\},\{0\},\{1\},\{0\},\ldots)
    \end{array} \right]
\end{eqnarray}
This dynamical system derived from Equation~\ref{eq:nbew-templ} has a
periodic orbit and no fixed points; therefore it is \emph{unsteady} in
the terminology of Definition~\ref{def:dynamical}.  In this dynamic
interpretation \goedel's self-denying formula \rqq{} specifies a
discrete dynamical system that oscillates between the states of truth
and falsity.  The objective $x$ describing the truth value of the
formula oscillates between being necessarily false (having the
solution-value set $\{0\}$) and being necessarily true (having the
solution-value set $\{1\}$).

Note that the set $A$ of constraints is empty in the dynamic
interpretation of \goedel's formula; in the parametric formulation
shown in Equation~\ref{eq:nbew-templ} there are no constraint
templates $\alpha$, only an assignment template $\delta_x$.  In this
respect \goedel's formula is more like the Fibonacci recurrence than
like the self-referential quadratic equations analyzed earlier.

\end{example}

\begin{example}
\exfmt

It is illustrative to consider the complement of \goedel's
construction: a formula that asserts its own \emph{provability}.
Using the parameter $y$ to represent the formula, we start with the
declaration $y \in \{\val{t},\val{f}\}$ and the axiom:
\begin{equation}
  \vdash y = \mathsf{Bew}(y)
  \label{eq:bew}
\end{equation}
This type declaration and axiom are translated into parametric
polynomial equations according to Definition~\ref{def:parametric}:
\begin{eqnarray}
  \begin{array}{rrcl}
    \tau_y: & y & \in & \{0,1\} \subset \mathbb{R} \\
    \delta_y : & y & \Leftarrow & 
    \left( \set{S}_{\{\}}(y)=\{1\} \right)
  \end{array}
  \label{eq:bew-templ}
\end{eqnarray}
Following Definition~\ref{def:dynamical}, the identity function
$F(y):y$ is extracted as the evolution function encoded by
Equation~\ref{eq:bew-templ}.  Incorporating the constraint $y=F(y)$
the static interpretation of Equation~\ref{eq:bew-templ} yields:
\begin{eqnarray}
  \label{eq:bew-inst}
  \begin{array}{rrcl}
    \tau_y: & y & \in & \{0,1\} \subset \mathbb{R} \\
    \delta_y : & y & = & y
  \end{array}
\end{eqnarray}
As the constraint $\delta_y$ is tautological the solution-value set
$\set{S}_{\{\tau_y,\delta_y\}}(y)$ evaluates to $\{0,1\}$.  Thus the
formula that asserts its own provability is \emph{ambiguous} in the
static interpretation.

In the dynamic interpretation the parametric system in
Equation~\ref{eq:bew-templ} specifies the following dynamical system
based on the recurrence relation $y_{t+1} \Leftarrow y_t$:
\begin{displaymath}
  \xymatrix{ 
    *++[o][F]{0} \ar@(dl,ul)[]^{\{0\}} &
    *++[o][F]{1} \ar@(ur,dr)[]^{\{1\}}
  }
\end{displaymath}
Either initial value $y_0$ gives a monotonous sequence of
solution-value sets for $y$ at successive times $t$:
\begin{eqnarray}
  \vec{\set{S}}_{\{\}}^*(y) & \Rightarrow &
  \left[ \begin{array}{rcl}
      y_0 = 0 & \mapsto & 
      (\{0\},\{0\},\{0\},\{0\},\ldots) \\
      y_0 = 1 & \mapsto & 
      (\{1\},\{1\},\{1\},\{1\},\ldots)
    \end{array} \right]
\end{eqnarray}
This dynamical system has two fixed points; hence it is
\emph{contingent} by Definition~\ref{def:dynamical}.  In contrast to
\goedel's specification in Equation~\ref{eq:nbew} that describes a
formula that cannot be realized consistently, the complementary
specification in Equation~\ref{eq:bew} describes a formula that can be
realized in two different ways, either of which is consistent with the
specification: in one case the formula is necessarily true (having the
solution-value set $\{1\}$), and in the other case the formula is
necessarily false (having the solution-value set $\{0\}$).  The
tautological vacuousness of self-affirmation mirrors the
viciously-circular contradiction of self-denial.

\end{example}

\begin{example}\exfmt
  Let us consider a set of axioms that gives a slightly more
  complicated dynamical system than \goedel's formula \rqq.  We now
  specify a logical formula $z$ with the following properties: $z$ is
  not provable; if $z$ is ambiguous, then it is not true; and if $z$
  is not provable, the negation of $z$ is not provable, and $z$ is not
  ambiguous, then $z$ is true.  Thus in logical notation we have the
  declaration $z \in \{\val{t},\val{f}\}$ and the following axioms:
  \begin{equation}
    \label{eq:zbew}
    \begin{array}{ll}
      \vdash & \neg \mathsf{Bew}(z) \\
      \vdash & \mathsf{Viel}(z) \rightarrow \neg z \\
      \vdash &
      (\neg \mathsf{Bew}(z) \wedge
      \neg \mathsf{Bew}(\neg z) \wedge
      \neg \mathsf{Viel}(z)) \rightarrow z
    \end{array}
  \end{equation}
  Here the notation $\mathsf{Viel}(z)$ is introduced for the
  proposition that $z$ is ambiguous (`viel\-deu\-tig' in German), in
  other words that its solution-value set $\set{S}_A(z)=\{0,1\}$.  

  We introduce the parameter $\theta$ to represent the solution-value
  set for $z$.  The domain of $\theta$ is the power set
  $\{\{\},\{0\},\{1\},\{0,1\}\}$ of the elementary values $\{0,1\}$.
  Following Definition~\ref{def:flat} we translate the axioms in
  Equation~\ref{eq:zbew} into the following parametric polynomial
  system (here we have delayed polynomial translation of the material
  implication operators and represented negations in a slightly
  different way for convenience):
  \begin{eqnarray}
    \begin{array}{rrcl}
      \tau_z: & z & \in & \{0,1\} \subset \mathbb{R} \\
      \tau_\theta: & \theta & \in & \{\{\},\{0\},\{1\},\{0,1\}\} \\
      \alpha_1: & (\theta=\{1\}) & = & 0 \\
      \alpha_2: & ((\theta=\{0,1\}) \rightarrow (1-z)) & = & 1 \\
      \alpha_3: & 
      ((\theta \neq \{1\})(\theta \neq \{0\})(\theta \neq \{0,1\}) 
      \rightarrow z) & = & 1 \\
      \delta_\theta: & \theta & \Leftarrow & 
      \set{S}_{\{\alpha_1,\alpha_2,\alpha_3\}}(z)
    \end{array}
    \label{eq:zbew-templ}
  \end{eqnarray}
  This parametric system of polynomial equations specifies the
  following dynamical system (each node is labeled with a state
  $\theta$, and each arc indicates the updated solution-value set
  calculated for $z$ assuming the state corresponding to the
  originating node):
  
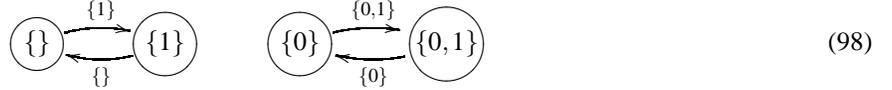
\begin{equation}
    \label{eq:zbew-graph}
    \xymatrix{
      *++[o][F]{\{\}}\ar@/^/[r]^{\{1\}} & 
      *++[o][F]{\{1\}}\ar@/^/[l]^{\{\}} &
      *++[o][F]{\{0\}}\ar@/^/[r]^{\{0,1\}} & 
      *++[o][F]{\{0,1\}}\ar@/^/[l]^{\{0\}}
    }
  \end{equation}
  Hypothetico-deductive analysis produces the following mappings for
  the evolution function $F$:
  \begin{equation}
    \label{eq:zbew-f}
    \{\} \mapsto \{1\}, \quad
    \{1\} \mapsto \{\}, \quad
    \{0\} \mapsto \{0,1\}, \quad
    \{0,1\} \mapsto \{0\}
  \end{equation}
  For example, given the hypothesis $\theta=\{\}$ the constraint
  templates in Equation~\ref{eq:zbew-templ} are instantiated as:
  \begin{equation}
    \begin{array}{rrcl}
      \alpha_1: & (\{\}=\{1\}) & = & 0 \\
      \alpha_2: & ((\{\}=\{0,1\}) \rightarrow (1-z)) & = & 1 \\
      \alpha_3: & 
      ((\{\} \neq \{1\})(\{\} \neq \{0\})(\{\} \neq \{0,1\}) 
      \rightarrow z) & = & 1
    \end{array}
  \end{equation}
  which simplify to:
  \begin{equation}
    (0) = 0, \quad
    ((0) \rightarrow (1-z)) = 1, \quad
    ((1)(1)(1) \rightarrow z) = 1
  \end{equation}
  Translating the material implication operators according to
  Table~\ref{tbl:translation} yields the polynomial constraints:
  \begin{equation}
    0 = 0, \quad 1 = 1, \quad z = 1
  \end{equation}
  These constraints lead to the solution-value set
  $\set{S}_{\{0=0,1=1,z=1\}}(z) \Rightarrow \{1\}$ for $z$, which is
  assigned as the new value of the parameter $\theta$ by the updating
  function in template $\delta_\theta$.  Thus is established the
  mapping $\{\} \mapsto \{1\}$ for the evolution function $F(\theta)$.
  The other hypotheses $\theta=\{1\}$, $\theta=\{0\}$, and
  $\theta=\{0,1\}$ are handled in a similar way.  The results are the
  same using coefficients in $\mathbb{R}$ or in $\mathbb{F}_2$.

  For this problem there is no direct representation of the evolution
  function $F(\theta)$ as a polynomial with coefficients in
  $\mathbb{R}$ or in $\mathbb{F}_2$, since the domain
  $\{\{\},\{0\},\{1\},\{0,1\}\}$ for $\theta$ is not a subset of
  either set of numbers.  Nonetheless it is evident from the mappings
  in Equation~\ref{eq:zbew-f} and the graph in
  Equation~\ref{eq:zbew-graph} that there are no fixed points in this
  dynamical system; there is no state of $\theta$ that maps to itself.
  Hence the formula $z$ defined by the axioms in
  Equation~\ref{eq:zbew} is unsteady in the dynamic interpretation and
  unsatisfiable in the static interpretation.  Comparing the graphs in
  Equations \ref{eq:zbew-graph} and \ref{eq:bew-graph} you can see
  that the dynamic behavior of this new formula $z$ is a bit different
  from that of \goedel's \rqq.  For example the new formula $z$
  produces two periodic orbits instead of one.  Also, assuming that
  $z$ is a theorem never leads to the conclusion that the negation of
  $z$ is a theorem (as it does for \goedel's formula); for the new
  formula $z$ these states are in different, disconnected orbits.
\end{example}

\subsection{Quantifier Elimination}

It is not necessary to use quantifiers to examine \goedel's
incompleteness argument; however for the sake of completeness let us
address briefly how quantifiers should be handled in this algebraic
formulation of logic.  Existential and universal quantifiers add
complexity to polynomial formulas and systems of constraints.  In the
field of computational algebraic geometry, there is a large body of
work on eliminating quantifiers from polynomial formulas with
real-number coefficients \cite{tarski-algebra,arnon}.  Since Boole's
translation method allows logical formulas to be converted into
polynomials with ordinary numeric coefficients, such algebraic
geometry methods for managing quantifiers are applicable to problems
that were originally specified in traditional logical notation.
Eliminating quantifiers from an equation may change the overall system
of equations from a simple conjunction (all equations are to be
satisfied simultaneously) to a disjunction (various subsets of
equations must be solved separately).

What is difficult about quantifiers is not the quantifiers themselves,
but the indeterminate functions of infinite domains that often
accompany them.  Indeterminate functions of \emph{finite} domains can
be treated parametrically using Lemma~\ref{lem:matrix}.  However
indeterminate predicates such as $\mathsf{Man}(x)$ and
$\mathsf{Mortal}(x)$ assumed to have an infinite universe of discourse
for $x$ (as $x\in\mathbb{N}$) are problematic.  Perhaps it would be
possible to choose a `quorum' (sufficiently-large but finite universe
of discourse) based on the number and arity of predicates declared in
a logical system.  Anyway there is more work to be done to connect the
treatment of quantifiers in computational algebraic geometry with the
treatment of quantifiers in classical first-order logic.

\section{Discussion}

I have argued two major propositions in this essay: that logic is an
application of mathematics, more specifically solving systems of
polynomial equations to yield sets of feasible elementary values; and
that allowing equations to refer to their own solutions creates
discrete dynamical systems.  In this analysis I have employed Boole's
groundbreaking algebraic formulation of logic, revised and extended
with concepts from algebraic geometry and computer science.  Within
this framework of dynamic polynomial logic, \goedel's special formula
that asserts its own unprovability is easily expressed and is not at
all paradoxical.  Furthermore, \goedel's formula is not so special
after all: it is easy to modify ordinary quadratic equations to
exhibit the same behavior, or to specify other logical formulas that
encode more elaborate dynamical systems.  Kurt \goedel{} provided
innovative but convoluted proofs of some normal properties of a simple
recurrence relation, accompanied by spectacular misinterpretations.  

A familiar result from geometry may provide insight into the subtle
cognitive error made by \goedel{} and many classical logicians: what I
shall call the Pythagorean fallacy.

\subsection{The Pythagorean Fallacy}

The Pythagorean theorem holds an important lesson about undecidability
and incompleteness.  We all know the famous equation $a^2+b^2=c^2$
relating the lengths of the sides and hypotenuse of a right triangle.
Perhaps less well known is the fact that Pythagoras and his followers
initially considered the only valid numbers to be ratios of natural
numbers.  To them, the solution for $c$ with some values of $a$ and
$b$ (such as $a=1$ and $b=1$) was \emph{alogos}
(\begin{greek}>al'ogos\end{greek}): though we typically say
  `irrational', a more literal translation would be `not reasonable'.
  Of course, to the modern thinker it is not a number like $\sqrt{2}$
  that seems unreasonable; it is the expectation that an algebraic
  equation with natural-number coefficients must have a solution that
  is the ratio of two natural numbers.

Let us first recognize the phenomenon of the \emph{anticipated-actual
  type mismatch}: the expectation that the value of a mathematical
expression should be one type of object when it is in fact another.
It is an anticipated-actual type mismatch to expect the expression
$\sqrt{a^2+b^2}$ to have a rational value for all natural-number
arguments $a$ and $b$.  The anticipated-actual type mismatch is
facilitated by expressions that use objects of one type to construct
objects of another type: for example whereas the sum or product of any
integers must itself be an integer, the quotient of integers may not
be.  We can recognize that certain \emph{mechanisms} of constructing
expressions lead to certain types of mathematical values: for example
if we are allowed the square root operation then we can make algebraic
numbers from integers, and if we are allowed infinite series then we
can make transcendental numbers from rational numbers.  Notation and
terminology may obscure the actual mechanisms in use.  For example
there is no square root sign in the equation $a^2+b^2=c^2$, yet we
must invoke this operation to solve for $c$.

Let us call it the \emph{Pythagorean fallacy} to encounter an
anticipated-actual type mismatch but then to misplace the blame: to
decide that the actual value, rather than the incorrect expectation,
is at fault.  This error is facilitated by formulas that have the
expected type for some but not all arguments (in other words the
expected type is a special case of the actual type): for example the
value of $\sqrt{3^2+4^2}$ is indeed a rational number ($5$) but the
value of $\sqrt{1^2+1^2}$ is not any rational number.  We may
recognize certain circumscribed \emph{exceptions} in actual values:
for example whereas dividing one rational number by another generally
yields a rational number, division by zero does not.  It is not
problematic to recognize such limited exceptions.  But it is
inappropriate to indulge the Pythagorean fallacy by declaring that
some general mathematical method is fundamentally defective because of
an anticipated-actual type mismatch, when the fault lies with an
incorrect expectation.

Ancient Greek mathematicians (and perhaps Pythagoras himself)
eventually resolved their Pythagorean fallacy in an exemplary way:
they broadened their concept of number to include the irrational
values necessary for the general solution of algebraic equations like
$a^2+b^2=c^2$.  Our modern concept of number has grown to include many
more `unreasonable' types of objects: not only algebraic numbers like
$\sqrt{2}$, but also transcendental numbers like $e$ and $\pi$ and
imaginary numbers like $\sqrt{-1}$.  We have come to understand that
the mechanisms of geometric and algebraic construction ought to
produce irrational numbers sometimes, even when the input values are
natural numbers.  We would find it absurd to accuse geometry or
algebra of undecidability or incompleteness because there are
expressions whose values are not ratios of natural numbers.

\subsection{\goedel's Error}

Returning to the main results from \cite{goedel}, what \goedel's
Theorem~VI says is technically true but not detrimental to formal
reasoning as he and many others have concluded.  However it is
emphatically not the case that \goedel's special formula is
semantically true; instead this formula that asserts its own
unprovability is exceptional in a particular way.  The types of
`undecidability' and `inconsistency' that have been identified as
exceptions---unsatisfiability, ambiguity, unsteadiness, and
contingency---are semantically and syntactically appropriate features
of systems of equations.  Moreover, these exceptions can just as well
be demonstrated in ordinary polynomial equations which had no origin
in traditional logical notation.  For example it is evident that the
formula $x$ has no decidable real-number value given the definition
that $x$ is a real number satisfying $x^2=-1$; and that the formula
$y$ has no definite real value when $y$ is defined as a real number
satisfying $y^2=y$; and that the formula $z$ has no definite integer
value when $z$ is defined as an integer satisfying $z_{t+1} \Leftarrow
1-z_t$.  Each of these formulas $x$, $y$, and $z$ is exceptional in a
particular way and should not have a definite elementary value.  There
is no need to divorce syntax from semantics.  \emph{Undecidability is
  not a crime!}

\appendix

\clearpage
\section{Executive Summary}

In his famous 1931 paper \goedel{} described a special logical formula
that asserts of itself that it cannot be proven, in the context of a
formal system PM that can express logical formulas, theorems, proof,
and natural numbers.  \goedel{} claimed that this special formula must
be semantically true but syntactically \emph{undecidable}: impossible
to prove or disprove by calculations within PM.  \goedel{} concluded
that the existence of such a true-but-undecidable formula renders any
formal system like PM essentially incomplete and incapable of proving
even its own consistency.  In this presentation I shall introduce a
new method of analysis called \emph{computational algebraic logic}
that proves most of these conjectures wrong.  Although \goedel's
special formula is indeed `undecidable' in the sense that it has a
value outside of the set $\{\textsc{true},\:\textsc{false}\}$, such
undecidability is a feature rather than a bug in formal reasoning: the
formula is not semantically true.

There are three main ideas to be discussed: exceptions, translation,
and dynamical systems.  First, mathematicians are unfailingly perplexed
by formulas that give values outside of the sets they had originally
expected; they give pejorative names to the unexpected
\emph{exceptional} values. \emph{Natural} numbers are usually welcome.
But others are called \emph{negative}, \emph{fractional} (broken!),
\emph{irrational}, or \emph{imaginary}.  If you knew nothing about
complex numbers, you might say that $\sqrt{-1}$ is \emph{undecidable}:
the value of this formula cannot be proven to be any \emph{real}
number, even though its argument is real.  Likewise $\sqrt{2}$, $\pi$,
$e$, and so on seem exceptional and undecidable if you are expecting
\emph{rational} numbers.  \goedel's \emph{undecidable} formula is just
another unexpected object in this litany of un-\emph{natural} numbers:
the best solution is to find the right data structure to represent it,
not to indict formal reasoning as incomplete and inconsistent.

Second, translation from logical to algebraic notation helps to
uncover the right data structure.  You would not want to do your
calculus homework or balance your checkbook in Roman numerals;
likewise it is difficult to perform logical inference in traditional
notation (our mash-up from Frege, Peano, Hilbert, and a few others).
Boole already showed in 1854 that formulas in binary logic can be
translated into ordinary polynomials, and that logical deduction can
be performed by solving systems of polynomial equations.  Boole's work
is widely misrepresented; for example, he used ordinary arithmetic in
which $1+1=2$, not what is now called `Boolean algebra'.  Anyway,
using Boole's original method, the logical axiom $p$ translates as the
equation $p=1$ and the axiom $p \rightarrow q$ translates as the
equation $pq-p+1=1$, using real-valued variables $p$ and $q$.  The
constraints $p^2=p$ and $q^2=q$ ensure that each variable must be
either $0$ or $1$.  The truth value of the formula $q$, for which I
introduce the notation $\set{S}(q)$, is given by the set of solutions
to these equations:
\begin{eqnarray*}
  \set{S}(q) & = & \left\{ \; q \; : \;
  (p,q) \in \mathbb{R}^2, \; p=1,\; pq-p+1=1,\; p=p^2,\; q=q^2 \; 
  \right\}
\end{eqnarray*}
This system of polynomial equations is easy enough to solve by
inspection; the only solution is $(p,q)=(1,1)$ which gives the
\emph{solution-value set} $\set{S}(q)=\{1\}$.  So the only feasible value
for $q$ is $1$; therefore, we call $q$ a \emph{theorem}.  Voil\`a,
\emph{modus ponens} by elementary algebra.  What was missing from
Boole's work was a general method to solve systems of multivariate
polynomial equations.  Luckily for us, Buchberger invented such a
method in the 1960s; today his Gr\"obner-basis algorithms are widely
implemented in computer algebra systems.

Next, what does \goedel's special formula look like in algebraic
notation?  To define a formula $x$ that is true if and only if $x$
itself is not provable, we can use a real-valued variable $x$ and the
equation $x = (\set{S}(x) \neq \{1\})$ in which the idea `$x$ is not
provable' is expressed by saying that the solution-value set
$\set{S}(x)$ is not equal to the set $\{1\}$ (we adopt the usual
computer-programming convention that the test of inequality $\neq$
returns the numeric value $0$ if its arguments are equal and $1$ if
they are not equal).  This self-referential equation turns out to be a
\emph{recurrence relation} which is equivalent to the recurrence
$x_{t+1} \Leftarrow 1 - x_{t}$ with each $x_{t}\in\{0,1\}$.  The
mathematical object that is defined is a \emph{discrete dynamical
  system} with two states ($x=0$ and $x=1$) and a periodic orbit that
oscillates between them:
\begin{displaymath}
  \xymatrix{
    *++[o][F]{0} \ar@/^/[r] &
    *++[o][F]{1} \ar@/^/[l]
  }
\end{displaymath}
The graph shows exactly what is exceptional about \goedel's special
formula: it defines system that has no \emph{fixed points} (steady
state).  This dynamical system encodes for $x$ an infinite sequence
$(\ldots,\{0\},\{1\},\{0\},\{1\},\ldots)$ of alternating
solution-value sets, with a phase determined by what initial condition
$x_0$ was assumed.  

Dynamical systems are the `complex numbers' that capture the value of
\goedel's special formula and similar self-referential systems of
equations.  These dynamical systems are perfectly computable from
self-referential formulas by ordinary means and their existence does
not jeopardize the consistency of the formal system within which they
are derived (elementary algebra).  I do not think that any
mathematician or philosopher is troubled that the Fibonacci numbers
are not all the same; we should understand that \goedel's special
formula is the same kind of construction and that it yields a
similarly dynamic sequence of elementary values.  It is no more
correct to say that \goedel's special formula is semantically true or
false than it is to say that Fibonacci's formula is semantically zero,
one, or two (or three, or five, or eight; you get the idea).

Within this framework of algebraic logic, many other types of
exceptions can be recognized through computation.  To complement
\goedel's special formula we might introduce a new formula $y$ that is
true if and only if it is provable: $y=(\set{S}(y)=\{1\})$.  This is
equivalent to the recurrence relation $y_{t+1} \Leftarrow y_{t}$ with
each $y_t\in\{0,1\}$.  What is defined is a discrete dynamical system
with two fixed points:
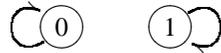
\begin{displaymath}
  \xymatrix{ 
    *++[o][F]{0} \ar@(dl,ul)[] &
    *++[o][F]{1} \ar@(ur,dr)[]
  }
\end{displaymath}
Such a dynamical system is \emph{contingent}.  It encodes for $y$ two
different infinite sequences of solution-value sets.  From the initial
condition $y_0=0$ that $y$ is not a theorem the sequence
$(\{0\},\{0\},\{0\},\ldots)$ confirms that $y$ is never a theorem.
And from the initial condition $y_0=1$ that $y$ is a theorem the
sequence $(\{1\},\{1\},\{1\},\ldots)$ confirms that $y$ is always a
theorem.  

Rounding out the host of exceptional formulas are those whose
solution-value sets are empty and those whose solution-value sets have
multiple members.  Using Boole's translation $1-y$ for $\neg y$, we
might model the liar paradox with the equations $z=1-z$ and $z^2=z$
for $z\in\mathbb{R}$; in this case $\set{S}(z)=\{\}$ and the formula
$z$ is \emph{unsatisfiable} (just like the equation $0=1$).  We might
add the complementary truth-teller problem using $w=w$ and $w^2=w$
with $w\in\mathbb{R}$; in this case $\set{S}(w)=\{0,1\}$ and the
formula $w$ is \emph{ambiguous}.

What \goedel{} saw as \emph{undecidability} does not show some
foundational defect in formal reasoning any more than irrational or
imaginary numbers invalidate arithmetic.  The myth of incompleteness
is dispelled when self-reference is unmasked as a recurrence relation
that defines a dynamical system.  \goedel's incompleteness argument is
commonly understood to prove that logic and mathematics are
fundamentally incompatible, but in fact he demonstrated quite the
opposite: for formal reasoning it is essential to recognize that
problems in logic are problems in algebra, just as Boole had
demonstrated.

\clearpage

\section{Russell's Paradox: Recursion Redux}

The methods presented in this analysis of \goedel's incompleteness
theorems can be applied to another well-known problem in logic,
Bertrand Russell's paradox (introduced in a 1902 letter to Frege that
is reproduced in \cite{heijenoort}).  Russell asked whether a special
set---the set of all sets that do not contain themselves---contains
itself or not; either answer seems paradoxical.  Although Russell's
problem does not translate directly as a system of polynomial
equations, it is nonetheless a system of mathematical constraints and
the same insights apply.  Russell's construction is unsatisfiable in
the static interpretation and unsteady in the dynamic interpretation;
it behaves exactly like \goedel's special formula.  Furthermore, a bit
of reflection reveals that the specific issue with Russell's paradox
is neither totality nor unrestricted comprehension but indefiniteness:
his construction requires the set being defined to be used as a free
variable in the predicate that defines it, which violates even the
unrestricted axiom schema of comprehension.  Regarding Russell's
paradox, it is sets of humans not axioms of sets that have a problem
with comprehension.

\subsection{Algebraic Analysis}

Let us define $\mathbb{S}$ as a set of sets that includes at least
some set $r$, and further define $r$ as the set of all members of
$\mathbb{S}$ that are not members of themselves.  We desire the truth
value of the formula $r \in r$ that asks whether the set $r$ is a
member of itself.  The definitions of $r$ and of $\mathbb{S}$ provide
two constraints that constitute a generalized (beyond polynomial
equations) system of constraints:
\begin{eqnarray}
  r & := & \{s \in \mathbb{S} : s \notin s\} \label{eq:r} \\
  r & \in & \mathbb{S} \label{eq:s}
\end{eqnarray}
In the special case that $\mathbb{S}$ is considered to be the set of
all possible sets (assuming for the moment that such a construct is
meaningful), this pair of constraints amounts to the classic
formulation of Russell's paradox.  But Equations \ref{eq:r} and
\ref{eq:s} are not limited to that special case; in particular
$\mathbb{S}$ may denote a countable and finite set of sets and it is
not important whether $\mathbb{S}$ contains itself.  It is, however,
essential to the problem that $\mathbb{S}$ contains $r$.

In order to ascertain whether the set $r$ thus defined contains
itself, we use the solution-value set for the formula $r \in r$
subject to the constraints in Equations \ref{eq:r} and \ref{eq:s}:
\begin{equation}
  \set{S}_{\left\{
    r = \{s \in \mathbb{S} : s \notin s\},\: 
    r \in \mathbb{S} 
    \right\}}
  (r \in r)
  \label{eq:russell}
\end{equation}
Using the methods in Section~\ref{sec:dynamic} there are static and
dynamic ways to interpret this solution-value set.  In terms of
Definition~\ref{def:parametric} we use the value of the set-membership
expression $r \in r$ as the parameter $\theta$ for this problem.  We
note that the expression $r \in r$ takes a value in the Boolean finite
field $\mathbb{F}_2:=\{0,1\}$ according to whether the relation in it
is false or true.  Thus if the relation $r \in r$ is false we have $(r
\in r)=0$ (equivalently $r \notin r$) and $\theta=0$; conversely if
the relation is true we have $(r \in r)=1$ (equivalently $r \in r$)
and $\theta=1$.

Using a hypothetico-deductive approach we consider these two possible
cases $\theta=0$ and $\theta=1$ for the parameter $\theta=(r \in r)$
as illustrated in Table~\ref{tbl:cct-russell}.  In each instantiation
$A(\theta_t)$ the conjectured value of the formula $r \in r$ is
specified as an explicit constraint (since it cannot be substituted
directly into an existing constraint as in the earlier examples).
From either conjectured value we compute the other.  That is, from the
conjecture $\theta=\theta_1=0$ indicating $r \notin r$ the constructed
set $r$ includes itself and hence the computed value $\theta_1=(r \in
r) \Rightarrow 1$.  Likewise from the conjecture $\theta=\theta_2=1$
indicating $r \in r$ the constructed set $r$ does not include itself
and so the computed value $\theta_2=(r \in r) \Rightarrow 0$.  Thus
from Equations \ref{eq:r} and \ref{eq:s} we derive the state
transitions $0 \mapsto 1$ and $1 \mapsto 0$ defining the following
dynamical system:
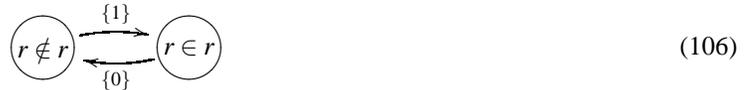
\begin{equation}
  \xymatrix{
    *++[o][F]{r \notin r} \ar@/^/[r]^{\{1\}} & 
    *++[o][F]{r \in r} \ar@/^/[l]^{\{0\}}
  }
  \label{eq:graph-russell}
\end{equation}
with the states labeled $r \notin r$ and $r \in r$ instead of
$\theta=0$ and $\theta=1$ for convenience.  Since this dynamical
system has no fixed points, we conclude in this dynamic interpretation
that the formula $r \in r$ is \emph{unsteady} subject to the
constraints in Equation \ref{eq:r} and \ref{eq:s}.  Moreover, in the
static interpretation the formula $r \in r$ is \emph{unsatisfiable}
given the constraints in Equation \ref{eq:r} and \ref{eq:s} since its
conjectured and computed values never agree.  

\begin{table}
  \caption{The state-transition worksheet for Russell's set-inclusion
    query using the parameter $\theta := (r \in r)$.}
  \label{tbl:cct-russell}
  \begin{tabular}{r@{\qquad}c@{\quad}c@{\qquad}l@{\qquad}l} \hline\hline
    $i$ & $\theta_i$ & $A(\theta_i)$ &
    $\set{S}_{A(\theta_i)}(r \in r)$ & $F(\theta_i)$ \\ \hline
    $1$ & $0$ & 
    $\{ r = \{s \in \mathbb{S} : s \notin s\},\: r \in \mathbb{S},\: 
    r \notin r\}$ &
    $\{1\}$  & $1$ \\
    $2$ & $1$ & 
    $\{ r = \{s \in \mathbb{S} : s \notin s\},\: r \in \mathbb{S},\:
    r \in r\}$ &
    $\{0\}$ & $0$ \\ \hline\hline
  \end{tabular}
\end{table}

Russell's special formula specifies exactly the same dynamical system
as \goedel's special formula: an oscillator with two states.  This
coincidence is not surprising: when the definition of a formula
includes a recursive reference to a solution-set feature that has two
possible states, the only way to make an unstable/unsatisfiable system
is to have the conjecture of either state lead to the computation of
the other.  Such instability is exactly what seems paradoxical about
formulas like \goedel's and Russell's (for some reason the
contingent/ambiguous dual of each problem is not perceived as
paradoxical).

\subsection{Russell's Free Variable Revealed}

In more conventional language, Russell's construction necessarily
employs the set being defined as a free variable in the predicate that
defines it---which renders that predicate indefinite.  It is important
to recognize such indefinite predicates, which Zermelo recognized as
problematic in his axiomatization of set theory \cite{zermelo}.  In
particular, when specifying the construction of a set by the action of
a predicate upon the members of some universe of discourse, the
members of that universe of discourse must be considered when
evaluating the predicate for free variables.

\begin{lem}[Revealing Free Variables in Set-Building]
  \label{lem:free}
  Consider the definition of a set $y$ using some universe of
  discourse $z$ and some predicate $\phi(x)$:
  \begin{eqnarray*}
    y & := & \{ x \in z : \phi(x) \}
  \end{eqnarray*}
  Every member of the set $z$ that is a variable is in fact a
  \emph{free variable} in the predicate $\phi(x)$.  In particular, if
  the set $y$ being defined is a member of the set $z$ used in its own
  definition, then the instantiation $\phi(y)$ of the predicate (with
  $y$ as the argument) must be evaluated during the construction of
  $y$.  It is evident that $y$ is a free variable in this instance
  $\phi(y)$.  Leaving the universe of discourse anonymous as in the
  unrestricted definition $y:=\{x:\phi(x)\}$ does not change this
  property; if $y$ is a possible value for $x$ then the variable $y$
  must still be considered free in the predicate $\phi(x)$.  Actually,
  it is not just the predicate $\phi(x)$ that is at issue; it is the
  entire right-hand-side expression $\{x \in z:\phi(x)\}$ in the
  definition.  Even with a definition such as $y:=\{x \in z: \val{t}\}$
  whose trivial predicate $\val{t}$ (the elementary value true) has no
  variables, the definiendum $y$ must be considered a free variable in
  the definiens $\{x \in z: \val{t}\}$ if and only if $y \in z$;
  similarly, the definiendum $y$ is a free variable in the definiens
  $\{x:\val{t}\}$ of the unrestricted definition $y:=\{x:\val{t}\}$.
\end{lem}
It is important to recognize that Zermelo's \cite{zermelo} axiom
schema of separation (Axiom~III) contains two clauses: a separation
clause and a definiteness clause.  That axiom is repeated here:
\begin{quote}
  Whenever the propositional function $\mathfrak{E}(x)$ is definite
  for all elements of a set $M$, $M$ possesses a subset
  $M_\mathfrak{E}$ containing as elements precisely those elements $x$
  of $M$ for which $\mathfrak{E}(x)$ is true.
\end{quote}
\noindent The idea of a \emph{definite} propositional function
(addressed in Lemma~\ref{lem:free}) already entails the principle of
separation, in the sense that unrestricted set comprehension always
implies an indefinite predicate.  For example, using the unrestricted
definition $y:=\{x:\phi(x)\}$ the variable $y$ must be free in the
predicate $\phi(x)$ by the argument above; hence the predicate is
indefinite.  To be sure, an expression that includes free variables
could still have a definite value if other constraints in the system
limit the values of those variables; for the moment we assume the
absence of such additional constraints.

However, requiring separation without requiring definiteness does not
close the loophole in naive set theory.  Simply introducing some
independently-defined set-of-sets $z$ to restrict set comprehension as
in $y:=\{x \in z: \phi(x)\}$ is not sufficient to ensure definiteness;
in the case $y \in z$ the variable $y$ remains free in $\phi(x)$ and
thus the predicate remains indefinite (at some point in the
construction of $y$ the predicate $\phi(y)$ must be evaluated to
determine whether to include $y$ in itself).  That is all to say that
an \emph{unrestricted} axiom schema of comprehension, retaining the
definiteness clause from Zermelo's Axiom~III but omitting its
separation clause, would be sufficient to prevent Russell's paradox.
Conversely a restricted axiom schema of comprehension retaining the
separation clause but omitting the definiteness clause would not
prevent the paradox.  The principle that a \emph{definition} requires
a \emph{definite} expression in its definiendum is not at all peculiar
to set theory!

Returning to Equation~\ref{eq:r}, in the definition $r:=\{s \in
\mathbb{S} : s \notin s\}$, the variable $r$ is free in the predicate
$s \notin s$ exactly if $r$ is a member of $\mathbb{S}$.  Hence using
this restricted definition it would violate the definiteness clause in
Zermelo's Axiom~III to consider $r$ a set exactly if $r \in
\mathbb{S}$.  In the unrestricted definition $r:=\{s: s \notin s\}$,
the variable $r$ is always free in the predicate $s \notin s$ since
$r$ is considered to be a member of the universe of discourse from
which the values of $s$ are drawn.  Hence both clauses of Zermelo's
Axiom~III would be violated by the unrestricted definition.  Using
either the restricted or the unrestricted formulation, at some point
the predicate to be evaluated must become the expression $r \notin r$
in which $r$ is plainly a free variable; thus the predicate is
indefinite.  The precise issue with Russell's paradox is neither
totality nor unrestricted comprehension; it is the attempt to define a
set using an indefinite predicate.

One can imagine a set theory in which Zermelo's definiteness clause is
relaxed and thereby it is allowed to use the set being defined as a
free variable in its own definition.  In such a liberalized set theory
it would be possible to use recurrence relations to define dynamical
systems that give infinite sequences of sets, just as it is possible
to use a recurrence relation to define the Fibonacci sequence in
elementary algebra.  Russell's construction would give one of these
dynamical systems, characterized by a periodic orbit that oscillates
between two states as shown in the graph in
Equation~\ref{eq:graph-russell}.  The unrestricted `set of all sets'
(as from $y:=\{x:\val{T}\}$) would define a different dynamical system,
as if by the recurrence relation $y_{t+1} \Leftarrow y_{t} \cup
\{y_{t}\}$ (this system would have no fixed points, but orbits that
are neither periodic nor convergent).  As in the earlier examples,
these dynamical systems could instead be interpreted as static
constraints (both Russell's special set and the `set of all sets'
$y:=\{x:\val{T}\}$ would be unsatisfiable in this interpretation).

The method of algebraic analysis detailed above illustrates exactly
which exceptions can occur if a set $y$ is allowed to be a free
variable in its own definition (either from the explicit circumstance
$y \in z$ or from unrestricted comprehension): unsteadiness and
contingency.  (Note that a self-referential set $y$ could also have a
steady, unique value and be unexceptional.)  In this sense we gain the
ability to discern more details about the various paradoxes in naive
set theory, for example to clarify how `the set of all sets that do
not contain themselves' differs mathematically from `the set of all
sets' not otherwise specified (the constructions give structurally
distinct dynamical systems, one with a periodic cycle and the other
without).  We can also discover similarities between exceptional
constructions in set theory and exceptional constructions in the rest
of logic and algebra, for example the congruence between Russell's
`set of all sets that do not contain themselves' and \goedel's
`logical formula that asserts its own unprovability.'

\subsection{Functional Programming}

As a final thought on Russell's non-paradox, it is interesting to me
as a computer programmer that the passage in Frege's
\emph{Begriffsschrift} that inspired Russell to create his
set-inclusion problem also has a computer-programming interpretation.
This is how Russell introduced his problem \cite{russell-letter}:
\begin{quote}
  There is just one point where I have encountered a difficulty.  You
  state (p.~17) that a function, too, can act as the indeterminate
  element.  This I formerly believed, but now this view seems doubtful
  to me because of the following contradiction.
\end{quote}
\noindent Indeed, in the referenced \S{}9 of \emph{Begriffsschrift}
Frege had written, ``On the other hand, it may be that the argument is
determinate and the function indeterminate.''  This is idea can be
explored in terms of predicates and sets, as Russell did.  But it is
also a fairly straightforward description of \emph{functional
  programming}---the idea that functions should be treated as
first-class mathematical objects (like numbers).  This is implemented
in several computer programming languages (most prominently Lisp
\cite{lisp}).  Even beyond the conceptual appeal of such an
object-oriented approach, it can be practical to solve systems of
equations with functions as unknowns, when the type of function is
suitably restricted (e.g.\ to polynomials over a given set of
variables with coefficients in the finite Boolean field
$\mathbb{F}_2$; Boole's function development method allows parametric
descriptions of arbitrary logical functions of any desired arity).
Some challenging problems in logic, such as Smullyan's puzzles about
liars and truth-tellers \cite{smullyan-name}, amount to solving for
unknown logical functions rather than unknown elementary truth values.
The calculus of infinitesimals is another place in mathematics where
functions are treated as first-class objects.  What is integration
except solving a system of constraints (on rates of change and perhaps
boundary values) in which a function is the unknown variable?

\clearpage
\section{Notes on Modal Logic: Thinking Outside the Box}

We can use dynamic polynomial logic to define modal operators in terms
of solution-value sets.  First the alethic modes corresponding to the
members of $2^{\{0,1\}}$ (here using $s$ to abbreviate the
solution-value set $\set{S}_A(p)$ of the objective
$p(x_1,x_2,\ldots,x_n)$ relative to the set $A$ of axioms):
\begin{itemize}
  \item $\square$ (or $\square_1$) necessarily true: $\square(p)$ means $s =
    \{1\}$ and $\square(\neg p)$ means $s = \{0\}$
  \item $\square_0$ necessarily false: $\square_0(p)$ means $s = \{0\}$ and
    $\square_0(\neg p)$ means $s = \{1\}$
  \item $\bowtie$ ambiguous (mnemonic $\infty$ for all truth values):
    $\bowtie(p)$ or $\bowtie(\neg p)$ mean $|s|>1$ thus $s = \{0,1\}$
  \item $\oslash$ unsatisfiable (mnemonic empty set $\emptyset$):
    $\oslash(p)$ or $\oslash(\neg p)$ mean $s = \{\}$
\end{itemize}
Then some hybrids, each corresponding to several members of
$2^{\{0,1\}}$:
\begin{itemize}
  \item $\lozenge$ (or $\lozenge_1$) possibly true: $\lozenge(p)$
    means $1 \in s$ thus $s \in \{\{1\},\{0,1\}\}$; \\ and
    $\lozenge(\neg p)$ means $0 \in s$ thus $s \in
    \{\{0\},\{0,1\}\}$
  \item $\lozenge_0$ possibly false: $\lozenge_0(p)$ means $0 \in s$
    thus $s \in \{\{0\},\{0,1\}\}$; \\ and $\lozenge_0(\neg p)$
    means $1 \in s$ thus $s \in \{\{1\},\{0,1\}\}$
  \item $\boxdot$ definite (mnemonic `get to the point'): $\boxdot(p)$ or
    $\boxdot(\neg p)$ mean $|s|=1$ thus $s\in\{\{0\},\{1\}\}$
\end{itemize}
Perhaps we should use $\square(p|A)$ or
$\square(p|q_1,q_2,\ldots,q_m)$ etc.\ to make the axioms explicit.

These operators and their negations cover $14$ of the $16$ members the
power set $2^{2^{\{0,1\}}}$.  Of the remaining $2$ members
$s\in\{\{\},\{0\},\{1\},\{0,1\}\}$ is tautological and $s\in\{\}$ is
impossible (even the empty set is not a member of itself).  Modal
operations are recurrence relations; they must be evaluated as
dynamical systems (though it is fine to report static
interpretations).  All these modal operators have the property that
they can be evaluated during hypothetico-deductive analysis: given a
computed solution-value set $s$ it is possible to tell whether any
operator is satisfied using standard set operations (equality,
inequality, membership, complement, cardinality).  For multiple
variables $(x_1,x_2,\ldots,x_n)$ we should probably use the solution
set of all variables as the state of the dynamical system.

For $d$-valued logic we could introduce $d$ distinct modal necessity
operators $\square_0,\ldots,\square_{d-1}$ with each $\square_k$
meaning $s=\{k\}$.  The remaining $2^d-d-1$ non-empty solution-value
sets in the power set of the set of elementary values (each
necessarily with cardinality $> 1$) would satisfy $\bowtie$ (thus
making it hybrid rather than singleton).  We could likewise introduce
$d$ distinct modal possibility operators where each $\lozenge_k$ means
$k \in s$.  Whatever special value $k^*$ is considered `true' would
customarily be omitted from the operator sign: thus the unadorned
$\square$ and $\lozenge$ instead of $\square_{k^*}$ and
$\lozenge_{k^*}$.

Note that the position of negation matters.  For example $\neg(\square
p)$ means $s \neq \{1\}$ thus the solution-value set $s$ could be
$\{\}$, $\{0\}$, or $\{0,1\}$; whereas $\square(\neg p)$ means the
solution-value set $s$ is definitely $\{0\}$.  In other words we
distinguish between `not necessarily true' and `necessarily false'.
Note also that `possible' here does not follow classical modal logic:
$\lozenge(p)$ is not equivalent to $\neg \square (\neg p)$ because the
latter includes the possibility $s=\{\}$ that $p$ is unsatisfiable
whereas the former does not.  That is, we distinguish between
`possible' and `not necessarily false'.  Also $\bowtie(\neg p)$ means
the same as $\bowtie(p)$: the negation of an ambiguous formula is
itself ambiguous.

In this framework $\square(p)$ does not mean the same thing as $p$.
In fact the axiom $\vdash \square(p)$ is quite a different assertion
from the axiom $\vdash p$; in isolation, the former specifies a
dynamical system with a fixed point at the state that $p$ is
unsatisfiable.  However the solution-value set $\set{S}_A(\square(p))$
evaluates to $\{1\}$ exactly if the solution-value set $\set{S}_A(p)$
also yields $\{1\}$; in this sense $\square(p)$ is similar to $p$ when
these expressions are treated as objectives.

Only unsatisfiability is `explosive': from a contradiction it follows
that all formulas are unsatisfiable (including the constants $0$ and
$1$: e.g.\ $\set{S}_{\{1=0\}}(1) := \{1:1=0\} \Rightarrow \{\}$).
Thus if $\oslash(p|A)$ is true for any formula $p$ then it must be
true for all formulas subject to the same axioms $A$.  Also, the sense
in which algebraic proof is monotonic is that adding constraints
cannot add members to solution sets; it can leave them unchanged or
remove members.  Thus for any sets $A$ and $B$ of axioms and any
formula $p$ we have $\set{S}_{A \cup B}(p) \subseteq \set{S}_{A}(p)$
and $\set{S}_{A \cup B}(p) \subseteq \set{S}_{B}(p)$.  It is not the
case that $\square(p|A)$ guarantees $\square(p|A,B)$ since the
additional axioms $B$ could render the whole system unsatisfiable;
thus we could have $\square(p|A)$ but $\oslash(p|A,B)$.  Algebraic
proof is global not local; all constraints must be taken into account,
and the solution set $\set{S}$ is what is left after all infeasible
solution values have been eliminated.

\clearpage
\section{By George! Proposed Logic Query Language}

Sketch of commands for a proposed computer system that implements
`Logic Query Language' incorporating dynamical systems, Boolean
translation from logic to algebra, and polynomial equations.

\begin{verbatim}
/* george: computational algebraic logic (LQL) */

parameter c in {0,1,2};
b = 2;      // simple macro assignment, resolved at compile time
c := |$x|;  // update rule: c gets size of solution-value set $x of x
x^2 + b*x + c == 0;     // equation template
x^2 + b*x + |$x| == 0;  // same constraint, implict parameter definition

parameter t in FF(2);
t := $y == {1};   // update rule: t gets true if y was proved a theorem
t := ?y;          // same update rule: t gets true if y was proved a theorem
y == 1 - t;       // template: y true if it was not proved a theorem
|- y <-> !t;      // same constraint in logical notation, |- for assert
y == ($y != {1}); // same constraint with implicit parameter definition
y == !?y;         // same with modal operator ? (necessarily true)
y = !?y;          // interpret as static constraint or recurrence?

// Fibonacci-like recurrence relation: update rules but no constraints
parameter t1, t2 in NN;
t1 := t2;
t2 := t1 + t2;
t1[0]=0; t2[0]=1;  // initial conditions give t1[t] : (0,1,1,2,3,5,8,...)

function F(2,2) indefinite<z>;  // arity 2, each argument in FF(2)
// creates variables z[0,0], z[0,1], z[1,0], z[1,1] for coefficients
// call as F(x,y), for which polynomial function of x, y, z[i,j] substituted

parameter d;
function G(d,d,d) indefinite<w>;  // arity 3, each argument in FF(d)
// creates variables w[0] .. w[d^3] for coefficients

function G(x,y) : x -> y || y -> x;  // a definite function
G(x,y) : x -> y || y -> x;           // a definite function

(:A x,y,z:  F(x,y) -> z)  // universal quantifier
( :E x :  F(x,x) -> 1 )   // existential quantifier

parameter `in in FF(2);  // keyword as identifier

% solve x;  // runtime query solution-value set for x
% $x;       // same: returns recurrence F(x), domain U={0,1,2}, objective x
% ?x;       // ask if dynamic theorem: exactly one fixed pt x == F(x) with x == 1?
% $ t1 @ (0,1);  // to get sequence (0,1,1,2,3,5,8,...)
% $ t1;          // to get collection { {0}@(0,0) }
% $ x @@  // collection of sequences $x@0, $x@1, $x@2
\end{verbatim}

\clearpage

\section{What the Tortoise Said to Achilles about What the Tortoise
  Said to Achilles}

\textsc{Achilles} had once again overtaken the tortoise, and had
seated himself comfortably on its back.

``So you've understood the Liar?'' said the Tortoise, ``even though he
told you he was lying?  I thought some mathematician or other proved
that the thing couldn't be done?''

``It \emph{can} be done,'' said Achilles, ``It \emph{has} been done.
\emph{Quod erat faciendum.}  You see, it was a simple matter of
solving a system of \emph{simultaneous} equations; and so---''

``But what if the equations did not happen all at the \emph{same}
time?'' the Tortoise interrupted, ``What then?''

``Then I shouldn't have solved them,'' Achilles modestly replied;
``and \emph{you} would have got several times round the world, by this
time!''

``You impress me---\emph{compress}, I mean,'' said the Tortoise; ``for
you are truly dense, and no mistake!  Well now, would you like to hear
of some logical propositions that are supposed to be true or false,
while they are \emph{really} neither of the two?''

``Very much indeed!'' said the Grecian warrior, whereupon he produced
a notebook, a fine quill pen, and a small pot of ink from his
shield-pocket protector.  (He wondered to himself, ``Am \emph{I} not
the protector of this pocket?'').

``Now listen carefully, for I shall reveal to you the Principle of the
Excluded Muddle.''




\qquad \vdots

Achilles rose to his feet and balanced atop the Tortoise, still
writing furiously.  ``So \emph{my} answer depends upon itself, rather
like a turtle standing on its own shell?''

``But just \emph{what} do you suppose \emph{that} answer depends upon,
my gallant Greek?  Surely you can tell me \emph{something} more
particular about it.''

``Very clever, young Tortoise, very clever'' said Achilles, ``but it's
turtles all the way down!  I dare say any \emph{particular}
answer would be nonsense.''

``Nonsense indeed,'' replied the Tortoise, ``But which \emph{kind} of
nonsense?  There are two different methods to this madness, the
\emph{unsteady} kind and the \emph{contingent} kind.  You can tell
them apart if you start from the \emph{beginning}, which is after all
why we call it so.''

``Let us start with the very first turtle---or \emph{ur}tle if you
like, since it has no beginning.''

\qquad \vdots

Achilles sighed, ``With such a long and wandering explanation, which
will surely vex some Logicians of the Twenty-First
Century---\emph{would} you mind renaming yourself \emph{Tor-tu-ous}?''

``As you please!'' replied the rational reptile, admiring the
warrior's tidy notes, ``Provided that \emph{you}, for \emph{your}
excellent penmanship, will call yourself \emph{A Quill Ease}!''

\nocite{carroll-tortoise}

\clearpage

\bibliographystyle{plain}
\nocite{heijenoort}
\begin{raggedright}
  \small
  \bibliography{goedel}

\begin{thebibliography}{10}

\bibitem{agudelo-carnielli}
Juan~C. Agudelo and Walter Carnielli.
\newblock Polynomial ring calculus for modal logics: A new semantics and proof
  method for modalities.
\newblock {\em Review of Symbolic Logic}, 4:150--170, 2011.

\bibitem{arnon}
Dennis~S. Arnon.
\newblock A bibliography of quantifier elimination for real closed fields.
\newblock {\em Journal of Symbolic Computation}, 5:267--274, 1988.

\bibitem{bell}
E.~T. Bell.
\newblock {\em The Development of Mathematics}.
\newblock McGraw-Hill, second edition, 1945.

\bibitem{boole-mal}
George Boole.
\newblock {\em The Mathematical Analysis of Logic, Being an Essay Towards a
  Calculus of Deductive Reasoning}.
\newblock Macmillan, London, 1847.

\bibitem{boole}
George Boole.
\newblock {\em An Investigation of the Laws of Thought, on Which Are Founded
  the Mathematical Theories of Logic and Probabilities}.
\newblock Macmillan, London, 1854.

\bibitem{buchberger}
Bruno Buchberger.
\newblock Gr{\"o}bner bases: {A} short introduction for systems theorists.
\newblock {\em Lecture Notes in Computer Science}, 2178:1--19, 2001.

\bibitem{carnielli-polynomizing}
Walter Carnielli.
\newblock {\em Polynomizing: Logical Inference in Polynomial Format and the
  Legacy of {B}oole}, volume~64 of {\em Studies in Computational Intelligence},
  pages 349--364.
\newblock Springer, 2007.

\bibitem{carroll-barbers}
Lewis Carroll.
\newblock A logical paradox.
\newblock {\em Mind, \em{New Series}}, 3(11):436--438, 1894.

\bibitem{carroll-tortoise}
Lewis Carroll.
\newblock What the tortoise said to {A}chilles.
\newblock {\em Mind}, New Series 4(14):278--280, 1895.

\bibitem{clark}
Allan Clark.
\newblock {\em Elements of Abstract Algebra}.
\newblock Wadsworth, Belmont, CA, 1971.
\newblock Reprinted by Dover in 1984.

\bibitem{ideals-varieties}
David Cox, John Little, and Donal O'{S}hea.
\newblock {\em Ideals, Varieties, and Algorithms: {A}n Introduction to
  Computational Algebraic Geometry and Commutative Algebra}.
\newblock Springer, New York, third edition, 2007.

\bibitem{eves}
Howard Eves.
\newblock {\em Foundations and Fundamental Concepts of Mathematics}.
\newblock PWS-Kent, Boston, third edition, 1990.
\newblock Reprinted by Dover in 1997.

\bibitem{frege}
Gottlob Frege.
\newblock {\em Begriffsschrift, eine der arithmetischen nachgebildete
  {F}ormelsprache des reinen {D}enkens}.
\newblock Louis Nebert, Halle (Saale), 1879.
\newblock English translation ``\emph{Begriffsschrift}, a formula language,
  modeled upon that of arithmetic, for pure thought'' appears in
  \cite{heijenoort}.

\bibitem{galor}
Oded Galor.
\newblock {\em Discrete Dynamical Systems}.
\newblock Springer, New York, 2010.

\bibitem{goedel}
Kurt G{\"o}del.
\newblock {\"U}ber formal unentscheidbare {S}{\"a}tze der {P}rincipia
  mathematica und verwandter {S}ysteme {I}.
\newblock {\em Monatshefte f{\"u}r Mathematik und Physik}, 38:173--198, 1931.
\newblock English translation ``On formally undecidable propositions of
  \emph{{P}rincipia {m}athematica} and related systems {I}'' appears in
  \cite{heijenoort}.

\bibitem{goldstein}
Rebecca Goldstein.
\newblock {\em Incompleteness: {T}he Proof and Paradox of Kurt G{\"o}del}.
\newblock Norton, New York, 2005.

\bibitem{concrete}
Ronald~L. Graham, Donald~E. Knuth, and Oren Patashnik.
\newblock {\em Concrete Mathematics}.
\newblock Addison-Wesley, second edition, 1994.

\bibitem{holmgren}
Richard~A. Holmgren.
\newblock {\em A First Course in Discrete Dynamical Systems}.
\newblock Springer, New York, second edition, 1996.

\bibitem{kapur-narendran}
Deepak Kapur and Paliath Narendran.
\newblock An equational approach to theorem proving in first-order predicate
  calculus.
\newblock {\em {ACM} {SIGSOFT} Software Engineering Notes}, 10(4):63--66, 1985.

\bibitem{nagel}
Ernest Nagel and James~R. Newman.
\newblock {\em G{\"o}del's Proof}.
\newblock New York University Press, revised edition, 2001.

\bibitem{post}
Emil~Leon Post.
\newblock Introduction to a general theory of elementary propositions.
\newblock {\em American Journal of Mathematics}, 43:163--185, 1921.
\newblock Reprinted in \cite{heijenoort}.

\bibitem{roanes-lozano}
Eugenio Roanes-Lozano, Luis~M. Laita, and Eugenio Roanes-Mac{\'\i}as.
\newblock A polynomial model for multi-valued logics with a touch of algebraic
  geometry and computer algebra.
\newblock {\em Mathematics and Computers in Simulation}, 45:83--99, 1998.

\bibitem{russell-letter}
Bertrand Russell.
\newblock Letter to {F}rege, 1902.
\newblock First published in \cite{heijenoort}.

\bibitem{seuren}
Pieter A.~M. Seuren.
\newblock Eubulides as a 20th-century semanticist.
\newblock {\em Language Sciences}, 27:75--95, 2005.

\bibitem{smullyan-name}
Raymond Smullyan.
\newblock {\em What is the Name of This Book?}
\newblock Prentice-Hall, Englewood Cliffs, NJ, 1978.

\bibitem{snapper}
Ernst Snapper.
\newblock The three crises in mathematics: {L}ogicism, intuitionism, and
  formalism.
\newblock {\em Mathematics Magazine}, 4:207--216, 1979.

\bibitem{lisp}
Guy~L. Steele.
\newblock {\em Common Lisp: {T}he Language}.
\newblock Digital Press, second edition, 1990.

\bibitem{stone}
M.~H. Stone.
\newblock The theory of representation for {B}oolean algebras.
\newblock {\em Transactions of the American Mathematical Society}, 40:37--111,
  1936.

\bibitem{tarski-algebra}
Alfred Tarski.
\newblock A decision method for elementary algebra and geometry.
\newblock Technical Report R-109, {RAND} Corporation, Santa Monica, CA, 1948.

\bibitem{heijenoort}
Jean van Heijenoort, editor.
\newblock {\em From {F}rege to {G}{\"o}del: {A} Source Book in Mathematical
  Logic, 1879--1931}.
\newblock Harvard University Press, 1967.

\bibitem{pm}
Alfred~North Whitehead and Bertrand Russell.
\newblock {\em Principia Mathematica}, volume 1--3.
\newblock Cambridge University Press, 1910--1913.

\bibitem{wittgenstein}
Ludwig Wittgenstein.
\newblock {\em Tractatus Logico-Philosophicus}.
\newblock Routledge \& Kegan Paul Ltd, London, 1922.
\newblock Translated from German by C. K. Ogden; reprinted by Dover in 1999.

\bibitem{mathematica}
{Wolfram Research, Inc.}
\newblock {\em Mathematica}.
\newblock Champaign, IL, 2010.
\newblock Version~8.0.

\bibitem{zermelo}
Ernst Zermelo.
\newblock Untersuchungen {\"u}ber die {G}rundlagen der {M}engenlehre {I}.
\newblock {\em Mathe\-ma\-tische An\-na\-len}, 65:261--281, 1908.
\newblock English translation ``Investigations in the foundations of set theory
  {I}'' appears in \cite{heijenoort}.

\end{thebibliography}
\end{raggedright}

\end{document}